\documentclass[10pt, reqno]{amsart}

\usepackage[margin=2.5cm]{geometry}
\usepackage[dvipsnames]{xcolor}
\usepackage[colorlinks=true, linkcolor=Blue, citecolor=Blue]{hyperref}
\usepackage{graphicx}
\usepackage{framed}
\usepackage{amsmath}
\usepackage{amssymb}
\usepackage{amsfonts}
\usepackage{mathrsfs}
\usepackage{array}
\usepackage{enumerate}

\usepackage{mathtools}
\usepackage{tikz-cd}
\usepackage{multirow}
\usepackage{float}
\usepackage{quiver}
\usepackage{subfiles}
\usepackage{subcaption} 
\usepackage{pbox,tkz-graph,xifthen}
\usetikzlibrary{positioning, calc}
\usepackage{url}
\usepackage[style=numeric, maxbibnames=5]{biblatex}
\usepackage{amsthm}

\newcommand{\ZZ}{\mathbb{Z}}

\newcommand{\CC}{\mathbb{C}}
\newcommand{\C}{\mathcal{C}}
\newcommand{\OO}{\mathcal{O}}
\newcommand{\I}{\mathcal{I}}

\newcommand{\D}{\mathcal{D}}

\newcommand{\Gal}{\operatorname{Gal}}
\newcommand{\Aut}{\operatorname{Aut}}

\newcommand{\weight}{\operatorname{weight}}
\newcommand{\map}{\rightarrow}
\newcommand{\iso}{\cong}
\newcommand{\A}{\mathcal{A}}
\newcommand{\B}{\mathcal{B}}

\theoremstyle{plain} 
\newtheorem{lemma}{Lemma}[section]

\newtheorem{prop}[lemma]{Proposition}
\newtheorem{cor}[lemma]{Corollary}

\newtheorem{claim*}{Claim}
\newtheorem{thm}[lemma]{Theorem}

\theoremstyle{definition}
\newtheorem{defn}[lemma]{Definition}
\newtheorem{fact}[lemma]{Fact}

\newtheorem{remark}[lemma]{Remark}
\newtheorem{example}[lemma]{Example}

\definecolor{purple}{HTML}{692f56} 
\def\edge#1#2#3{\frac{#1{-}#2}{#3}}

\newcommand{\dx}{0.8}       
\newcommand{\dy}{1.25}      
\newcommand{\ov}{0.18}      
\newcommand{\vup}{1.8}     
\newcommand{\vdown}{-1.3}

\newcommand{\seg}[2]{\draw ($(#1)!-\ov!(#2)$) -- ($(#2)!-\ov!(#1)$);}

\newcommand{\stub}[2]{\draw ($(#1)!-\ov!(#2)$) -- ($(#1)!0.30!(#2)$);}

\newcommand{\zigzag}{
  
  \coordinate (P0) at (0,0);                       
  \coordinate (A1) at ({\dx},{\dy});               
  \coordinate (V1) at ({2*\dx},0);                 
  \coordinate (A2) at ({3*\dx},{\dy});             
  \coordinate (Sa) at ({4*\dx},0);                 
  \coordinate (Ta) at ({5*\dx},{\dy});             
  \coordinate (V2) at ({6*\dx},0);                 
  \coordinate (A3) at ({7*\dx},{\dy});             
  \coordinate (P1) at ({8*\dx},0);                 
  
  \draw[very thick] (0,\vup) -- (0,\vdown);        
  \seg{P0}{A1}                                     
  \seg{A1}{V1}                                     
  \seg{V1}{A2}                                     
  \stub{A2}{Sa}                                    
  \stub{V2}{Ta}                                    
  \seg{V2}{A3}                                     
  \node at ({4.5*\dx},{0.5*\dy}) {$\cdots$};
  
  \node[left]                  at (0,{1.20*\dy})        {$v_0$};
  \node[blue,left]             at (0,{-0.70*\dy})       {$d_0$};
  \node[blue,above left=-1pt]  at ($(P0)!0.78!(A1)$)    {$d_1$};
  \node[below right=-1pt]      at ($(P0)!0.50!(A1)$)    {$v_1$};
  \node[right]                 at ($(A1)!0.45!(V1)$)    {$v_2$};
  \node[right]                 at ($(V1)!0.55!(A2)$)    {$v_3$};
  \node[left]                  at ($(V2)!0.60!(A3)$)    {$v_{r-1}$};
}

\begin{document}
\title{Reduction types of curves in tame extensions}
\author{Júlia Martínez-Marín}
\address{School of Mathematics, University of Bristol, Bristol, BS8 1UG, UK}
\email{julia.martinezmarin@bristol.ac.uk}

\begin{abstract}
Let $C/K$ be a curve of genus $g\geq 1$ over a discretely valued field $K$ with algebraically closed residue field. In this paper, we study the behaviour of the reduction type of $C$ under finite tame extensions. If $L/K$ is a tame extension of degree $n$, we give sufficient conditions on the reduction type $R$ of $C$ so that the reduction type of the base change $C_L$ is uniquely determined by $R$ and $n$. 
In particular, if $C$ has abelian and toric ranks zero, then $R$ satisfies these conditions.
We also give an algorithm to compute the possibilities for the reduction type of $C_L$ given only $R$ and $n$. The main method involves studying harmonic covers of graphs with a cyclic action.

As a result, we obtain a characterisation of tame curves with potentially good reduction in terms of their reduction type, and a new proof of a result by Liu-Lorenzini-Raynaud regarding the non-realisability of certain reduction types. 
\end{abstract}

\maketitle
\tableofcontents
\section{Introduction}
Let $K$ be a discretely valued field with complete valuation ring $O_K$ and algebraically closed residue field $k$ of characteristic $p\geq 0$. Let $C/K$ be a smooth, proper, geometrically connected curve of genus $g(C)\geq 1$. Let $\C/O_K$ be the minimal regular normal crossings model of $C$. The reduction type of $C/K$ encodes the discrete invariants of the special fibre of $\C$: multiplicities of the irreducible components, their geometric genera, and the intersection number for every pair of components. For a given curve $C/K$, its reduction type depends on the ground field $K$. 
In this paper we study the behaviour of the reduction type of a curve under tame base change. 

When $p\nmid n$, let $C(n)$ denote the base change of $C$ by the unique degree $n$ extension $K(n)$ of $K$ in a separable closure $K^s$. If $C$ is an elliptic curve, then the reduction type of $C(n)$ is completely determined by the reduction type of $C$ and $n$, which can be proved via Tate's algorithm \cite[Theorem 3]{Tim-Vladimir2013}. If $g(C)=2$, this is no longer true. In \cite{Namikawa-Ueno}, Namikawa and Ueno give examples of curves $C_1,C_2/\CC((t))$ with the same reduction type, such that their base changes $C_1(2),C_2(2)$ have different reduction types. Thus, it is natural to ask under what conditions is the reduction type of $C(n)$ uniquely determined by the reduction type of $C$ and $n$. This is answered in the following corollary, which holds for curves of genus $g\geq 1$. It depends on a simple invariant, which we call the \emph{weight}. Informally, the weight of a set of components in a reduction type is the gcd of the multiplicities of all the components intersecting those in the set (see Section \ref{sec: weight}).

\begin{cor}[{Corollary \ref{cor: reduction type}}] \label{unique reduction type}
    Let $p\nmid n$, and let $R$ be a reduction type. Suppose that each prime $q\mid n$ satisfies one of the following: 
   \begin{enumerate}[(a)]
        \item \label{a} $q\nmid \operatorname{weight}(R)$, and any component $v$ with geometric genus $g(v)\geq1$ and any cycle $Z$ in $R$ satisfy
        \[
        q\nmid \operatorname{weight}(v)  \text{ and }   q\nmid \operatorname{weight}(Z).\]
        
        \item \label{b} \(q\mid \operatorname{weight}(R)\), and $R$ has no cycles and exactly one component with geometric genus $g\geq 1$.
        \item \label{c} \( q\mid \operatorname{weight}(R)\), and $R$ has exactly one cycle and no components with geometric genus $g\geq 1$.
    \end{enumerate}
    Then for any curve $C$ with reduction type $R$, the reduction type of $C(n)$ is uniquely determined, and it depends only on $R$ and $n$, and not on $C$ or $K$. In particular, the result holds for curves with abelian and toric ranks zero. 
\end{cor}

From Namikawa and Ueno's examples in \cite{Namikawa-Ueno}, we can check that the conditions in the corollary are not only sufficient but also necessary in genus 2. We expect this to be the case in any genus $g\geq 1$, but there are some realisability issues that prevent us from proving this (see Remark \ref{realisability}). Thus, we first look at reduction types as combinatorial objects (not necessarily arising from curves) and their \emph{reduction type covers}, for which the conditions turn out to be both sufficient and necessary. 

\begin{thm} [{Theorem \ref{thm: reduction type cover}}]\label{unique reduction type cover}
    Let $n\geq 1$. A reduction type $R$ has a unique reduction type cover of degree $n$ if and only if for each prime $q\mid n$, one of \eqref{a}, \eqref{b}, or \eqref{c} of Corollary \ref{unique reduction type} holds.

\end{thm}

In genus 2, if $n$ is odd, then all 104 families of reduction types satisfy the conditions in Corollary \ref{unique reduction type}. If $n$ is even, then exactly 102 families satisfy the conditions. The 2 remaining families, each have two possibilities after a degree $n$ extension. We give one example for each these families below (multiplicities of irreducible components are marked \textcolor{blue}{1}, \textcolor{blue}{2}, ..., and their geometric genera are marked \textcolor{blue}{g1}, \textcolor{blue}{g2}, ..., and omitted if zero).
\vspace{-2pt}
\begin{figure}[htbp]
    \centering
    \begin{subfigure}[b]{0.45\textwidth}
        \centering
        \begin{tikzpicture}[xscale=1,yscale=0.9]
        \draw[thick,shorten >=2] (-1/3,0)--(107/30,0) node[inner sep=2pt,blue,above 
        left,scale=0.8] {2} node[white, inner sep=2pt,below left,scale=0.8] {$v$};
        \draw[shorten <=-3pt] (0,0)--node[inner sep=2pt,scale=0.8,blue,right] {1} 
        (0,3/5);
        \draw[shorten <=-3pt] (4/5,0)--node[inner sep=2pt,scale=0.8,blue,right] {1} 
        (4/5,3/5);
        \draw[shorten <=-3pt, shorten >=-3pt] (29/15,0)--node[inner 
        sep=2pt,scale=0.8,blue,right] {2} (29/15,4/5);
        \draw[thick,shorten >=2] (8/5,4/5)--(113/30,4/5) node[inner sep=2pt,blue,above 
        left,scale=0.8] {2 \smash{g}1} node[inner sep=2pt,below left,scale=0.8] 
        {$v$};
        \end{tikzpicture}
        \caption*{(A)} 
        \label{fig:diagram1}
    \end{subfigure}
    \hspace{0.3cm}
    \begin{subfigure}[b]{0.45\textwidth}
        \centering
        \begin{tikzpicture}[xscale=1,yscale=0.9]
        \draw[thick,shorten >=2] (-1/3,0)--(113/30,0) node[inner sep=2pt,blue,above 
        left,scale=0.8] {2} node[inner sep=2pt,below left,scale=0.8] {$v$};
        \draw[shorten <=-3pt, shorten >=-3pt] (0,0)--node[inner 
        sep=2pt,scale=0.8,blue,left] {2} (0,3/5);
        \draw[shorten <=-3pt, shorten >=-3pt] (0,3/5)--node[inner 
        sep=2pt,scale=0.8,blue,above] {2} (3/5,3/5);
        \draw[shorten <=-3pt, shorten >=-3pt] (3/5,0)--node[inner 
        sep=2pt,scale=0.8,blue,right] {2} (3/5,3/5);
        \draw[shorten <=-3pt, shorten >=-3pt] (29/15,0)--node[inner 
        sep=2pt,scale=0.8,blue,right] {2} (29/15,4/5);
        \draw (4/3,4/5)--node[inner sep=2pt,scale=0.8,blue,above,pos=0.6] {2} 
        (41/15,4/5);
        \draw[shorten <=-3pt] (23/15,4/5)--node[inner sep=2pt,scale=0.8,blue,left] {1} 
        (23/15,7/5);
        \draw[shorten <=-3pt] (38/15,4/5)--node[inner sep=2pt,scale=0.8,blue,right] {1} 
        (38/15,7/5);
        \end{tikzpicture}
        \caption*{(B)} 
        \label{fig:diagram2}
    \end{subfigure}
\end{figure}
\vspace{-4pt}

\noindent Then the reduction types A1 and A2 (resp. B1 and B2) below arise as reduction type covers of A (resp. B) of degree 2. They are semistable, and so they stay different for every even degree $n$ as well. These covers are realised by base changing the curves listed in \cite[pages 159,172,181,184]{Namikawa-Ueno} under the labels $2\text{I}_{0-1}$ and $\text{II}^*_{2-0}$, ($2\text{I}_4-2$ and $\text{II}^*_{2-4}$), respectively.
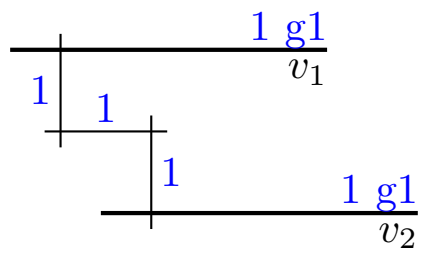
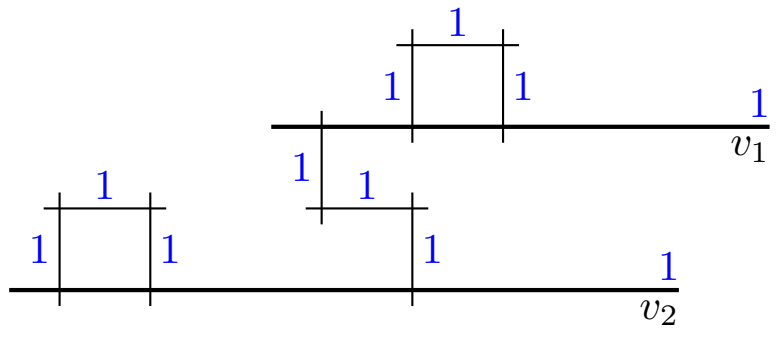
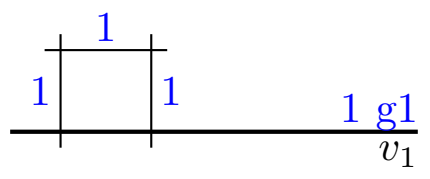
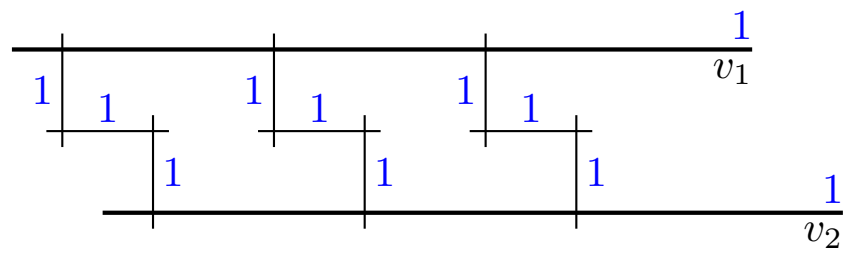
\vspace{-13pt}
\begin{figure}[htbp]
    \centering
    \begin{subfigure}[b]{0.45\textwidth}
        \centering
        \begin{tikzpicture}[xscale=1,yscale=0.9]
        \draw[thick,shorten >=2] (4/15,0)--(73/30,0) node[inner sep=2pt,blue,above 
        left,scale=0.8] {1 \smash{g}1} node[inner sep=2pt,below left,scale=0.8] 
        {$v_2$};
        \draw[shorten <=-3pt, shorten >=-3pt] (3/5,0)--node[inner 
        sep=2pt,scale=0.8,blue,right] {1} (3/5,3/5);
        \draw[shorten <=-3pt, shorten >=-3pt] (3/5,3/5)--node[inner 
        sep=2pt,scale=0.8,blue,above] {1} (0,3/5);
        \draw[shorten <=-3pt, shorten >=-3pt] (0,3/5)--node[inner 
        sep=2pt,scale=0.8,blue,left] {1} (0,6/5);
        \draw[thick,shorten >=2] (-1/3,6/5)--(11/6,6/5) node[inner sep=2pt,blue,above 
        left,scale=0.8] {1 \smash{g}1} node[inner sep=2pt,below left,scale=0.8] 
        {$v_1$};
        \end{tikzpicture}
        \caption*{(A1)}
        \label{fig:A1}
    \end{subfigure}
    \hspace{0.3cm}
    \begin{subfigure}[b]{0.45\textwidth}
        \centering
        \begin{tikzpicture}[xscale=1,yscale=0.9]
        \draw[thick,shorten >=2] (-1/3,0)--(25/6,0) node[inner sep=2pt,blue,above 
        left,scale=0.8] {1} node[inner sep=2pt,below left,scale=0.8] {$v_2$};
        \draw[shorten <=-3pt, shorten >=-3pt] (0,0)--node[inner 
        sep=2pt,scale=0.8,blue,left] {1} (0,3/5);
        \draw[shorten <=-3pt, shorten >=-3pt] (0,3/5)--node[inner 
        sep=2pt,scale=0.8,blue,above] {1} (3/5,3/5);
        \draw[shorten <=-3pt, shorten >=-3pt] (3/5,0)--node[inner 
        sep=2pt,scale=0.8,blue,right] {1} (3/5,3/5);
        \draw[shorten <=-3pt, shorten >=-3pt] (7/3,0)--node[inner 
        sep=2pt,scale=0.8,blue,right] {1} (7/3,3/5);
        \draw[shorten <=-3pt, shorten >=-3pt] (7/3,3/5)--node[inner 
        sep=2pt,scale=0.8,blue,above] {1} (26/15,3/5);
        \draw[shorten <=-3pt, shorten >=-3pt] (26/15,3/5)--node[inner 
        sep=2pt,scale=0.8,blue,left] {1} (26/15,6/5);
        \draw[thick,shorten >=2] (7/5,6/5)--(143/30,6/5) node[inner sep=2pt,blue,above 
        left,scale=0.8] {1} node[inner sep=2pt,below left,scale=0.8] {$v_1$};
        \draw[shorten <=-3pt, shorten >=-3pt] (7/3,6/5)--node[inner 
        sep=2pt,scale=0.8,blue,left] {1} (7/3,9/5);
        \draw[shorten <=-3pt, shorten >=-3pt] (7/3,9/5)--node[inner 
        sep=2pt,scale=0.8,blue,above] {1} (44/15,9/5);
        \draw[shorten <=-3pt, shorten >=-3pt] (44/15,6/5)--node[inner 
        sep=2pt,scale=0.8,blue,right] {1} (44/15,9/5);
        \end{tikzpicture}
        \caption*{(B1)}
        \label{fig:B1}
    \end{subfigure}
\end{figure}
    \vspace{-0.7cm}
\begin{figure}[htbp]
        \centering
        \begin{subfigure}[b]{0.45\textwidth}
        \centering
        \begin{tikzpicture}[xscale=1,yscale=0.9]
        \draw[thick,shorten >=2] (-1/3,0)--(73/30,0) node[inner sep=2pt,blue,above 
        left,scale=0.8] {1 \smash{g}1} node[inner sep=2pt,below left,scale=0.8] 
        {$v_1$};
        \draw[shorten <=-3pt, shorten >=-3pt] (0,0)--node[inner 
        sep=2pt,scale=0.8,blue,left] {1} (0,3/5);
        \draw[shorten <=-3pt, shorten >=-3pt] (0,3/5)--node[inner 
        sep=2pt,scale=0.8,blue,above] {1} (3/5,3/5);
        \draw[shorten <=-3pt, shorten >=-3pt] (3/5,0)--node[inner 
        sep=2pt,scale=0.8,blue,right] {1} (3/5,3/5);
        \end{tikzpicture}
        \caption*{(A2)}
        \label{fig:A2}
        \end{subfigure}
        \hspace{0.3cm}
        \begin{subfigure}[b]{0.45\textwidth}
        \centering
        \begin{tikzpicture}[xscale=1,yscale=0.9]
        \draw[thick,shorten >=2] (4/15,0)--(157/30,0) node[inner sep=2pt,blue,above 
        left,scale=0.8] {1} node[inner sep=2pt,below left,scale=0.8] {$v_2$};
        \draw[shorten <=-3pt, shorten >=-3pt] (3/5,0)--node[inner 
        sep=2pt,scale=0.8,blue,right] {1} (3/5,3/5);
        \draw[shorten <=-3pt, shorten >=-3pt] (3/5,3/5)--node[inner 
        sep=2pt,scale=0.8,blue,above] {1} (0,3/5);
        \draw[shorten <=-3pt, shorten >=-3pt] (0,3/5)--node[inner 
        sep=2pt,scale=0.8,blue,left] {1} (0,6/5);
        \draw[shorten <=-3pt, shorten >=-3pt] (2,0)--node[inner 
        sep=2pt,scale=0.8,blue,right] {1} (2,3/5);
        \draw[shorten <=-3pt, shorten >=-3pt] (2,3/5)--node[inner 
        sep=2pt,scale=0.8,blue,above] {1} (7/5,3/5);
        \draw[shorten <=-3pt, shorten >=-3pt] (7/5,3/5)--node[inner 
        sep=2pt,scale=0.8,blue,left] {1} (7/5,6/5);
        \draw[shorten <=-3pt, shorten >=-3pt] (17/5,0)--node[inner 
        sep=2pt,scale=0.8,blue,right] {1} (17/5,3/5);
        \draw[shorten <=-3pt, shorten >=-3pt] (17/5,3/5)--node[inner 
        sep=2pt,scale=0.8,blue,above] {1} (14/5,3/5);
        \draw[shorten <=-3pt, shorten >=-3pt] (14/5,3/5)--node[inner 
        sep=2pt,scale=0.8,blue,left] {1} (14/5,6/5);
        \draw[thick,shorten >=2] (-1/3,6/5)--(139/30,6/5) node[inner sep=2pt,blue,above 
        left,scale=0.8] {1} node[inner sep=2pt,below left,scale=0.8] {$v_1$};
        \end{tikzpicture}
        \caption*{(B2)}
        \label{fig:B2}
    \end{subfigure}
\end{figure}
\vspace{-5pt}

The proof of Theorem \ref{unique reduction type cover} is constructive, and we give an algorithm in Section \ref{algorithm} to compute the set of reduction type covers of degree $n$ of a given reduction type $R$. We also give a formula for the abelian and toric ranks of $C(n)$ (and its Jacobian), in terms only of the reduction type of $C$ (Proposition \ref{abelian and toric rank}).

Our method leads to a new proof of the following realisability result by Liu, Lorenzini and Raynaud \cite[Corollary 7.4]{Liu-Lorenzini-Raynaud}.

\begin{cor} [{Corollary \ref{cor: Liu-Lorenzini-Raynaud}}]\label{Liu-Lorenzini-Raynaud}
Suppose a curve $C/K$ has abelian and toric ranks zero. Let $d_1,\dots,d_r$ be the multiplicities of the components in the reduction type $R$ of $C$. If $d:=\gcd(d_1,\dots,d_r) =\weight(R)> 1$, then $p>0$ and $d$ is a power of $p$.
\end{cor}

We call a curve \emph{tame} if it attains semistable reduction over a tame extension. In \cite[Theorem 2.1]{Lorenzini90}, Lorenzini gives a necessary condition for a tame curve to have potentially good reduction. We extend this result and give a characterisation of tame curves with potentially good reduction in terms of their reduction type. 

\begin{thm}[{Theorem \ref{thm: tame curve}}]\label{tame curve}
    A tame curve of genus $g\geq 1$ has potentially good reduction if and only if its reduction type has exactly one principal component and no cycles. 
\end{thm}

\begin{remark}
    One can still apply all of the results above with slightly different assumptions on $K$. If the valuation ring $O_K$ is not assumed to be complete, but it is excellent (e.g. if $\text{char } K=0)$, and the residue field $k$ is not necessarily algebraically closed, but it is perfect, then we may argue as follows. By \cite[Proposition 9.3.28]{Liu}, the reduction type of $C/K$ is preserved under base change by an unramified extension or by the completion of $K$. Thus, we may base change to the maximal unramified extension $K^{\text{nr}}$ of $K$ and assume it is complete. Then the residue field of $K^{\text{nr}}$ is algebraically closed and $C/K$ and $C_{K^{\text{nr}}}/K^{\text{nr}}$ have the same reduction type, so we can apply our results to study the behaviour of the reduction type of $C/K$ under base change by a tame extension of ramification degree $n$ (no longer unique).
\end{remark}

\begin{remark} \label{realisability}
As explained above, the converse of Corollary \ref{unique reduction type} holds for curves of genus 2. However, generalising this to higher genus is complicated. There are three difficulties:
\begin{enumerate}
    \item Realisability of reduction types: in \cite[Corollary 3.7]{Winters}, Winters showed that if $\text{char }k=0$ and $\weight(R)=1$ for a reduction type $R$ (as defined in Section \ref{sec: reduction type}), then there exists a smooth, projective, geometrically integral curve $C/K$ with reduction type $R$. However, in mixed or positive characteristic and genus $g>1$, the realisability of reduction types remains open. 
    \item Realisability of reduction type covers: 
    even in the cases covered by Winters' theorem we still do not know that every possible reduction type cover is realisable by a curve in a tame extension. In genus 2, this is possible by Namikawa and Ueno's examples in \cite{Namikawa-Ueno}.
    \item Symmetry: suppose that the conditions \eqref{a}, \eqref{b} or \eqref{c} in Corollary \ref{unique reduction type} do not hold for $R$ and some prime $q\mid n$. Then, Theorem \ref{unique reduction type cover} implies that there is more than one reduction type cover of degree $n$. However, these could be equal as reduction types, after forgetting the ``covering map''. To show that there are at least two different reduction types arising as degree $n$ reduction type covers of $R$, it is enough to show that there are at least two degree $n$ cyclic covers $\Gamma_1,\Gamma_2\to \Gamma$ of the dual graph $\Gamma$ of $R$ with some prescribed dilation, and such that $\Gamma_1\not\iso \Gamma_2$ as graphs. We expect this to hold due to Remark \ref{rem: symmetry}. 
\end{enumerate}
\end{remark}

\subsection{Methodology and relation to the literature}
\mbox{}
\vspace{5pt}

In \cite{Halle-Nicaise}, Halle and Nicaise study the local behaviour of strict normal crossings (snc) models under base change by an extension of degree $n$, with $p\nmid n$. Let $C/K$ be a curve with reduction type $R$, and let $\C/O_K$ be the minimal snc (msnc) model of $C$. Let $\C_n$ be the normalisation of the base change $\C\times_{O_K}O_{K(n)}$, and denote the normalisation map by $\nu$. Let $\C(n)$ denote the minimal desingularisation of $\C_n$, and let $\C(n)_{\operatorname{min}}/O_{K(n)}$ denote the msnc model of $C(n)$. We are interested in obtaining the combinatorial information of the special fibre of $\C(n)_{\operatorname{min}}$ from that of $\C$, since this is equivalent to understanding the reduction type of $C(n)$ given $R$, due to the correspondence between msnc and mrnc reduction types \cite{Jakab} (see Section \ref{sec: reduction type}). Halle and Nicaise show that $\C(n)$ is an snc model, so we have maps
\[\begin{tikzcd}
	{\C(n)} & {\C_n} & {\C\times_{O_K}O_{K(n)}} & \C \\
	{\C(n)_{\operatorname{min}}}
	\arrow["{\text{desing.}}", from=1-1, to=1-2]
	\arrow[from=1-1, to=2-1]
	\arrow["\nu", from=1-2, to=1-3]
	\arrow[from=1-3, to=1-4]
\end{tikzcd}\]

\noindent Since $\C(n)\to \C(n)_{\operatorname{min}}$ is obtained by contracting specific $(-1)$-curves, it is enough to understand the combinatorial information of the special fibre $\C(n)_k$, given that of $\C_k$. Halle and Nicaise give a full description of the induced map $\C(n)_k\to (\C_n)_k$. They also give results on the multiplicities of the components in $(\C_n)_k$, as well as the number of components in $(\C_n)_k$ lying above each non-principal component in $\C_k$. They give a global description of the map $(\C_n)_k\to \C_k$ when $n$ is coprime to $p$ and to 
\[
e(C):=\text{lcm}\{m_v:v \text{ a principal component of } R\},
\]
where $m_v$ denotes the multiplicity of $v$. In this case, for each principal component $v$ in $\C_k$, there is a unique component in $(\C_n)_k$ lying above $v$. From this, they can show that the reduction type of $C(n)$ is determined by the reduction type of $R$ \cite[Proposition 4.2.3.3]{Halle-Nicaise}. 
However, we are often interested in $n$ not coprime to $e(C)$. For example, if $C$ is a tame curve, then the minimal extension $L/K$ over which $C_L$ becomes semistable has degree $e(C)$ \cite{Halle}. 

Using Halle and Nicaise's local results, we obtain results on the number of components in $(\C_n)_k$ lying above any given component in $\C_k$, as well as their genera. From this, we can reduce the problem of understanding $(\C_n)_k\to \C_k$ to studying certain \emph{covers} of the dual graph $\Gamma$ of $R$. In graph theory, these are called harmonic covers. Moreover, our assumptions on $K$ and $n$ imply that $K(n)/K$ is a cyclic extension, and this gives rise to a cyclic action on these harmonic covers. Such covers (with associated \emph{dilation datum} $D$), are classified by a sheaf cohomology group called the dilated cohomology of $(\Gamma,D)$ \cite{Len-Ulrisch-Zakharov}. Thus, our problem reduces to understanding the possible dilation data $D$, and studying the dilated cohomology of $(\Gamma,D)$. 

In \cite{Lorenzini90}, Lorenzini studies a closely related problem in the ``opposite direction''. He assumes that $C$ is tame and its Jacobian has potentially good reduction. Writing $n=e(C)$, this implies that the dual graph of the reduction type of $C(n)$ is a tree. Under these assumptions, he obtains bounds on the number of vertices with at least 3 adjacent vertices in the dual graph of the msnc model of $C$. 

Lastly, we note that there is a notion of \emph{branched cover} of a graph with an action, studied e.g. in \cite{Gambheera-Vallieres}, \cite{Kataoka}, \cite{Mizuno}. Such covers are special examples of harmonic covers with an action that have trivial dilation along edges. Since in our case, edges might have non-trivial dilation (by Fact \ref{number of points above intersection point}), we work with the more general notion of harmonic covers.

\subsection{Outline}
\mbox{}
\vspace{5pt}

We begin in Section \ref{sec: preliminaries} with preliminaries on models of curves, reduction types, dual graphs and harmonic covers of graphs. We also review results on tame base change of models \cite{Halle-Nicaise} and dilated cohomology, which classifies harmonic covers of graphs with an action \cite{Len-Ulrisch-Zakharov}. In Section \ref{sec: genus formula}, we obtain a formula for the genus of a component in a reduction type after tame base change, and look at some consequences. In Section \ref{sec: main results}, we show the main results. In Sections \ref{sec: results on dilated cohomology} and \ref{sec: results on cyclic covers} we give a formula for the dilated cohomology of a graph, and use this to prove results relating the dilated cohomology to cyclic covers of graphs with dilation. In Sections \ref{sec: defn of reduction type covers} and \ref{sec: results on reduction type covers} we define reduction type covers and prove Theorem \ref{unique reduction type cover} and Corollary \ref{unique reduction type}. In Section \ref{sec: other consequences}, we use similar methods to show Corollary \ref{cor: Liu-Lorenzini-Raynaud} and Theorem \ref{thm: tame curve}. In Section \ref{sec: computing}, we give an algorithm that takes as input a reduction type $R$ and a positive integer $n$ and outputs the set of reduction type covers of degree $n$ of $R$. In Section \ref{sec: abelian and toric ranks}, we obtain a formula for the abelian and toric ranks of $C(n)$, given only the reduction type of $C$. We end in Section \ref{sec: examples} with examples of reduction type covers.

\subsection*{Acknowledgements}
I would like to thank my advisor Tim Dokchitser for suggesting this problem and for many helpful discussions throughout the development of the paper. I am also grateful to Qing Liu and Yoav Len for useful conversations. Lastly, I would like to thank Gergely Jakovác for conversations regarding Section \ref{sec: Galois theory of valuations}. I am supported by the Heilbronn Institute for Mathematical Research.

\subsection*{AI statement}
Free versions of ChatGPT and Claude were used for literature search and drawing figures. The ideas in the paper and the writing are due to the author. 

\subsection*{Notation}
Throughout, $C/K$ denotes a smooth, proper, geometrically connected curve with genus $g(C)\geq 1$. For $p\nmid n$, we write $K(n)$ for the unique degree $n$ extension of $K$ in a separable closure $K^s$ and denote by $C(n)$ the base change $C\times_K K(n)$. If $\C/O_K$ is a normal model of $C$, we denote by $\C_k$ its special fibre. We write $\C_n$ for the normalisation of $\C\times_{O_K} O_{K(n)}$, and $\C(n)$ for the minimal desingularisation of $\C_n$. We call an irreducible component $v$ in $\C_k$ simply a \emph{component}, and we denote by $p(v)$ and $g(v)$ its arithmetic and geometric genus, respectively. We say that a component $u$ in $(\C_n)_k$ or $\C(n)_k$ \emph{lies above} a component $v$ in $\C_k$ if $u$ maps surjectively to $v$ under the corresponding maps $\C(n)\to \C_n \to \C$. We extend this notion to unions of components such as chains.

Throughout, $\Gamma$ denotes a graph, and $X(\Gamma)$ denotes the union of the \emph{vertices} $V(\Gamma)$ and \emph{edges} $E(\Gamma)$. We often identify components in a reduction type or a special fibre with the corresponding vertices in the dual graph, and denote them by $v$. Unless otherwise stated, $\A$ denotes a cyclic group of order $n$.

To ease notation, we sometimes write $(a,b)$ and $(a_1,\dots,a_r)$ for $\gcd(a,b)$ and $\gcd(a_1,\dots,a_r)$, respectively. 

\subsection*{Convention} Throughout the paper, we will not consider the reduction type $[m]I_r$, i.e. a single looping chain of $r$ components of genus 0 (or a node if $r=1$) of multiplicity $m$. This is the only reduction type with no principal components, so ignoring it simplifies certain definitions and proofs. The only other Kodaira type that can become $I_{r}$ after a tame extension is $I_{r'}^*$, so whenever a result is about base change, we also ignore the reduction type $[m]I_r^*$. We note that the main results in this paper are already well understood for these reduction types. For example, if an elliptic curve $E/K$ has reduction type $I_r$ (resp. $I_r^*$), then its base change $E(n)$ has reduction type $I_{rn}$ (resp. $I_{rn}^*$ if $n$ is odd, or $I_{rn}$ if $n$ is even) \cite{Tim-Vladimir2013}. 

\section{Preliminaries} \label{sec: preliminaries}
\subsection{Regular models of curves}
\mbox{}
\vspace{5pt}

A \emph{model} $\C/O_K$ of $C/K$ is a proper, flat scheme with generic fibre $\C_K\iso C$ and connected, $1$-dimensional special fibre $\C_k$, not necessarily reduced or irreducible. If $\C$ is a regular scheme, we call it a \emph{regular model}. We call it a \emph{regular model with normal crossings (rnc)} if moreover the special fibre $\C_k\subset \C$ is a divisor with normal crossings. We call $\C$ a \emph{strict normal crossings model (snc)} if it is an rnc model such that $\C_k$ has no self-intersecting components. Every curve has a minimal regular model, dominated by any other regular model, a minimal rnc model (\emph{mrnc}), dominated by any other rnc model, and a minimal snc model (\emph{msnc}), dominated by any other snc model. These are unique up to isomorphism. See \cite[Section 3]{Saito} for an overview and references.

The special fibre of any model $\C/O_K$ of $C/K$ gives a divisor $[\C_k]=\sum_v m_vv$ on $\C$. The sum is finite, taken over irreducible curves $v$ that we refer to as \emph{components}. Let 
\begin{itemize}
    \item $m_v$ be the \emph{multiplicity} of the component $v$ in $\C_k$,
    \item $g(v)$ be the geometric genus of $v$, and
    \item $p(v)$ be the arithmetic genus of $v$. 
\end{itemize}
Suppose that $\C$ is an rnc model. Then $v$ can only have nodes (transversal self-intersection points) as singularities, and $p(v)-g(v)\geq 0$ is the number of nodes on $v$. 

\begin{defn}\cite[Definition 3.1]{Tim2025} \label{incident components}
    The \emph{incident components} to a component $v$ is a multiset $I(v)$ with
    \begin{itemize}
        \item every component $v'\neq v$ taken $|v\cap v'|$ times,
        \item $v$ itself taken 2($p(v)$-$g(v)$) times (so every node on $v$ contributes $v$ to $I(v)$ twice).
    \end{itemize}
\end{defn}

\subsection{Reduction types and dual graphs} \label{sec: reduction type}
\mbox{}
\vspace{5pt}

\begin{defn} \cite[Definition 3.6]{Tim2025}
    An \emph{rnc reduction type} is a tuple $R=(S;m;p;g; \cdot)$ with
    \begin{itemize}
        \item $S$ a non-empty finite set, whose elements $v$ are called components.
        \item $m=(m_v)$, $p=(p(v))$, $g=(g(v))$ lists of integers indexed by $v\in S$, all with $m_v\geq 1$ (multiplicity) and $p(v)\geq g(v)\geq 0$ (arithmetic/geometric genus).
        \item $\cdot : \ZZ S \times \ZZ S \to \ZZ$ a symmetric bilinear pairing with $v\cdot v'\geq 0$ for $v\neq v'$, and such that the left kernel of $\cdot$ is of rank one and contains $\sum_{v\in S} m_v v$. We say $v$ and $v'$ \emph{intersect} if $v\cdot v'\neq 0$. 
    \end{itemize}
    If the tuple $R$ additionally satisfies
    \begin{itemize}
        \item Every $v\in S$ with $p(v)=0$ and $v\cdot v=-1$ has $|I(v)|\geq 3$, with $I(v)$ as in Definition \ref{incident components}.
    \end{itemize}
    then we call $R$ an \emph{mrnc reduction type}, or simply a \emph{reduction type}.
\end{defn}

\begin{defn}
    An \emph{snc reduction type} is an rnc reduction type $R=(S;m;p;g; \cdot)$ such that $p(v)=g(v)$ for every $v\in S$. A \emph{msnc} reduction type is an snc reduction type that moreover satisfies 
    \begin{itemize}
        \item Every $v\in S$ with $p(v)=0$ and $v\cdot v=-1$ has either $|I(v)|\geq 3$ or $I(v)=\{v',v'\}$ for some $v\neq v'\in S$.
    \end{itemize}
\end{defn}

The special fibre of an rnc (resp. snc) model of a curve $C$ naturally gives an rnc (resp. snc) reduction type. One can also define a reduction type associated to the minimal regular model of $C$, and there is a one-to-one correspondence between these and mrnc reduction types \cite[Theorem 4.10]{Jakab}. Since msnc reduction types of curves are obtained from blowing up mrnc reduction types at finitely many points, there is also a one-to-one correspondence between mrnc reduction types and msnc reduction types \cite[Theorem 4.6]{Jakab}.

\begin{remark} \label{rem: msnc mrnc}
    The msnc and mrnc models of a curve only differ by blow ups at the points of self-intersection of components in the special fibre. For each such blow up at a point of a component $v$, we get a new component $u$ intersecting $v$ transversally twice, and of multiplicity $m_u=2m_v$.
\end{remark}

Given a reduction type, one can associate to it a (total) genus. If $R$ is the reduction type of a curve $C$ of genus $g(C)$, then the total genus of $R$ is $g(C)$ (see \cite{Tim2025} or Theorem \ref{Euler characteristic}). 

\begin{defn}\cite[Definition 4.1]{Tim2025}
A component $v$ in an rnc reduction type is called \emph{principal} if $g(v)\geq 1$ or the number of incident components satisfies $|I(v)|\geq 3$. 
\end{defn}

Non-principal components have $p(v)=g(v)=0$, and meet the rest of the special fibre in 1 or 2 points. They form \emph{chains}. A chain $v_1,\dots,v_r$ is \emph{outer} if it meets one principal component $v_0$, and \emph{inner} if it meets two principal components $v_0$ and $v_{r+1}$, as in the figure below. 
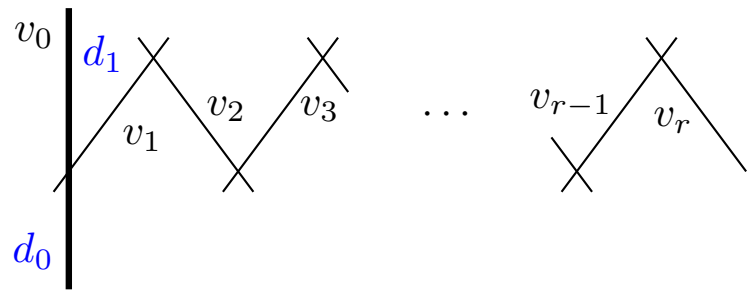
\begin{figure}[htbp]
    \centering
    \begin{subfigure}[b]{0.45\textwidth}
    \centering
    \begin{tikzpicture}[xscale=0.7, yscale=0.6, every node/.style={scale=0.8}]
  \zigzag
  \draw ($(A3)!-\ov!(P1)$) -- (P1);                       
  \node[left] at ($(A3)!0.55!(P1)$) {$v_r$};
\end{tikzpicture}
    \caption*{Outer chain of type $d_0\frac{d_1}{}$}
    \label{fig:placeholder}
    \end{subfigure}
        \begin{subfigure}[b]{0.45\textwidth}
        \centering
\begin{tikzpicture}[xscale=0.7, yscale=0.6, every node/.style={scale=0.8}]
  \zigzag
  \draw[very thick] ({8*\dx},\vup) -- ({8*\dx},\vdown);   
  \seg{A3}{P1}                                            
  \node[blue,above right=-1pt] at ($(A3)!0.30!(P1)$) {$d_r$};
  \node[left]                  at ($(A3)!0.55!(P1)$) {$v_r$};
  \node[right]      at ({8*\dx},{1.20*\dy})  {$v_{r+1}$};
  \node[blue,right] at ({8*\dx},{-0.70*\dy}) {$d_{r+1}$};
\end{tikzpicture}
    \caption*{Inner chain of type $d_0\edge{d_1}{d_r}{n}d_{r+1}$}
    \label{fig:placeholder}
    \end{subfigure}
\end{figure}
  
\noindent The multiplicities of $v_0$, $v_{r+1}$ are part of the chain data, and we view an intersection point of two (not necessarily distinct) principal components as an inner chain from $v_0$ to $v_1$ with $r=0$. We call $v_1$ the \emph{initial component} and $m_{v_1}$ the \emph{initial multiplicity} of the chain at $v_0$. In the inner case, we can also view the chain as going in the opposite direction, and call $m_{v_r}$ its initial multiplicity at $v_{r+1}$. Chains are determined by their \emph{types}. Write $d_i=m_{v_i}$ for the multiplicities of chain components. We call the chains in the figure an \emph{outer chain of type} $d_0\frac{d_1}{}$ and an \emph{inner chain of type} $d_0\edge{d_1}{d_r}{n}d_{r+1}$, where $n$ is the depth, defined as follows. 
    \[
    n=-1+\text{ number of times } \gcd(d_0,\dots,d_{r+1}) \text{ occurs in the sequence } d_0,\dots,d_{r+1}.
    \]
We refer to \cite[Definitions 4.2, 4.3 and Proposition 5.3]{Tim2025} for more details.

Given an rnc reduction type $R=(S;m;p;g;\cdot)$, one can associate a \emph{dual graph} $\Gamma$ as follows. The set $S$ is the set of vertices of $\Gamma$, and there are $v\cdot v'$ edges between the vertices $v$ and $v'$, whenever $v\neq v'$. Moreover, to each vertex $v$ we attach $p(v)-g(v)$ loops (i.e. edges from $v$ to itself). 

\begin{defn}
    Let $\Gamma$ be a graph. A \emph{path} is a sequence of edges $(e_1,\dots,e_r)$ for which there is a sequence of vertices $(v_1,\dots,v_{r+1})$ such that $e_i$ connects $v_i$ and $v_{i+1}$, and such that the $v_i$ are all distinct, except possibly $v_{r+1}=v_1$. A \emph{cycle} is a path such that $v_{r+1}=v_1$.
\end{defn}

We may view a cycle $\Delta$ given by $(e_1,\dots,e_r)$ and  with vertices $(v_1,\dots,v_{r+1})$ as a subgraph by letting $V(\Delta)=\{v_i\}_i$ and $E(\Delta)=\{e_i\}_i$. 

\begin{defn}
    A \emph{cycle} in a reduction type $R$ is a set of components in $R$ corresponding to the vertices of a cycle in the dual graph of $R$.
\end{defn}

The abelian and toric ranks of a curve $C/K$ (or its Jacobian) can be read from any rnc reduction type of $C$ \cite[Lemma 7.5.18, Theorem 7.5.19]{Liu}. Indeed, suppose $R=(S;m;p;g;\cdot)$ is an rnc reduction type arising from the special fibre of an rnc model of $C$, and let $\Gamma$ denote its dual graph. We denote by $H_1(\Gamma)$ the simplicial homology of $\Gamma$ with integer coefficients. The toric rank of $C$ is 
\[
t(C)=\operatorname{rk} H_1(\Gamma),
\]
and its abelian rank is
\[
a(C)=\sum_{v\in S} g(v).
\]
We note that this is independent of the choice of rnc model.

\subsection{Weight of components} \label{sec: weight}
\mbox{}
\vspace{5pt}

\begin{defn}[{\cite[Definition 7.3]{Tim2025}}]
    Let $v$ be a component of multiplicity $m_v$ in an rnc reduction type, and let $d_1,\dots,d_r$ denote the multiplicities of the components that intersect $v$. The \emph{weight} of $v$ is defined as 
    \[
    \operatorname{weight}(v)=\gcd(m_v,d_1,\dots,d_r).
    \]
\end{defn} 

We extend this notion to dual graphs and subgraphs of them (e.g. cycles, trees, a specific vertex...).

\begin{defn}
    Let $Z$ be a subgraph of the dual graph $\Gamma$ of an rnc reduction type. Let $v_1,\dots,v_r$ denote the vertices of $Z$. We define the \emph{weight} of $Z$ to be 
    \[
    \operatorname{weight}(Z)=\gcd(\operatorname{weight}(v_1),\dots,\operatorname{weight}(v_r)).
    \] 
    We define the \emph{weight} of an edge $e\in E(\Gamma)$ connecting the vertices $v_1$ and $v_2$ to be 
    \[
    \operatorname{weight}(e)=\gcd(m_{v_1}, m_{v_2}),
    \] 
    where $m_{v_i}$ denotes the multiplicity of $v_i$.
\end{defn}

\begin{remark} \label{weight of a chain}
    The weight of a chain is the $\gcd$ of the multiplicities of any two consecutive components, including the principal components it meets (see e.g. \cite[Lemma 5.4]{Tim2025} or \cite{Lorenzini2000}). 
\end{remark}

 \subsection{Base change and normalisation} \label{sec: base change}
\mbox{}
\vspace{5pt}

Let $\C/O_K$ be an snc model of $C/K$. 
Let $p\nmid n$, and let $K(n)/K$ be the unique degree $n$ extension of $K$ in $K^s$. Let $\C_n$ denote the normalisation of $\C\times_{O_K}O_{K(n)}$ and $\C_n\rightarrow \C$ the canonical morphism. Let $\C(n)$ denote the minimal desingularisation of $\C_n$.

\begin{prop}[{\cite[Proposition 4.1.3.2]{Halle-Nicaise}}] \label{Halle-Nicaise 1.3.2}
Let $v$ be an irreducible component of $\C_k$ and let $d$ denote its multiplicity in $\C_k$. The following properties hold:
\begin{enumerate}
    \item \label{HN 1} The scheme $u=\C_n\times_{\C} v$ is a disjoint union of smooth irreducible curves $u_{i}$. The multiplicity $d'$ of each component $u_{i}$ in $(\C_n)_k$ is given by
    \[
    d'=\frac{d}{\gcd(d,n)},
    \]
    and the morphism $\C_n\times_{\C} v \rightarrow v$ is a Galois étale cover of degree $\gcd(d,n)$ above $v-\{x_1,\dots,x_r\}$, where the $x_j$ are the points of intersection of $v$ with other irreducible components of $\C_k$. 
    \item \label{HN 2} If $v$ is a rational curve that intersects the other components of $\C_k$ in precisely one (resp. two) points, then each $u_{i}$ is a rational curve that intersects the other components of $(\C_n)_k$ in precisely one (resp. two) points. In both cases, the number of connected components of $u$ is equal to $\ell=\gcd(\operatorname{weight}(v),n)$.
    \item \label{HN 3} Each singular point of $\C_n$ is an intersection point of two distinct irreducible components of $(\C_n)_k$. Moreover, let $x$ be a point that belongs to the intersection of two distinct irreducible components $u_1$ and $u_2$ of $(\C_n)_k$ which dominate irreducible components $v_1$ and $v_2$ of $\C_k$, respectively. Let $d_1$ and $d_2$ be the respective multiplicities of $v_1$ and $v_2$ in $\C_k$. Then the special fibre of the minimal desingularisation $\D$ of the local germ $\operatorname{Spec}(\mathcal{O}_{\C_n,x})$ is a divisor with strict normal crossings whose combinatorial data (i.e. multiplicities of components and their intersection numbers; see Proposition \ref{desingularisation}) only depend on $d_1,d_2$ and $n$. Moreover, each exceptional component of $\D_k$ is a rational curve that meets the other irreducible components of $\D_k$ in precisely two points. 
    \end{enumerate}
    In particular, the $O_{K(n)}$-scheme $\C(n)$ is an snc model.
\end{prop}

We say that the irreducible components $u_{i}$ defined in the proposition \emph{lie above} $v$. We say that a point in $(\C_n)_k$ mapping to $x\in \C_k$ \emph{lies above} $x$.

\begin{fact} [{\cite[Section 4.1.2.1]{Halle-Nicaise}}] \label{number of points above intersection point}
    Let $x$ be an intersection point of two distinct components $v_1$ and $v_2$ in $\C_k$, and let $d_i$ denote the multiplicity of $v_i$. Then there are exactly $\gcd(d_1,d_2,n)$ points in $(\C_n)_k$ lying above $x$. 
\end{fact}

\begin{thm}[{\cite[Theorem 7.5]{Halle}}] \label{Halle}
    Let $C/K$ be a tame curve and suppose that $L/K$ is the minimal extension over which $C_L$ has semistable reduction. Let $v_1,\dots,v_r$ be the principal components in the reduction type of $C$ and denote by $m_1,\dots, m_r$ their respective multiplicities. Then $L/K$ is the unique extension of degree $\operatorname{lcm}(m_1,\dots,m_r)$. 
\end{thm}

\subsection{Hirzebruch-Jung continued fractions and desingularisation}
\mbox{}
\vspace{5pt}

We now review how to desingularise the normalisation $\C_n$ of $\C\times_{O_K} O_{K(n)}$ in order to get a regular model $\C(n)$. Note that even though the desingularisation described below is locally minimal, the model $\C(n)$ might not be minimal.

Let $n'$ and $r$ be positive integers, such that $0<r<n'$ and $(r,n')=1$. The Hirzebruch-Jung continued fraction expansion $[b_1,b_2,\dots,b_\lambda]_{\operatorname{HJ}}$ of $n'/r$ is given by
\[\frac{n'}{r}=b_1-\frac{1}{b_2-\frac{1}{\ldots-\frac{1}{b_\lambda}}},\] 
where $\lambda\geq 1$ and $b_i\geq 2$. 
Let $m_1,m_2\in \ZZ_{>0}$ be such that $(m_1,n')=(m_2,n')=1$ and $rm_1+m_2=0 \ \operatorname{ mod } n$. The condition $\lambda=1$ happens if and only if $r=1$, in which case we set $\mu_1=(m_1+m_2)/n'$. Otherwise, let $\mu_i$ be the unique solution to the system of equations 
\[
\begin{pmatrix}
b_1 & -1 & & & & &  \\
-1 & b_2 & -1 & & & & \\
 & -1 & b_3 & -1 & & & \\
 & & & & \ddots & &   \\
  & & & &  & -1 & b_{\lambda-1} & -1 \\
    & & & & &  & -1 & b_\lambda
\end{pmatrix}
\begin{pmatrix}
\mu_1 \\
\mu_2 \\
\vdots \\
\vdots   \\
\mu_{\lambda-1} \\
\mu_\lambda
\end{pmatrix} =
\begin{pmatrix}
m_1 \\
0\\
\vdots \\
\vdots   \\
0\\
m_2
\end{pmatrix}.
\]
We call $(\mu_1,\dots,\mu_\lambda)$ the \emph{multiplicity vector} associated to the tuple $(n',r,m_1,m_2).$

We use the formulation from \cite[Proposition 2.10.2]{Padma} for the following proposition.

\begin{prop}[{\cite[Proposition 4.4.2.5]{Halle-Nicaise}}] \label{desingularisation}
Retain the notation from Proposition \ref{Halle-Nicaise 1.3.2}(\ref{HN 3}). The special fibre of the minimal desingularisation $\D$ of the local germ $\operatorname{Spec}(\mathcal{O}_{\C_n,x})$  consists of a chain of $(-2)$-curves that connect the strict transforms of $u_1$ and $u_2$. The multiplicities of the $(-2)$-curves in the chain can be computed by the following procedure. For $i=1,2$, let $d_i'=d_i/\operatorname{gcd}(d_i,n)$ and let 
\[
n'=\frac{n\operatorname{gcd}(d_1,d_2,n)}{\operatorname{gcd}(d_1,n)\operatorname{gcd}(d_2,n)}.
\] 
If $n'=1$, set $\lambda =0$. Otherwise, let $r$ be the unique solution to the equation $rd_1'+d_2'=0 \ \operatorname{mod} n'$ such that $0<r<n'$. Let $[b_1,\dots,b_\lambda]_{\operatorname{HJ}}$ be the Hirzebruch-Jung continued fraction expansion of $n'/r$. Let $(\mu_1,\dots,\mu_\lambda)$ be the multiplicity vector associated to the tuple $(n',r,d_1',d_2')$. Then there are $\lambda$ components in the chain joining the strict transform of $u_1$ and the strict transform of $u_2$ in $\C(n)_k$. Their multiplicities in order are $\{\mu_1,\dots, \mu_\lambda\}$.  
\end{prop}

\subsection{Galois action on the special fibre $(\C_n)_k$}\label{sec: Galois theory of valuations}
\mbox{}
\vspace{5pt}

Let $p\nmid n$, and let $L=K(n)$. Then, the extension $L/K$ is Galois (see e.g. \cite{monodromy}). Let $\C$ and $\C_n$ be as in Section \ref{sec: base change}. Let $\sigma\in \Gal(L/K)$, and consider the isomorphism
\[
\C \times_{O_K} O_L \xrightarrow[]{id\times\sigma} \C\times_{O_K} O_L.
\]
Since $\C_n$ is the normalisation of $\C\times_{O_K} O_L$, this map lifts to an isomorphism $\C_n\xrightarrow{\sigma} \C_n$ by the universal property of normalisation, such that the following diagram commutes.

\[\begin{tikzcd}
	{\C_n} && {\C_n} \\
	& \C
	\arrow["\sigma", from=1-1, to=1-3]
	\arrow[from=1-1, to=2-2]
	\arrow[from=1-3, to=2-2]
\end{tikzcd}\]
Given an irreducible component $v$ in $\C_k$, Proposition \ref{Halle-Nicaise 1.3.2} says that $u=\C_n \times_{\C} v$ is a disjoint union of smooth irreducible curves $u_i$. From the diagram, we get an induced isomorphism $u\xrightarrow{\sigma} u$. Thus, the group $\Gal(L/K)$ acts on the set of irreducible components in $(\C_n)_k$ lying above $v$. Next, we show that this action is in fact transitive.

Since the models $\C$ and $\C_n$ are normal, the local rings $\mathcal{O}_{\C,\eta_v}$ and $\mathcal{O}_{\C_n,\eta_{u_i}}$ at the generic points of $v$ and $u_i$ are discrete valuation rings. Denote the corresponding normalised valuations by $w_v$ and $w_{u_i}$, and note that $w_v$ and $w_{u_i}$ are valuations on the function fields $K(\C)$ and $L(\C_n)$, respectively. Moreover, since $\C_n\to \C$ is dominant and $u_i$ maps surjectively to $v$, we have an injection of local rings $\mathcal{O}_{\C,\eta_v}\hookrightarrow \mathcal{O}_{\C_n,\eta_{u_i}}$. This means that $w_{u_i}$ extends $w_v$. Since $\Gal(L/K)$ acts on the components $u_i$, it also acts on the valuations $w_{u_i}$. 

Let $C_L=C(n)=C\times_K L$, and note that 
\[
K(C)\iso K(\C) \text{ and } L(C_L)\iso L(\C_n).
\]
Since $C$ is assumed to be smooth and geometrically connected, it is geometrically irreducible. By  (\cite[Corollary 3.2.14(d)]{Liu}), this implies that $K(C)\cap K^s=K$, so in particular $K(C)\cap L=K$. Then, the extension $L(C_L)/K(C)$ is Galois, and
\[
G:=\Gal(L(C_L)/K(C))\iso \Gal(L/K).
\]
If $w$ is a valuation on $K(C)$, recall that $G$ acts transitively on the set of extensions $w_i\mid w$. 
This implies that the action of $\Gal(L/K)$ on the valuations $w_{u_i}$ extending $w_v$ is transitive. We state this as a lemma for future reference.

\begin{lemma} \label{lem: Galois action on components}
   Let $v$ be a component in $\C_k$. Then, the group $\Gal(K(n)/K)\iso \ZZ/n\ZZ$ acts transitively on the set of components in $(\C_n)_k$ lying above $v$.
\end{lemma}

\subsection{Harmonic and cyclic dilated covers of graphs} \label{harmonic morphisms section}
\mbox{}
\vspace{5pt}

Let $\Gamma$ be a finite graph. We denote by $X(\Gamma)$ the union of the \emph{vertices} $V(\Gamma)$ and \emph{edges} $E(\Gamma)$. We endow $\Gamma$ with an orientation so that we have \emph{source} and \emph{target} maps $s,t: E(\Gamma)\to V(\Gamma)$. For simplicity, we assume that $\Gamma$ has no edge $e$ such that $s(e)=t(e)$. The \emph{tangent space} of a vertex is \[T_v(\Gamma)=\{e\in E(\Gamma) : v=s(e) \text{ or } t(e)\}.\]
A \emph{morphism} of graphs $f:\Gamma'\to \Gamma$ is a map of sets $f:X(\Gamma')\to X(\Gamma)$ that sends vertices to vertices, edges to edges, and preserves adjacency. 
A \emph{harmonic morphism} $(f,d_f)$ consists of a graph morphism $f: \Gamma' \to \Gamma$ and a \emph{dilation} assignment $d_f: X(\Gamma')\to \ZZ_{>0}$ such that 
\[
d_f(v')=\sum_{e'\in T_{v'}\Gamma'\cap f^{-1}(e)}d_f(e')
\]
for every $v'\in V(\Gamma')$ and every $e\in T_{f(v')}\Gamma$. If $\Gamma$ is connected, then the \emph{(global) degree} of $f$ is defined as 
\[
\deg (f)=\sum_{v'\in f^{-1}(v)} d_f(v')=\sum_{e'\in f^{-1}(e)}d_f(e')
\]
for any choice of $v\in V(\Gamma)$ and $e\in E(\Gamma)$. A harmonic morphism $(f,d_f)$ satisfying $d_f(x)=1$ for all $x\in X(\Gamma)$ is a covering space in the topological sense.

Let $\A$ be a finite group. An $\A$-action on $\Gamma$ is a homomorphism from $\A$ to the automorphism group $\operatorname{Aut}(\Gamma)$ such that, for every $g\in \A$ and every $e\in E(\Gamma)$, we have $g(e)=e$ only if $g(s(e))=s(e)$ and $g(t(e))=t(e)$.   

\begin{defn}[{\cite[Definition 1.4]{Len-Ulrisch-Zakharov}}] 
    Let $\Gamma$ be a graph and let $\A$ be a finite group of order $n$. A \emph{harmonic $\A$-cover} of $\Gamma$ is a harmonic morphism $f:\Gamma'\to \Gamma$ of degree $n$ together with an $\A$-action on $\Gamma'$ such that the following axioms hold:
    \begin{enumerate}[(i)]
        \item The map $f$ is $\A$-invariant, i.e. $f(g(x'))=f(x')$ and $d_f(g(x'))=d_f(x')$ for all $x'\in X(\Gamma')$ and all $g\in \A$.
        \item For all $x\in X(\Gamma)$, the group $\A$ acts transitively on the fibre $f^{-1}(x)$.
    \end{enumerate}
\end{defn}

\begin{remark}
    If $f:\Gamma'\to \Gamma$ is a harmonic $\A$-cover, then $f$ is the quotient morphism $\Gamma'\to \Gamma'/\A$, where $\Gamma'/\A$ is the quotient graph. We refer to \cite{Len-Ulrisch-Zakharov} for more details.
\end{remark}

Forgetting the $\A$-action on $\Gamma'$ in a harmonic $\A$-cover $\Gamma'\to \Gamma$ motivates the following definition. 

\begin{defn}
     Let $\Gamma$ be a graph. A \emph{cyclic dilated cover} (of \emph{degree} $n$) of $\Gamma$ is a graph $\Gamma'$ with a harmonic morphism $f:\Gamma'\to \Gamma$ arising from a harmonic $\ZZ/n\ZZ$-cover. 
\end{defn}

\subsection{Dilated cohomology}\label{sec: dilated cohomology}
\mbox{}
\vspace{5pt}

Retain the notation from Section \ref{harmonic morphisms section}. For simplicity, assume that $\Gamma$ has no multi-edges (i.e. there is at most one edge between any two vertices), and that $\A$ is a finite abelian group. It turns out that harmonic $\A$-covers with associated \emph{dilation datum} are classified by the \emph{dilated cohomology group}, which we introduce in this section. 

Choose an assignment $\ell: E(\Gamma)\to \mathbb{R}_{>0}$ of lengths to the edges of $\Gamma$, and equip $\Gamma$ with the shortest-path metric. We define the \emph{star cover}
\[
\mathcal{U}(\Gamma)=\{U_v\}_{v\in V(\Gamma)}
\]
of $\Gamma$ as follows. For each vertex $v\in V(\Gamma)$, let $U_v$ be the union of $v$ and the interiors of all edges incident to $v$. For distinct vertices $v,u$, the intersection $U_v\cap U_u$ is either the open edge connecting $v$ and $u$, if there is such an edge, or is empty otherwise. 

\begin{defn}[{\cite[Definition 2.1]{Len-Ulrisch-Zakharov}}] 
    An \emph{$\A$-dilation datum} $D$ on $\Gamma$ is a choice of a subgroup $D(x)\subseteq \A$ for every $x\in X(\Gamma)$, such that $D(e)\subseteq D(v)$ if $v=s(e) \text{ or } t(e)$. 
\end{defn}

Given a harmonic $\A$-cover $f:\Gamma'\to \Gamma$ and $x\in X(\Gamma)$, let $D(x)\subseteq \A$ denote the stabiliser of any element in $f^{-1}(x)$. Note that $D(x)$ is well-defined since $\A$ is abelian. We call $D(x)$ the \emph{dilation group} of $x$. The dilation groups of $x\in X(\Gamma)$ give rise to an $\A$-dilation datum $D$ associated to $\Gamma'\to \Gamma$. 
We say that the cyclic dilated cover arising from $\Gamma'\to\Gamma$ also has dilation datum $D$. We may sometimes omit the word \emph{dilated}, and simply call it a cyclic cover with dilation datum $D$.  

\begin{defn}[{\cite[Definition 2.2]{Len-Ulrisch-Zakharov}}] \label{defn: codilation sheaf}
    Let $D$ be an $\A$-dilation datum on $\Gamma$. We define the \emph{codilation sheaf} $\A_D$ on $\Gamma$ as follows. For a vertex $v\in V(\Gamma)$, we denote $C(v)=D(v)$. For an edge $e\in E(\Gamma)$, we denote $C(e)=D(s(e))+D(t(e))\subseteq \A$. We note that $D(e)\subseteq C(e)$ for any edge $e\in E(\Gamma)$. Now let $\mathcal{U}(\Gamma)=\{U_v\}_v$ be the star cover of $\Gamma$. The sections of $\A_D$ over the open cover and the induced intersections are
    \[
    \A_D(U_v)\iso \A/C(v) \ \text{ and } \ \A_D(U_e)\iso \A/C(e), 
    \]
    where $U_e=U_v\cap U_u$ if $e$ is the (unique) edge between $v$ and $u$. The restriction maps are induced by the inclusions $D(v)=C(v)\subseteq C(e)$ for an edge $e$ such that $v=s(e) \text{ or } t(e)$. Given a connected open set $U\subseteq \Gamma$, we set  $\A_D(U)=\A_D(U_v)$ if $v\in U\subseteq U_v$ for some vertex $v$, while $\A_D(U)=\A_D(U_e)$ if $U\subseteq U_e$ for a vertex $e$. For larger open sets, we define the space of sections via the sheaf axioms. 
    
    The \emph{dilated cohomology group} of a pair $(\Gamma, D)$ is the sheaf cohomology group $H^1(\Gamma, \A_D)$. 
\end{defn}

\begin{remark}
    The sheaf $\A_D$ depends on the choice of $\ell:E(\Gamma)\to \mathbb{R}_{>0}$, but the group $H^1(\Gamma,\A_D)$ does not.
\end{remark}

\begin{thm}[{\cite[Theorem 2.3]{Len-Ulrisch-Zakharov}}]  \label{classification of harmonic covers} 
    Let $D$ be an $\A$-dilation datum on $\Gamma$. There is a natural one-to-one correspondence between $\A_D$-torsors on $\Gamma$ and harmonic $\A$-covers of $\Gamma$ with associated $\A$-dilated datum $D$. Therefore, harmonic $\A$-covers of $\Gamma$ with $\A$-dilation datum $D$ are classified by the dilated cohomology group $H^1(\Gamma, \A_D)$.
\end{thm}

We will also use the following classification of the connected harmonic $\A$-covers of $\Gamma$.

\begin{prop} [{\cite[Proposition 4.6]{v1}}] \label{prop: connected covers}
Let $\Gamma$ be a connected graph, and let $D$ be an $\A$-dilation datum on $\Gamma$. If $\sum_{v\in V(\Gamma)}D(v)=\A$, then every harmonic $\A$-cover with dilation datum $D$ is connected. If not, then the set of disconnected $\A$-covers with dilation datum $D$ is the union of the images of the induced maps $H^1(\Gamma,\B_{D_\B})\hookrightarrow  H^1(\Gamma, \A_D)$ over all proper subgroups $\B\subsetneq \A$ such that $D(v)\subseteq \B$ for all $v\in V(\Gamma)$, where for each such $\B$, the dilation datum $D_\B$ denotes $D$ viewed as a $\B$-dilation datum. 
\end{prop}

\begin{remark}\label{rem: more general graphs}
    Everything in Sections \ref{harmonic morphisms section} and \ref{sec: dilated cohomology} holds, in more generality, for graphs with \emph{legs} and graphs with loops and multi-edges. However, one has to be more careful with certain definitions. For example, the star cover is defined via a so-called \emph{simple model} of the (metric) graph. The sheaf $\A_D$ depends on the choice of such a model, but $H^1(\Gamma,\A_D)$ does not.
    We refer to \cite{Len-Ulrisch-Zakharov} for more details. 
\end{remark}

\section{Genus formula and consequences} \label{sec: genus formula}

Let $\C$ be an snc model of $C/K$, and let $\C_n$ and $\C(n)$ be as in Section \ref{sec: base change}. Let $\varphi:\C_n\to \C$ denote the canonical map. Recall that, since $\C$ is snc, any component $v$ in the special fibre $\C_k$ has arithmetic and geometric genus $p(v)=g(v)$, so we call $g(v)$ the genus of $v$. To ease notation, we often denote $\gcd(a_1,\dots,a_r)$ by $(a_1,\dots,a_r)$. 

\begin{prop}[Genus formula] \label{genus formula}
Suppose $v$ is a component in $\C_k$ of multiplicity $m_v$ and genus $g(v)$ intersecting the rest of $\C_k$ at $r$ points. Let $d_1,\dots,d_r$ denote the multiplicities of the corresponding components. Suppose there are $\ell$ components in $(\C_n)_k$ lying above $v$, and denote them by $u_1,\dots u_\ell$. Then we obtain the following formula for the genus $g(u_j)$ for any $1\leq j \leq \ell$.
\begin{equation} \label{eqn: genus formula}
    2g(u_j)-2=\frac{(m_v,n)}{\ell}(2g(v)-2)+\sum_{i=1}^r  \frac{(m_v,d_i,n)}{\ell}\left(\frac{(m_v,n)}{(m_v, d_i, n)}-1 \right).
\end{equation}
In particular, the genus of $u_j$ is independent of $j$.
\end{prop}
\begin{proof}
    Let $u:= (\varphi^{-1}(v))_{\operatorname{red}}$. Since normalisation is a finite morphism, the map $\varphi$ induces a finite 
    morphism 
    \[
    \varphi_v:u \rightarrow v.
    \] 
    Let $x_1,\dots,x_r$ be the points of intersection of $v$ with other irreducible components of $\C_k$, and let $d_1,\dots,d_r$ denote the multiplicities of the corresponding components. Let 
    \[
    x\in v - \{x_1,\dots , x_r\}
    \] and let $x'\in (\C_n)_k$ be a point mapping to $x$. From Proposition \ref{Halle-Nicaise 1.3.2} we have that $\varphi_v:u \map v$ is étale of degree $\operatorname{deg}(\varphi_v)=(m_v,n)$ above $v-\{x_1,\dots, x_r\}$. 
    
    Assume now that $x_i\in \{x_1,\dots,x_r\}$. By Fact \ref{number of points above intersection point}, there are $(m_v, d_i, n)$ points in $(\C_n)_k$ lying above $x_i$, and let $x_i'\in (\C_n)_k$ be any such point. From the computations in \cite[Section 2.7]{Halle} and the local computations in \cite[Section 1.2]{Halle-Nicaise}, the map \[ \widehat{\OO}_{v, x_i}\map \widehat{\OO}_{u, x_i'}\] can be described by 
    \[k[[t]]\map k[[s]],\] where $t\mapsto s^{c_i}$ and $c_i=\frac{(m_v,n)}{(m_v,d_i,n)}$. In particular we see that $u$ is also regular at $x_i'$, and $\varphi_v:u \map v$ is tamely ramified with index $c_i$ at each of the $(m_v,d_i,n)$ points in $u$ lying above $x_i$. 

    Let $u_1,\dots u_\ell$ be the irreducible components of $u$. Since $\varphi_v:u\map v$ is étale of degree $\operatorname{deg}(\varphi_v)=(m_v,n)$ above $v-\{x_1,\dots, x_r\}$, 
    it follows that
    \[(m_v,n)=\sum_{1\leq j \leq \ell} [K(u_j):K(v)].\] 
    Since $\Gal(K(n)/K)$ permutes the components $u_j$, we have $(m_v,n)=\ell[K(u_j):K(v)]$ for any $1\leq j \leq \ell$. This implies that $\varphi|_{u_j}:u_j\map v$ is of degree $(m_v,n)/\ell$ and there are $(m_v,d_i,n)/\ell$ points of $u_j$ lying above $x_i$, each with ramification index $c_i$. In particular we have $\ell \mid (m_v,d_i,n)$ for every $1\leq i\leq r$.
        
    Then, the Riemann-Hurwitz formula gives that 
    \[
    2g(u_j)-2=\frac{(m_v,n)}{\ell}(2g(v)-2)+\sum_{i=1}^r \frac{(m_v,d_i,n)}{\ell}(c_i-1).
    \]
\end{proof}
\begin{remark}
    This genus formula generalises the one in \cite[Remark 10.4.8]{Liu}, where it is assumed that the multiplicities of all components in the reduction type divide $n$. 
\end{remark} 

\subsection{Number of components lying above a given component}
\mbox{}
\vspace{5pt}

In the proof of the genus formula we also showed the following result, which we state as a separate proposition for later referencing.

\begin{prop} \label{number of components divides}
    Let us keep the same hypotheses from Proposition \ref{genus formula}. The number $\ell$ of components in $(\C_n)_k$ lying above $v$ divides $(m_v,d_i,n)$ for every $1 \leq i \leq r$. Thus, we have $\ell\mid \gcd(\weight(v),n)$.
\end{prop}

\begin{prop} \label{number of components g=0}
    Let $v$ be a component with $g(v)=0$ and multiplicity $m_v$ in $\C_k$. Suppose $v$ intersects the rest of $\C_k$ at $r$ points, and let $d_1,\dots ,d_r$ denote the multiplicities of the corresponding components. Then there are 
    \(
    \gcd(m_v,d_1,\dots ,d_r,n)=\gcd(\operatorname{weight}(v),n)
    \)
    components in $(\C_n)_k$ lying above $v$.
\end{prop}
\begin{proof}
Factorise $(m_v,d_1,\dots,d_r,n)=p_1p_2\cdots p_s$ into not necessarily distinct primes. Then $p_1p_2\cdots p_s\mid n$, so consider the tower of subextensions \[K_0:=K\subset K_1\subset K_2 \subset \cdots \subset K_s \subset K(n)\] of $K(n)/K$, where for each $1\leq i \leq s$ the extension $K_i/K_{i-1}$ is defined as the unique extension of $K_{i-1}$ of degree $p_{i}$, i.e. $K_i=K(p_1\cdots p_i)$. One can check that $\C_{p_1\cdots p_{i+1}}$ is also the normalisation of $\C_{p_1\cdots p_{i}}\times_{O_{K_{i}}}O_{K_{i+1}}$.

For any $1\leq i \leq s$, we claim that for each component $v'$ in $(\C_{p_1\cdots p_{i-1}})_k$ lying above $v$, there are $p_i$ components in $(\C_{p_1\cdots p_{i}})_k$ lying above $v'$, and that all such components have genus 0. We prove this inductively. Let $\ell_1$ denote the number of components in $(\C_{p_1})_k$ lying above $v$ and let $u$ be any such component. Then the genus formula (Proposition \ref{genus formula}) gives

\[
2g(u)-2=\frac{p_1}{\ell_1}(-2)+\sum_{j=1}^r \frac{p_1}{\ell_1}\left(\frac{p_1}{p_1}-1\right)=\frac{-2p_1}{\ell_1}.
\]
Since $g(u)\geq 0$, this implies $\ell_1\geq p_1$. But $\ell_1\mid p_1$ by Proposition \ref{number of components divides}, so $\ell_1=p_1$, which implies $g(u)=0$.

For $2\leq i \leq s$, let $v'$ be a component in $(\C_{p_1\cdots p_{i-1}})_k$ lying above $v$. Let $\ell_i$ denote the number of components in $(\C_{p_1\cdots p_i})_k$ lying above $v'$, and let $u$ be any such component. By inductive assumption, we have $g(v')=0$. By Proposition \ref{Halle-Nicaise 1.3.2}, the multiplicity of $v'$ is 
\[
m_{v'}=\frac{m_v}{(m_v,[K_{i-1}:K])}=\frac{m_v}{p_1\cdots p_{i-1}},
\] 
but since $p_1\cdots p_s \mid m_v$, we have that $(m_{v'},p_i)=(m_v, p_i)=p_i$. A similar argument involving the multiplicities $d_j'$ of the components intersecting $v'$ shows $(m_{v'},d_j',p_i)=(m_v, d_j,p_i)=p_i$. The genus formula then gives
\[
2g(u)-2=\frac{p_i}{\ell_i}(-2)+\sum_j \frac{p_i}{\ell_i}\left(\frac{p_i}{p_i}-1\right)=\frac{-2p_i}{\ell_i}.
\]
Arguing similarly as before gives $\ell_i=p_i$ and $g(u)=0$, and this proves the claim. 
Let $\ell$ denote the number of components in $(\C_n)_k$ lying above $v$. Then $\ell_1\cdots \ell_s \mid \ell$, 
and Proposition \ref{number of components divides} implies $\ell\mid (m_v,d_1,\dots,d_r,n)$. Since $\ell_1\cdots \ell_s=p_1\cdots p_s=(m_v,d_1,\dots,d_r,n)$ it follows that $\ell= (m_v,d_1,\dots,d_r,n)$.\\
\end{proof}

\subsection{Contractible principal components}
\mbox{}
\vspace{5pt}

We can also use the genus formula to understand the behaviour of principal components under tame extensions. 

We will need a lemma regarding how the weight changes under tame base change. Let $\C$ be an rnc model of $C$, and let $\C_n$ and $\C(n)$ be as before. Let $\varphi: \C_n\to \C$ and $\rho: \C(n)\to \C_n$ be the canonical morphisms. Let $\Gamma$ and $\Gamma(n)$ denote the dual graphs of $\C_k$ and $\C(n)_k$, respectively. We identify subgraphs of $\Gamma$ and $\Gamma(n)$ with the corresponding sets of components and intersection points in $\C_k$ and $\C(n)_k$, respectively.

\begin{lemma}
    Let $w$ denote the weight of a subgraph $Z$ of $\Gamma$. Let $Z(n)$ denote the subgraph lying above $Z$ under the map $\varphi \circ\rho$.
    Then, 
    \[
    \weight(Z(n))=w/(w,n).
    \]
    Moreover, the weight of $Z(n)$ is preserved under contraction of $(-1)$-curves. \label{weight lemma}
\end{lemma}
\begin{proof}
    Suppose that $\C$ is an msnc model. 
    Let $w$ be the weight of a subgraph $Z$ of $\Gamma$, and let $d_1,\dots,d_r$ be the multiplicities of the vertices of $Z$. Then $w=(d_1,\dots,d_n)$. Let $e_i=d_i/w\in \ZZ$ for each $1\leq i\leq r$. We have $(e_1,\dots,e_r)=1$. Identifying $Z$ with the corresponding set of components and intersection points in $\C_k$, let $Z_n$ denote its preimage under the map $\varphi$. By Proposition \ref{Halle-Nicaise 1.3.2}\eqref{HN 1}, the multiplicity of any component in $Z_n$ lying above a component in $Z$ of multiplicity $d_i$ is $d_i/(d_i,n)$. Suppose there is a singular point $x$ at the intersection of two irreducible components of multiplicities $m_i$ and $m_{i+1}$. Then the minimal desingularisation of the local germ $\operatorname{Spec}(\mathcal{O}_{\C_n,x})$ consists of a chain of $(-2)$-curves of multiplicities $\mu_1,\dots,\mu_\lambda$ given by Proposition \ref{desingularisation}. By \cite[Lemma 2.9.2]{Padma} we have $(m_i,\mu_1,\dots,\mu_\lambda,m_{i+1})=(m_i,m_{i+1})$. Therefore, we have
    \[ \weight(Z(n))=\left(\frac{d_1}{(d_1,n)},\dots,\frac{d_r}{(d_r,n)}\right)=\left(\frac{e_1w}{(e_1w,n)},\dots,\frac{e_rw}{(e_rw,n)}\right).\]
    Let $a=(w,n)$ and write $w'=w/a$ and $n'=n/a$. Then 
    \[\weight(Z(n))=\left(\frac{e_1w}{(e_1w'a,n'a)},\dots,\frac{e_rw}{(e_rw'a,n'a)}\right)=\frac{w}{a}\left(\frac{e_1}{(e_1w',n')},\dots,\frac{e_r}{(e_rw',n')}\right)=\frac{w}{a},\]
    where the last equality follows from $(e_1,\dots,e_r)=1$. 
    
    Suppose $v$ is a $(-1)$-curve that intersects the rest of the special fibre at $s$ points, and let $a_1,\dots,a_s$ denote the multiplicities of the corresponding components. Then the multiplicity $d$ of $v$ satisfies $d=a_1+\cdots+a_s$. Then $(a_1,\dots,a_s,d)=(a_1,\dots,a_s)$, so the weight is also preserved after contracting $(-1)$-curves. Since the weight is preserved under contraction (or addition) of $(-1)$-curves, the result also holds if $\C$ is an rnc model. \\
\end{proof}

Next, let us recall the \emph{Euler characteristic} $\chi_v$ of a component $v$, defined in the following theorem.

\begin{thm}[{\cite[Theorem D]{Tim2025}}] \label{Euler characteristic} For each principal component $v$ in the reduction type of a curve $C/K$ there is a unique invariant $\chi_v\in \ZZ$ depending on the multiplicity $m_v$ of $v$, its geometric genus $g(v)$ and the sets $\OO_v$, $\I_v$ of initial multiplicities of the outer chains and inner chains from $v$, respectively. It is given by
\[
\chi_v=m_v(2-2g(v)-|\OO_v|-|\I_v|)+\sum_{o_i\in \OO_v}\gcd(m_v, o_i).
\]
If $v_1,\dots,v_n$ are the principal components in the reduction type of $C$, then we have
\[
\chi_{v_1}+\cdots +\chi_{v_n}=2-2g(C).
\]
\end{thm}

\begin{remark} \label{rem: euler characteristic}
     By Remark \ref{rem: msnc mrnc}, a component $v$ in the mrnc model of a curve $C$ is principal if and only if the unique component $v'$ in the msnc model of $C$ mapping surjectively to $v$ is principal, and moreover $\chi_v=\chi_{v'}$. 
\end{remark}

Using the genus formula we will see how the Euler characteristic of a principal component behaves under tame base change.

Let $\C(n)_{\operatorname{min}}$ be the msnc model of $C(n)/K(n)$. If $\C$ is the msnc model of $C$, then we have a map $\C(n)\to \C(n)_{\operatorname{min}}$ by Proposition \ref{Halle-Nicaise 1.3.2}. Let $u$ be a component in $(\C(n)_{\operatorname{min}})_k$, and let $u'$ be the unique component in $\C(n)_k$ mapping surjectively to $u$ under the map $\C(n)\to \C(n)_{\operatorname{min}}$. Let $v$ be the image of $u'$ under the map $\C(n)\to \C$. Then, either $v$ is a point or it is a component in $\C_k$. If $v$ is the latter, by abuse of notation we say that $u$ \emph{lies above} $v$, even though there is not necessarily a map $\C(n)_{\operatorname{min}}\to \C$. We note that, if there are $\ell$ components in the reduction type of $C(n)$ lying above $v$, then there are $\ell$ components in $(\C_n)_k$ lying above $v$. 

\begin{prop} \label{euler characteristic under tame extns}
    Let $C/K$ be a curve. Let $u$ be a principal component in the reduction type of $C(n)$ lying above a principal component $v$ in $\C_k$, and denote by $\ell$ the number of such components. Then $\chi_v\leq\ell\chi_u$. 
\end{prop}
\begin{proof}
    First, note that it is enough to prove the result for $n$ prime. Moreover, by Remark \ref{rem: euler characteristic}, we may assume that $\C$ and the reduction type $R_n$ of $C(n)$ are msnc. Let $u_1,\dots,u_\ell$ be the components in the reduction type of $C(n)$ lying above $v$, and assume that at least one of them is principal. Then, by Proposition \ref{Halle-Nicaise 1.3.2}, the component $v$ is also principal and there are $\ell$ components $v_1,\dots,v_\ell$ in $(\C_n)_k$ lying above $v$. Up to rearranging the components, each $u_i$ is the image under $\C(n)\to \C(n)_\text{min}$ of a component $v_i'$ in $\C(n)_k$ that maps surjectively to $v_i$ under the map $\C(n)\to \C_n$. Then, by Lemma \ref{lem: Galois action on components}, each $u_i$ is a principal component and $\chi_{u_i}=\chi_{u_j}$ for any $1\leq i,j\leq \ell$, so let $u=u_1$. By Proposition \ref{Halle-Nicaise 1.3.2}, the multiplicity of $u$ is $m_u=\frac{m_v}{(m_v,n)}$. Hence the Euler characteristic of $u$ is
    \[
    \chi_{u}=\frac{m_v}{(m_v,n)}(2-2g(u)-|\OO_{u}|-|\I_{u}|)+\sum_{o_i\in \OO_u}(m_u,o_i).
    \]  
    Let $\mathcal{O}_v$, $\mathcal{I}_v$ be as in Theorem \ref{Euler characteristic}. Note that $g(u)=g(v_1)$. By the genus formula,
\begin{align*}
2-2g(u) = & \ \frac{(m_v,n)}{\ell}(2-2g(v))-\sum_{d_i\in \OO_v\cup\I_v}\frac{(m_v,d_i,n)}{\ell}\left(\frac{(m_v,n)}{(m_v,d_i,n)}-1\right) \\
= & \ \frac{(m_v,n)}{\ell}(2-2g(v)-|\OO_v|-|\I_v|)+\frac{1}{\ell}\sum_{d_i\in \OO_v\cup\I_v} (m_v,d_i,n).
\end{align*}
Thus, we can rewrite the Euler characteristic of $u$ as follows.
\begin{equation}\label{eqn: euler characteristic of u wrt v}
    \chi_u=\frac{1}{\ell}\chi_v-\frac{1}{\ell}\sum_{d_i\in \OO_v}(m_v,d_i)+\frac{m_v}{(m_v,n)}\left(\frac{1}{\ell}\sum_{d_i\in \OO_v\cup \I_v} (m_v,d_i,n)-|\OO_u|-|\I_u|\right)+\sum_{o_i\in \OO_u}(m_u,o_i).
\end{equation}
Let $S\subseteq \OO_v\cup \I_v$ be the subset of initial multiplicities of chains $Z$ such that there exists a chain in $R_n$ lying above $Z$. By Proposition \ref{Halle-Nicaise 1.3.2} (or Fact \ref{number of points above intersection point} in the case of inner chains obtained from the intersection of two principal components), for each $d_i\in S$, there are $(m_v,d_i,n)$ chains in $R_n$ lying above the chain $Z_i$ corresponding to $d_i$ in the reduction type $R$ of $C$. Thus, there are $(m_v,d_i,n)/\ell$ chains in $R_n$ from $u$ lying above $Z_i$. This means that
\[
|\OO_u|+|\I_u|=\sum_{d_i\in S} \frac{(m_v,d_i,n)}{\ell}.
\]
Let $o_i\in \OO_u$. The corresponding chain in $R_n$ lies above a chain $Z_i$ in $R$, with initial multiplicity $d_i\in S$. By Lemma \ref{weight lemma} and Remark \ref{weight of a chain}, we have that 
\[
(m_u,o_i)=\frac{\weight(Z_i)}{(\weight(Z_i),n)}=\frac{(m_v,d_i)}{(m_v,d_i,n)}.
\]
This yields the following inequality, which is not necessarily an equality since, a priori, we do not know whether $Z_i$ is outer or inner. 
\[
\sum_{o_i\in \OO_u} (m_u,o_i)\geq \sum_{d_i\in S\cap \OO_v} \frac{(m_v,d_i,n)}{\ell}\frac{(m_v,d_i)}{(m_v,d_i,n)} = \sum_{d_i\in S\cap \OO_v} \frac{(m_v,d_i)}{\ell}
\]
Then, substituting the above expressions in \eqref{eqn: euler characteristic of u wrt v}, simplifying, and multiplying through by $\ell$ gives
\begin{align*}
    \ell\chi_u \geq & \  \chi_v -\sum_{d_i\in \OO_v\backslash S} (m_v,d_i)+\frac{m_v}{(m_v,n)}\sum_{d_i\in (\OO_v\cup \I_v) \backslash S} (m_v,d_i,n) \\
    \geq &  \ \chi_v +\sum_{d_i\in \OO_v\backslash S} \left(\frac{m_v}{(m_v,n)}(m_v,d_i,n)-(m_v,d_i) \right) 
\end{align*}
Let $d_i\in \OO_v$. Recall that we assumed that $n$ is prime. If $n\nmid m_v$, or $n\mid m_v$ and $n\mid d_i$, then 
\[
\frac{m_v}{(m_v,n)}(m_v,d_i,n)-(m_v,d_i) = m_v - (m_v,d_i) \geq 0.
\]
If $n\mid m_v$ and $n\nmid d_i$, then $(m_v,d_i)=(m_v/n,d_i)$, so we get again
\[
\frac{m_v}{(m_v,n)}(m_v,d_i,n)-(m_v,d_i) = \frac{m_v}{n}-\left(\frac{m_v}{n},d_i\right)\geq 0.
\]
Hence we obtain $\ell \chi_u \geq \chi_v$. \\
\end{proof}

\begin{defn}[\cite{Tim2025}] \label{d-tails}
A component $v$ in a reduction type is called a \emph{D-tail (of weight m)} if it has $g(v)=0$, $m_v=2m$, $\OO_v=\{m,m\}$ and $\I_v=\{2m\}$ for some $m\geq 1$. 
\end{defn}

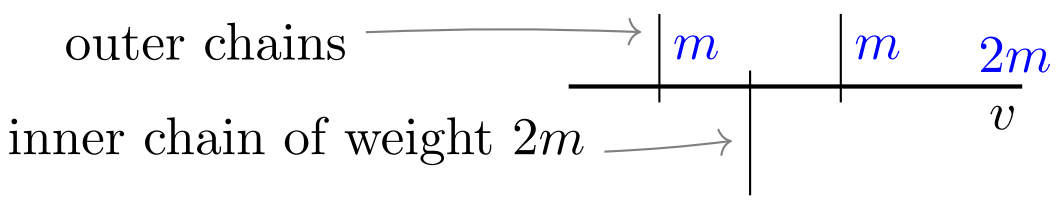
\begin{center}
\begin{tikzpicture}[scale=0.6]
  \draw[thick] (0,0) -- (5,0) node[below right=0pt and -10pt] {$v$} node[above right=-1pt and -12pt,blue] {$2m$};
  \draw[shorten <=-3pt] (1,0) -- (1,0.8) node[above right=-12pt and -1pt,blue] {$m$};
  \draw[shorten <=-3pt] (3,0) -- (3,0.8) node[above right=-12pt and -1pt,blue] {$m$};
  \draw[shorten <=-3pt] (2,0) -- (2,-1.2) node[blue,below right=-12pt and -1pt] {};
  \node[align=center] (inner) at (-3,-0.6) {inner chain of weight $2m$};
  \draw[->,black!50] (inner) edge[bend right=2] (1.8,-0.6);
  \node (outer) at (-4,0.5) {outer chains};
  \draw[->,black!50] (outer) edge[bend left=2] (0.8,0.6);
  \node at (10,0) {(D-tail of weight $m$)};
\end{tikzpicture}
\end{center}

\begin{thm} \label{non-contractibility of principal components} 
    Let $C/K$ be a curve with $g(C)\geq 2$. If $v$ is a principal component in the reduction type of $C$, then there is a component $u$ in the reduction type of $C(n)$ lying above $v$. Moreover, if the component $u$ is not principal, then $v$ is a D-tail.
\end{thm}
\begin{proof}
    By Remark \ref{rem: msnc mrnc}, it is enough to show the result for msnc reduction types. Let $R$ and $R_n$ denote the msnc reduction types of $C$ and $C(n)$, respectively. Let $u_1^1,\dots,u_{\ell_1}^1,\dots,u_1^s,\dots,u_{\ell_s}^s$ be the principal components in $R_n$, and for each $1\leq i\leq s$ and each $1\leq j\leq \ell_i$, suppose $u_j^i$ lies above a component $v_i$ in $R$. Then, by Proposition \ref{Halle-Nicaise 1.3.2}, the component $v_i$ is principal. By Lemma \ref{lem: Galois action on components}, we have $\chi_{u^i_j}=\chi_{u^i_{j'}}$ for any $1\leq j,j'\leq \ell_i$. Then, by Theorem \ref{Euler characteristic} and Proposition \ref{euler characteristic under tame extns} we have
    \[
    2-2g(C)=\sum_{i=1}^s\sum_{j=1}^{\ell_i}\chi_{u_j^i}=\sum_{i=1}^s\ell_i\chi_{u_1^i}\geq \sum_{i=1}^s \chi_{v_i}.
    \]
    Let $v_1,\dots,v_s,\dots,v_r$ denote the principal components in $R$. Then, for every $s+1\leq i\leq r$, there are no principal components in $R_n$ lying above $v_i$. Moreover, by Theorem \ref{Euler characteristic} again, we have
    \[
    2-2g(C)=\chi_{v_1}+\cdots+ \chi_{v_s}+\dots+\chi_{v_r}.
    \]
    Thus $\sum_{i=s+1}^r\chi_{v_i}\geq 0$. By \cite[Theorem E]{Tim2025}, we have $\chi_v\leq 0$ for any principal component $v$, and that the only $v$ with $\chi_v=0$ in the reduction type of a curve of genus at least 2 are D-tails. This implies that $v_i$ is a D-tail for every $s+1\leq i\leq r$. Note that this also shows that, if there are no components in $R_n$ lying above a principal component $v_i$ in $R$, then $v_i$ is a D-tail. 
    
    Thus, it remains to show that any D-tail in $R$ also has a component in $R_n$ lying above it. Note that it suffices to show the result for $n$ prime. Suppose that $v$ is a D-tail of weight $m$ in $R$. Since $g(C)\geq 2$, the inner chain in $\C_k$ from $v$ intersects another principal component $v'$ that is not a D-tail. We have already shown that this implies that there are components in $R_n$ lying above $v'$. Moreover, by Remark \ref{weight of a chain}, the weight of the inner chain $Z_0$ from $v$ to $v'$ is $2m$. Let us denote the outer chains from $v$ by $Z_1$ and $Z_2$, each of which consist of only one component $v_1$ or $v_2$. 
    
    If $n\nmid 2m$, then by Proposition \ref{number of components g=0}, there is exactly one component $u$ in $(\C_n)_k$ lying above $v$, one inner chain $Z_0'$ lying above $Z_0$, and one outer chain $Z_i'$ lying above $Z_i$ for each $i=1,2$. By Proposition \ref{Halle-Nicaise 1.3.2}, the multiplicity of $u$ or any component in $Z_0',Z_1'$ or $Z_2'$ is the same as that of the component in $\C_k$ they lie above of. If the component in $\C(n)_k$ lying above $u$ gets contracted under the map $\C(n)\to \C(n)_\text{min}$, then the chains $Z_1'',Z''_2$ in $\C(n)_k$ lying above $Z_1',Z_2'$ under the map $\C(n)\to \C_n$ must get contracted under $\C(n)\to \C(n)_\text{min}$ first. For $i=1,2$, the chain $Z_i'$ consists of only one component $u_i$ of multiplicity $m$ intersecting $u$. By Proposition \ref{desingularisation}, the chain $Z_i''$ consists of $u_i$ and a chain of $(-2)$-curves between $u_i$ and $u$. Thus, if $Z_i''$ gets contracted, then the component $u_i$ must be a $(-1)$-curve. This implies that all components in the chain $Z_i''$ have multiplicity $m$. Then, all components get contracted except for one, since $u$ has multiplicity $2m$. This contradicts that $Z_i''$ gets contracted, so there is a component in $R_n$ lying above $v$.

    Suppose that $n\mid m$. By Proposition \ref{number of components g=0}, there are $n$ components in $(\C_n)_k$ lying above each of $v,v_1$ and $v_2$ and above each component in the chain $Z_0$. The multiplicities of those lying above $v_1$ and $v_2$ are $m/n$, and the multiplicities of those lying above $v$ and the components in $Z_0$ are $2m/n$. Moreover, by Proposition \ref{desingularisation}, the intersections of all these components in $(\C_n)_k$ are already regular. Together with Lemma \ref{lem: Galois action on components}, this implies that there are $n$ D-tails in $\C(n)_k$ lying above $v$. In particular, they do not get contracted under the map $\C(n)\to \C(n)_\text{min}$.

    Lastly, suppose that $n=2$ and $n\nmid m$. Then there is one component $u$ in $(\C_n)_k$ lying above $v$, one component lying above each of $v_1$ and $v_2$, and two inner chains $Z_0^1, Z_0^2$ lying above $Z_0$. The components in $Z_0^1,Z_0^2$ and $v$ have multiplicity $m$. For $i=1,2$, the inner chain $Z_0^i$ goes from $u$ to a component $v_i'$ in $(\C_n)_k$ lying above $v'$, where $v_1'$ and $v_2'$ need not be distinct. By Proposition \ref{desingularisation}, the intersections of $u$ and any components in $Z_0^i$ are regular, so we identify these components with their preimages under $\C(n)\to \C_n$. Then, all components in $Z_0^i$, as well as $u$, are $(-2)$-curves. Since $v_i'$ do not get contracted under $\C(n) \to \C(n)_\text{min}$, neither do $Z_0^i$ or $u$. \\ 
\end{proof} 

\begin{remark}
    The above result implies that, outer (resp. inner) chains in $R_n$ lie above outer (resp. inner) chains in $R$.
\end{remark}

\section{Main results} \label{sec: main results}

\subsection{Results on dilated cohomology} \label{sec: results on dilated cohomology}
\mbox{}
\vspace{5pt}

In this section, we compute the dilated cohomology $H^1(\Gamma,\A_D)$ of a graph $\Gamma$, assuming $\A$ is cyclic. By Remark \ref{rem: more general graphs}, we may assume without loss of generality that $\Gamma$ has no loops or multi-edges (i.e. we take a simple model of $\Gamma$). 

Let $\A$ be a cyclic group of order $n$. We write $P(n)$ for the set of prime divisors of $n$. Then, by the Chinese remainder theorem we have a decomposition
\[
\A\iso \bigoplus_{q\in P(n) } \A_q,
\]
where $\A_q$ is a cyclic $q$-group for each prime $q\mid n$. Similarly, for each vertex or edge $x\in X(\Gamma)$, the dilation group $D(x)\subseteq \A$ also decomposes into (possibly trivial) cyclic $q$-groups $D(x)_q\subseteq \A_q$ as follows.
\[
D(x)\iso \bigoplus_{q\in P(n)} D(x)_q.
\]
For $C(x)$ as in Definition \ref{defn: codilation sheaf}, we define the cyclic $q$-groups $C(x)_q$ similarly. Consider the codilation sheaf $\A_{q,D_q}$ on $\Gamma$, where $D_q$ is the $\A_q$-dilation datum given by $D_q(x):=D(x)_q$ for each $x\in X(\Gamma)$. To ease notation, we also denote this sheaf by $\A_{q,D}$. Therefore, 
\[
\A_D(U_x)=\A/C(x)\iso \bigoplus_{q\in P(n) } \A_q/C(x)_q \iso \bigoplus_{q\in P(n) } \A_{q,D}(U_x),
\]
for every $x\in X(\Gamma)$, so the sheaves $\A_D$ and $\bigoplus_{q\in P(n) } \A_{q,D}$ are isomorphic on the level of stalks. Thus, we have an isomorphism of sheaves $\A_D\iso \bigoplus_{q\in P(n) } \A_{q,D}$. This gives rise to the following  decomposition of the dilated cohomology group.
\begin{equation} \label{decomposition of cohomology}
    H^1(\Gamma,\A_D)\iso H^1(\Gamma,\bigoplus_{q\in P(n) } \A_{q, D})\iso \bigoplus_{q\in P(n) } H^1(\Gamma, \A_{q,D}),
\end{equation}
and each $H^1(\Gamma, \A_{q,D})$ is a $q$-group.

If $\Gamma$ is a tree, then $H^1(\Gamma, \A_D)=0$ for any $\A$-dilation datum $D$ (see \cite[Proposition 3.22]{v1}). The next proposition gives a formula for the dilated cohomology of any graph that is not a tree. 

\begin{prop}\label{cohomology of a graph}
    Suppose that $\Gamma$ is a connected graph and $D$ an $\A$-dilation datum on $\Gamma$. Let $\beta$ be the Betti number of $\Gamma$, and assume $\beta\neq 0$. For each prime $q\mid n$, choose $e^q_1,\dots, e^q_{\beta}\in E(\Gamma)$ iteratively as follows. Let $e^q_i$ be an edge such that the subgraph $\Gamma^q_i=\Gamma\backslash\{e^q_1,\dots,e^q_i\}$ is connected and has Betti number $\beta-i$, and
    \[
    C(e^q_i)_q=\operatorname{max}\{C(e)_q: e\in E(\Gamma) \text{ is contained in a cycle of } \Gamma^q_{i-1}\},
    \] 
    where $\Gamma_0^q=\Gamma$. Then $T^q:=\Gamma^q_\beta$ is a spanning tree for $\Gamma$. For each $1\leq i\leq \beta$, let $\Delta^q_i$ be the unique cycle in the subgraph $T^q\cup {e^q_i}\subseteq \Gamma$.
    Then the dilated cohomology group of the pair $(\Gamma, D)$ is 
    \[
    H^1(\Gamma, \A_D)\iso \bigoplus_{q\in P(n) } \bigoplus_{i=1}^\beta \A_q/C(e_i^q)_q\iso \bigoplus_{q\in P(n) } \bigoplus_{i=1}^\beta \A_q/\sum_{v\in V(\Delta^q_i)} D(v)_q.
    \]
    In particular, if $\Gamma$ has a unique cycle $\Delta$, then
    \[
    H^1(\Gamma, \A_D)\iso \A/\sum_{v\in V(\Delta)} D(v).
    \]
\end{prop}
\begin{proof}
    By the decomposition of the dilated cohomology in Equation \eqref{decomposition of cohomology}, it suffices to show the result for $\A$ a cyclic $q$-group. To ease notation, we drop the subindeces and superindeces $q$. In this case, for each $1\leq i\leq \beta$, the group $\sum_{v\in V(\Delta_i)} D(v)$ is the largest $q$-group among the $D(v)$ for $v\in V(\Delta_i)$, so in particular there is $u_i\in V(\Delta_i)$ such that 
    \[
    D(u_i)=\sum_{v\in V(\Delta_i)} D(v).
    \] 
    Then $e_i$ is an edge such that $C(e_i)=D(u_i)$, so $\A_D(U_{e_i})= \A/D(u_i)$. Consider the spanning tree $T=T^q$ defined above. Then, \cite[Remark 4.3]{v1} shows that any element $[\eta]\in H^1(\Gamma, \A_D)$ can be represented by a cochain $(\eta^e)_e\in C^1(\Gamma, \A_D)=\bigoplus_{e\in E(\Gamma)} \A/C(e)$ such that $\eta^e=0$ for all $e\in E(T)$. Since $\eta^{e_i}\in \A_D(U_{e_i})=\A/D(u_i)$, this gives rise to a surjective map 
    \begin{equation} \label{cycle map}
    \bigoplus_i\A/D(u_i) \twoheadrightarrow H^1(\Gamma, \A_D).
    \end{equation}
    It remains to show that this map is injective. Suppose $\eta=(\eta^e)_e,\tau=(\tau^e)_e\in C^1(\Gamma, \A_D)$ satisfy $\eta^e=\tau^e=0$ for all $e\in E(T)$ and $[\eta]=[\tau]\in H^1(\Gamma, \A_D)$. We have 
    \[
    \tau - \eta \in \operatorname{im}(d:C^0(\Gamma,\A_D)\to C^1(\Gamma,\A_D)),
    \]
    so there is $(g_v)_v\in C^0(\Gamma, \A_D)$ such that $d(g_v)_v=\tau-\eta$. In particular, for all $e\in E(T)$, we have 
    \[
    g_{t(e)}-g_{s(e)}=\tau^e-\eta^e=0 \mod  D(t(e))+D(s(e)).
    \]
    For each $1\leq i\leq \beta$, this implies that $g_{t(e)}-g_{s(e)}=0 \mod  D(u_i)$ for all $e\in E(\Delta_i)\backslash\{e_i\}$. 
   Let $(g^e)_e:=d(g_v)_v$ and by abuse of notation, for $e\in E(\Delta_i)\backslash\{e_i\}$ we denote the image of $g^e$ under the quotient map $$\A/(D(t(e))+D(s(e)))\to \A/D(u_i)$$ also by $g^e$. Consider the following telescoping sums in $\A/D(u_i)$.
    \[
    \sum_{e\in E(\Delta_i)}g^e=\sum_{e\in E(\Delta_i)}g_{t(e)}-g_{s(e)}=0 \mod D(u_i).
    \]
    It follows that $g_{t(e_i)}-g_{s(e_i)}=0 \mod D(u_i)$, so $\tau^{e_i}-\eta^{e_i}=0 \mod D(u_i)$. Thus $\tau=\eta\in \bigoplus_i \A/D(u_i)$.
    
    The result in the case where $\Gamma$ has a unique cycle $\Delta$ follows from the isomorphism
    \[
    \A/\sum_{v\in V(\Delta)} D(v)\iso \bigoplus_{q\in P(n) } \A_q/\sum_{v\in V(\Delta)} D(v)_q.
    \]
    
\end{proof}

\begin{cor} \label{superdilated}
    Let $\Gamma$ be a connected graph and let $D$ be an $\A$-dilation datum on $\Gamma$. If every cycle in $\Gamma$ contains a vertex $v$ with dilation group $D(v)=\A$, then $H^1(\Gamma, \A_D)=0$. Moreover, if we assume that $\Gamma$ is not a tree and $\A$ is a cyclic $q$-group for some prime $q$, then the converse is true. 
\end{cor}
 \begin{proof}
 This follows directly from Proposition \ref{cohomology of a graph}.

\end{proof}

\subsection{Results on cyclic dilated covers of graphs} \label{sec: results on cyclic covers}
\mbox{}
\vspace{5pt}

In this section, we characterise when there is a unique cyclic dilated cover of degree $n$ of a graph in terms of the dilated cohomology or the homology of the graph. We continue to assume $\A=  \ZZ/n\ZZ$ and that $\Gamma$ has no loops or multi-edges. 

\begin{lemma} \label{lem: cyclic submodule}
    Let $\A= \ZZ/n\ZZ$ and let $H$ be a finite $\A$-module. Two elements $g,h\in H$ generate the same cyclic submodule of $H$ if and only if there is a unit $u\in \A^\times$ such that $g=uh$. 
\end{lemma}
\begin{proof}
    Suppose there is a unit $u\in \A^\times$ such that $g=uh$. Then $g$ is contained in the submodule $\{ah:a\in \A\}$ generated by $h$. Similarly, the element $h=u^{-1}g$ is contained in the submodule generated by $g$, so $g$ and $h$ generate the same cyclic submodule of $H$.
    For the converse, suppose that $g$ and $h$ generate the same cyclic submodule of $H$. Then there are $a,b\in \A$ such that $g=ah=abg$, so $(ab-1)g=0$. Since $H$ is an abelian group, let $d$ be the order of $g$ and $h$. Then $ab\equiv 1 \text{ mod } d$. Since $d\mid n$, let us write $\bar{a}$ for the image of $a$ in $\ZZ/d\ZZ$. Since the induced map on units
    \[
    (\ZZ/n\ZZ)^\times\to (\ZZ/d\ZZ)^\times
    \]
    is surjective, there is a lift $u\in \A^\times$ of $\bar{a}$. Then $u\equiv a \text{ mod } d$, so $(u-a)h=0$. Therefore $g=ah=uh$. 
    
    \end{proof}

\begin{prop} \label{cyclic subgroups give isomorphic covers}

    Let $\A=\ZZ/n\ZZ$. Then, there is a surjection of sets
    \[
    \{\text{cyclic subgroups of $H^1(\Gamma, \A_D)$}\} \twoheadrightarrow \{\text{cyclic covers of $\Gamma$ of degree $n$ with dilation datum $D$}\}
    \]
    arising from the correspondence between $H^1(\Gamma,\A_D)$ and harmonic $\A$-covers of $\Gamma$ with dilation datum $D$. 
\end{prop}
\begin{proof}
    By Theorem \ref{classification of harmonic covers} there is a one-to-one correspondence
    \begin{equation}\label{correspondence}
    H:=\{\A_D\text{-torsors on $\Gamma$}\}\xleftrightarrow{1:1} \{\text{harmonic $\A$-covers of $\Gamma$ with $\A$-dilation datum $D$}\}=:\text{Cov}
    \end{equation}
    as sets. An element $f:\Gamma'\to \Gamma\in \text{Cov}$ is endowed with an $\A$-action preserving and acting transitively on fibres, and let us denote it by $a:\A\to \Aut(\Gamma')$. Then, any $u\in \Aut(\A)$ endows $\Gamma'$ with an action $a\circ u: \A\to \Aut(\Gamma')$, which also preserves and acts transitively on fibres. Therefore, the group $\Aut(\A)$ acts on $\text{Cov}$, and the correspondence \eqref{correspondence} induces an action of $\Aut(\A)$ on $H$. Thus, there is a one-to-one correspondence between the sets of orbits
    \[
    H/\Aut(\A) \xleftrightarrow{1:1} \text{Cov}/\Aut(\A).
    \]
    The different actions $a$ and $a\circ u$ on $\Gamma'$ may give rise to different harmonic $\A$-covers, but they give the same cyclic dilated cover $\Gamma'\to \Gamma$. Therefore we have a surjection of sets
    \[
    \phi: \text{Cov}/\Aut(\A)\twoheadrightarrow  \{\text{cyclic covers of $\Gamma$ of degree $n$ with dilation datum $D$}\}.
    \]
    Lastly, we show that there is a one-to-one correspondence
    \begin{equation*} \label{cyclic subgroups of H1}
        \{\text{cyclic subgroups of $H^1(\Gamma, \A_D)$}\} \xleftrightarrow{1:1} H/\Aut(\A).
    \end{equation*}
    We recall that the set $H$ is in one-to-one correspondence with the sheaf cohomology group $H^1(\Gamma, \A_D)$. More explicitly, as explained in \cite[page 9]{Len-Ulrisch-Zakharov}, any $\A_D$-torsor $\mathcal{T}$ determines an element $[(\eta^e)_{e\in E(\Gamma)}]$ of the 1st-Čech cohomology group 
    \[
    \check{H}:=\check{H}^1(\mathcal{U}(\Gamma),\A_D)\iso H^1(\Gamma, \A_D),
    \]
    and the construction can be reversed to obtain $\mathcal{T}$ from $(\eta^e)_e$. Thus, there is a surjective map 
    \[
    H\twoheadrightarrow\{\text{cyclic subgroups of $\check{H}$}\} \iso \{\text{cyclic subgroups of $H^1(\Gamma, \A_D)$}\}, 
    \]
    given by sending $\mathcal{T}$ to the subgroup generated by $[(\eta^e)_e]\in \check{H}$.
    Note that $\check{H}$ is an $\A$-module. 
    For an element $u\in \Aut(\A)\iso \A^\times$, let us denote also by $u$ the corresponding element in $\A^\times$. 
    Then, one can check that $u[(\eta^e)_e]=[u(\eta^e)_e]$ is an element in $\check{H}$ corresponding to the torsor $u\cdot \mathcal{T}$ obtained from the induced action of $\Aut(\A)$ on $H$. By Lemma \ref{lem: cyclic submodule}, the elements $u[(\eta^e)_e]$ and $[(\eta^e)_e]$ generate the same cyclic subgroup of $\check{H}$. This induces a surjective map 
    \begin{equation}\label{H/Aut and cyclic subgroups}
       H/\Aut(\A)\twoheadrightarrow\{\text{cyclic subgroups of $\check{H}$}\} \iso \{\text{cyclic subgroups of $H^1(\Gamma, \A_D)$}\}. 
    \end{equation}
    Again by Lemma \ref{lem: cyclic submodule}, if two elements $[\eta],[\tau]\in \check{H}$ generate the same cyclic subgroup, then there is $u\in \A^\times$ such that $[\tau]=u[\eta]$. Thus, the torsors $\mathcal{T},\mathcal{T'}\in H$ corresponding to $[\eta]$ and $[\tau]$, respectively, satisfy $\mathcal{T'}=u\cdot \mathcal{T}$. This proves injectivity of \eqref{H/Aut and cyclic subgroups}.
\end{proof}

If $\A_D$ is a codilation sheaf for $\Gamma$ and $\Delta\subseteq \Gamma$ is a subgraph, we denote by $\A_D\vert_\Delta$ the codilation sheaf obtained by restricting $\A_D$ to $\Delta$.  As explained in \cite[Section 3.3]{v1}, for $i=0,1$ there are surjective maps $C^i(\Gamma,\A_D)\to C^i(\Delta, \A_D\vert_\Delta)$ given by restricting $V(\Gamma)$ to $V(\Delta)$ in the case of $i=0$ and restricting $E(\Gamma)$ to $ E(\Delta)$ for $i=1$. This induces a surjection 
\[
\vert_\Delta: H^1(\Gamma,\A_D)\twoheadrightarrow H^1(\Delta,\A_D\vert_\Delta).
\]

If $\B\subseteq \A$ is a subgroup such that $D(v)\subseteq \B$ for every $v\in V(\Gamma)$ for an $\A$-dilation datum $D$, we denote by $\B_D$ the codilation sheaf obtained by viewing $D$ as a $\B$-dilation datum.

\begin{lemma} \label{lem: restriction of covers to subgraphs}
     Let $\Gamma$ be a connected graph, and suppose there exists a cycle $\Delta$ in $\Gamma$ with $H^1(\Delta,\A_D\vert_\Delta)\neq 0$. Then $H^1(\Delta,\A_D\vert_\Delta)\iso \ZZ/r\ZZ$ for some $r\geq 2$. Let $[\eta_0],[\eta_1]\in H^1(\Gamma,\A_D)$ be elements mapping respectively to a non-generator and a generator of $\ZZ/r\ZZ$ under the map 
    \[
    H^1(\Gamma,\A_D)\twoheadrightarrow H^1(\Delta,\A_D\vert_\Delta)\iso \ZZ/r\ZZ.
    \]
    For $i=0,1$, let $f_i: \Gamma_i\to \Gamma$ be the harmonic $\A$-cover of $\Gamma$ corresponding to $[\eta_i]$. Then $f_0^{-1}(\Delta)$ is disconnected and $f_1^{-1}(\Delta)$ is connected. In particular, the graphs $\Gamma_0$ and $\Gamma_1$ are non-isomorphic as cyclic dilated covers of $\Gamma$.    
    
\end{lemma}
\begin{proof}
    For $i=0,1$, let $(\eta_i^e)_{e\in E(\Gamma)}$ be a representative for $[\eta_i]\in H^1(\Gamma,\A_D)$. Then $(\eta_i^e)_{e\in E(\Delta)}$ is a representative for the image of $[\eta_i]$ under the restriction map $\vert_\Delta$. From the explicit correspondence in the proof of Theorem \ref{classification of harmonic covers} (\cite[Theorem 2.3]{Len-Ulrisch-Zakharov}) it follows that the $\A$-cover of $\Delta$ corresponding to $[\eta_i]\vert_\Delta$ is precisely the restriction map $f_i\vert_\Delta:f_i^{-1}(\Delta)\to \Delta$.

    Since $\Delta$ is a cycle, Proposition \ref{cohomology of a graph} implies that $H:=H^1(\Delta,\A_D\vert_\Delta)\iso \A/\sum_{v\in V(\Delta)}D(v)$, so $H\iso \ZZ/r\ZZ$ for some $r\geq 1$. Since $H$ is nontrivial by assumption, we have $r\geq 2$ and $\sum_{v\in V(\Delta)}D(v)\neq \A$. Let $\B\subsetneq \A$ be subgroup such that $D(v)\subseteq \B$ for every $v\in V(\Delta)$. Then $H^1(\Delta,\B_D\vert_\Delta)\iso \B/\sum_{v\in V(\Delta)}D(v)$, so the image of 
    \[
    i_\B: H^1(\Delta, \B_D\vert_\Delta) \hookrightarrow  H^1(\Delta,\A_D\vert_\Delta)=H
    \]
    is a proper subgroup of $H$. Since $[\eta_1]\vert_\Delta$ is a generator of $H$, it is not contained in the image of $i_\B$ for any such $\B$. Let $d$ denote the order of $[\eta_0]\vert_\Delta$ in $H$. By assumption, we have $d\mid r$ and $d<r$. Let $\B'$ be the unique subgroup of $H$ of order $d$. By the third isomorphism theorem, we have $\B'=\B/\sum_{v\in V(\Delta)}D(v)$ for some $\B\subsetneq \A$. Then $[\eta_0]\vert_\Delta$ is contained in the image of $i_\B$. By Proposition \ref{prop: connected covers}, the graph $\Delta_0=f_0^{-1}(\Delta)$ is disconnected and $\Delta_1=f_1^{-1}(\Delta)$ is connected. \\ 
\end{proof}

\begin{thm} \label{uniqueness vs dilated cohomology}
    Let $\Gamma$ be a connected graph and $\A\iso \ZZ/n\ZZ$. Then there is a unique connected cyclic cover of $\Gamma$ of degree $n$ with dilation datum $D$ if and only if the dilated cohomology group $H^1(\Gamma,\A_D)$ satisfies one of the following conditions.
    \begin{enumerate}
        \item $H^1(\Gamma,\A_D)=0$ and $\sum_{v\in V(\Gamma)}D(v)=\A$ (equivalently, $\operatorname{lcm}_{v\in V(\Gamma)}(|D(v)|)=n$), or \label{unique trivial}
        \item $H^1(\Gamma, \A_D)\iso \ZZ/m\ZZ$ for some $m\geq 2$ and $\operatorname{lcm}_{v\in V(\Gamma)}(\frac{n}{d},|D(v)|)\neq n$ for every $d\mid m$, $d\geq 2 $. \label{unique nontrivial}
    \end{enumerate}
    Moreover, if $\Gamma$ satisfies Condition \eqref{unique nontrivial}, then $\Gamma$ has exactly one cycle. 
\end{thm}
\begin{proof}
    We will use Theorem \ref{classification of harmonic covers} repeatedly. Let $\beta$ denote the first Betti number of $\Gamma$.
    
    \textbf{Case 1:} $H^1(\Gamma,\A_D)=0$.
    
    In this case there is a unique harmonic $\A$-cover with $\A$-dilation datum $D$, so there is a unique cyclic dilated cover $\Gamma'\to \Gamma$ of degree $n$ with dilation datum $D$. By Proposition \ref{prop: connected covers}, if $\B:=\sum_{v\in V(\Gamma)}D(v)=\A$, then $\Gamma'$ is connected. Since $H^1(\Gamma,\A_D)=0$, if $\B\subsetneq \A$, then Proposition \ref{prop: connected covers} implies that $\Gamma'$ is disconnected.

    \textbf{Case 2:} $H^1(\Gamma,\A_D)\iso \ZZ/m\ZZ$ for some $m\geq 2$.
    
    In particular,the graph $\Gamma$ is not a tree, so $\beta\geq 1$. For each prime $q\mid n$ and each $1\leq i\leq \beta$ let $\Delta_i^q$ be cycles in $\Gamma$ defined as in Proposition \ref{cohomology of a graph}. Then, the $q$-group
    \[
    H^1(\Gamma,\A_{q,D})\iso \bigoplus_{i=1}^\beta \A_q/\sum_{v\in V(\Delta^q_i)} D(v)_q
    \]
    is cyclic. For a fixed prime $q$, the cycles $\Delta_1^q,\dots,\Delta_\beta^q$ form a basis for $H_1(\Gamma)$. 
    Suppose that $\beta\geq 2$. Then, for every prime $q\mid n$ there is a cycle $\Delta^q_j$ such that 
    \[
    \sum_{v\in V(\Delta^q_j)} D(v)_q=\A_q.
    \]
    This implies that $\sum_{v\in V(\Gamma)}D(v)=\A$, or equivalently, $\operatorname{lcm}_{v\in V(\Gamma)}(|D(v)|)=n$. By Proposition \ref{prop: connected covers}, it follows that the harmonic $\A$-cover corresponding to any element of $H^1(\Gamma,\A_D)$ is connected. Let $[\eta_0],[\eta_1]\in H^1(\Gamma,\A_D)$ be elements corresponding to 0 and a generator of $\ZZ/m\ZZ$, respectively. Since $m\geq 2$, by Proposition \ref{cohomology of a graph} there is a prime $q\mid n$ and a cycle $\Delta$ in $\Gamma$ for which $\sum_{v\in V(\Delta)}D(v)_q \neq \A_q$, which implies that $H^1(\Delta,\A_D\vert_\Delta)\neq 0$. Then, under the restriction map 
    \[
    H^1(\Gamma,\A_D)\twoheadrightarrow H^1(\Delta,\A_D\vert_\Delta)\iso \ZZ/r\ZZ,
    \]
    the elements $[\eta_0]$ and $[\eta_1]$ map to $0$ and a generator of $\ZZ/r\ZZ$, respectively. 
    Lemma \ref{lem: restriction of covers to subgraphs} implies that $[\eta_0]$ and $[\eta_1]$ give rise to 2 different connected degree $n$ cyclic covers of $\Gamma$ with dilation datum $D$. 
    
    Now suppose that $\beta=1$. Let $\Delta$ denote the unique cycle in $\Gamma$  and let $z=\vert \sum_{v\in V(\Delta)} D(v)\vert $. Then we have $|H^1(\Gamma,\A_D)|=m=n/z$ by Proposition \ref{cohomology of a graph}. For any $d\mid m$, let 
    \[
    b_d:=\operatorname{lcm}_{v\in V(\Gamma)}(n/d,|D(v)|).
    \]
    Let $[\eta]\in H^1(\Gamma,\A_D)$ be an element of order $a$. Then there is $d_a\geq 1$ such that $m=d_aa$. Write $b=b_{d_a}$ and let $\B$ be the unique subgroup of $\A$ of order $b$. Since $a=\frac{m}{d_a}=\frac{n}{d_az}$, the group $\B\subseteq \A$ is the smallest subgroup such that $D(v)\subseteq \B$ for every $v\in V(\Gamma)$ and $a\mid |H^1(\Gamma,\B_D)|=b/z$. Since $H^1(\Gamma,\A_D)$ is cyclic, this implies
    \[
    [\eta]\in \operatorname{im}(H^1(\Gamma, \B_D)\hookrightarrow H^1(\Gamma,\A_D)).
    \]
    By Proposition \ref{prop: connected covers}, the element $[\eta]$ gives rise to a disconnected cover if and only if $\B\subsetneq \A$, i.e. $b\neq n$. If $a=m$, then $b=n$, so $[\eta]$ gives rise to a connected harmonic $\A$-cover. Since $H^1(\Gamma,\A_D)$ is cyclic, any two elements of the same order generate the same cyclic subgroup, so by Proposition \ref{cyclic subgroups give isomorphic covers} they give rise to the same cyclic dilated cover. Therefore, in this case there is a unique connected cyclic cover of $\Gamma$ of degree $n$ with dilation datum $D$ if and only if for each $d\mid m$ with $d\geq 2$, we have $b_d\neq n$. 

    \textbf{Case 3:} $H^1(\Gamma,\A_D)$ is not cyclic.

    Let $q$ be a prime such that the $q$-Sylow subgroup of $H^1(\Gamma,\A_D)$ is not cyclic. Let $T$ be the spanning tree $T^q$ defined in Proposition \ref{cohomology of a graph}, and let $E(\Gamma)\backslash E(T)=\{e_1,\dots,e_\beta\}$. 
    By \cite[Remark 4.3]{v1}, any element $[\eta]\in H^1(\Gamma, \A_D)$ can be represented by a cochain $\eta=(\eta^e)_{e}\in C^1(\Gamma, \A_D)=\bigoplus_{e\in E(\Gamma)} \A/C(e)$ such that $\eta^e=0$ for all $e\neq e_i$. This gives rise to a surjective map 
    \[
    \bigoplus_{i=1}^\beta \A/C(e_i)\twoheadrightarrow H^1(\Gamma,\A_D),
    \]
    where an element $(a_i)_i\in \bigoplus_{i=1}^\beta \A/C(e_i)$ is mapped to the class $[\eta]\in H^1(\Gamma,\A_D)$ given by $\eta^e=0$ for all $e\in E(T)$ and $\eta^{e_i}=a_i$ for each $i$. Since the $q$-Sylow subgroup of $H^1(\Gamma,\A_D)$ is not cyclic and each quotient $\A/C(e_i)$ is cyclic, we may assume (after possibly rearranging the $e_i$) that the the $q$-Sylow subgroups $\A_q/C(e_1)_q$ and $\A_q/C(e_2)_q$ of $\A/C(e_1)$ and $\A/C(e_2)$ are not trivial. Assume moreover that
    \begin{equation}\label{eqn: inequality of sylow orders}
      |\A_q/C(e_1)_q|\leq |\A_q/C(e_2)_q|.
    \end{equation}
    Let $\Delta$ denote the unique cycle in the subgraph $T\cup \{e_1\}$ of $\Gamma$. Then, the subgraph obtained by intersecting $T\cap \Delta$ is a spanning tree for $\Delta$, and $E(\Delta)\backslash E(T\cap \Delta)=\{e_1\}$, so similarly as for $\Gamma$, we have a surjective map
    \[
    \A/C(e_1)\twoheadrightarrow H^1(\Delta,\A_D\vert_\Delta).
    \]
    We note that, by definition of $\Delta$ and Proposition \ref{cohomology of a graph}, this map is an isomorphism when restricting to the $q$-Sylow subgroups, i.e.  
    \begin{equation}\label{eqn: iso of q-Sylows}
            \A_q/C(e_1)_q\xrightarrow{\sim} H^1(\Delta,\A_{q,D}\vert_\Delta).
    \end{equation}
    Since the restriction map $\vert_\Delta: H^1(\Gamma,\A_D)\twoheadrightarrow H^1(\Delta,\A_D\vert_\Delta)$ is given by sending an element $[\eta]\in H^1(\Gamma,\A_D)$ with a representative $\eta=(\eta^e)_{e\in E(\Gamma)}$ to the class with representative $(\eta^e)_{e\in E(\Delta)}$, we have a commutative diagram
    \[\begin{tikzcd}
	{\bigoplus_{i=1}^\beta \A/C(e_i)} & {H^1(\Gamma,\A_D)} \\
	{\A/C(e_1)} & {H^1(\Delta,\A_D\vert_\Delta)}
	\arrow[two heads, from=1-1, to=1-2]
	\arrow[two heads, from=1-1, to=2-1]
	\arrow["{\vert_\Delta}", two heads, from=1-2, to=2-2]
	\arrow[two heads, from=2-1, to=2-2]
\end{tikzcd}\]
where the vertical map on the left is the projection. Let $d=|\A/C(e_1)|$, and let $a\in \A/C(e_i)$ be an element of order $d/q^{v_q(d)}<d$, where $v_q(d)$ is the $q$-adic valuation of $d$. Consider the classes $[\eta_0], [\eta_1] \in H^1(\Gamma,\A_D)$ given by the image of 
\[
(a,1,\dots,1),(1,1,\dots,1)\in \bigoplus_{i=1}^\beta\A/C(e_i),
\]
respectively. Then, up to reordering the edges in $E(\Delta)$, the classes $[\eta_0]\vert_\Delta,[\eta_1]\vert_\Delta\in H^1(\Delta,\A_D\vert_\Delta)$ are represented by cochains 
\[
(a,0,\dots,0),(1,0,\dots,0)\in C^1(\Delta,\A_D\vert_\Delta),
\]
respectively. By Equation \eqref{eqn: iso of q-Sylows} and since $q$ does not divide the order of $a$ but $q\mid |H^1(\Delta,\A_D\vert_\Delta)|$, we have that $[\eta_1]\vert_\Delta$ generates $H^1(\Delta,\A_D\vert_\Delta)$ while $[\eta_0]\vert_\Delta$ does not.
    Then, Lemma \ref{lem: restriction of covers to subgraphs} implies that the harmonic $\A$-covers of $\Gamma$ corresponding to $[\eta_0]$ and $[\eta_{1}]$ are not isomorphic as covers, so they give rise to distinct cyclic covers of $\Gamma$ of degree $n$ with dilation datum $D$. 

    It remains to show that the $\A$-covers corresponding to $[\eta_0]$ and $[\eta_1]$ are connected. First, note that Equation \eqref{eqn: inequality of sylow orders} implies that $[\eta_0]$ and $[\eta_1]$ are of maximal order.
    For any proper subgroup $\B\subsetneq \A$ such that $D(v)\subseteq \B$ for all $v\in V(\Gamma)$, there is a prime $\ell\mid n$ such that $\B_\ell\subsetneq \A_\ell$. For $1\leq i\leq \beta$, let $\Delta_i^\ell$ be as in Proposition \ref{cohomology of a graph} and let us write
    \[
    H^1(\Gamma,\A_{\ell,D})\iso \bigoplus_i \A_\ell/\sum_{v\in V(\Delta^\ell_i)} D(v)_\ell \iso \bigoplus_i \ZZ/\ell^{c_i}\ZZ,
    \]
    and 
    \[
     H^1(\Gamma,\B_{\ell,D})\iso \bigoplus_i \B_\ell/\sum_{v\in V(\Delta^\ell_i)} D(v)_\ell \iso \bigoplus_i \ZZ/\ell^{b_i}\ZZ,
    \]
    where $c_i,b_i\geq 0$. Since $\B_\ell\subsetneq \A_\ell$, we have $b_i<c_i$ for all $1\leq i\leq \beta$. This implies that the order of any element in $H^1(\Gamma, \B_D)$ is strictly less than the maximal order in $H^1(\Gamma,\A_D)$. Thus, the elements $[\eta_0],[\eta_1]$ cannot be in the image of the map $H^1(\Gamma, \B_D)\hookrightarrow H^1(\Gamma,\A_D)$ for any such $\B$. It follows from Proposition \ref{prop: connected covers} that $[\eta_0],[\eta_1]$ give rise to connected covers. \\
\end{proof}

\begin{remark}\label{rem: symmetry}
    By looking at the explicit correspondence in the proof of \cite[Theorem 2.3]{Len-Ulrisch-Zakharov} (Theorem \ref{classification of harmonic covers}), one can check that elements the $[\eta_0]$ and $[\eta_1]$ in Cases 2 and 3 of the proof give rise to harmonic covers $\Gamma'\to \Gamma$, and $\Gamma''\to \Gamma$ such that $\Gamma'\not\iso\Gamma''$ as graphs. 
\end{remark}

In the above proof we have also shown the following result.

\begin{prop} \label{non-existence of a connected cover}
    Let $\Gamma$ be a connected graph and $\A=\ZZ/n\ZZ$ for $n\geq 2$. There is no connected cyclic cover of $\Gamma$ of degree $n$ with dilation datum $D$ if and only if $H^1(\Gamma,\A_D)=0$ and $\sum_{v\in V(\Gamma)}D(v)\neq\A$.
\end{prop}

For a subgraph $\Delta$ of a graph $\Gamma$, we denote by $\Gamma\backslash\Delta$ the topological subspace obtained by viewing $\Gamma$ and $\Delta$ as topological spaces. Then, if $T=\Gamma \text{ or } \Gamma\backslash \Delta$, we write $H_1(T)$ for the simplicial homology of $T$ with integer coefficients. For a prime $q$, let $v_q(n)$ denote the $q$-adic valuation of $n$.

\begin{thm} \label{unique cyclic graph cover}
    Let $\A$ and $\A_q$ be as above. Let $\Gamma$ be a connected graph and $D$ an $\A$-dilation datum on $\Gamma$. For each prime $q\mid n$, we define a (possibly empty or disconnected) subgraph $B_q$ of $\Gamma$ by 
    \[
    V(B_q) = \{v\in V(\Gamma) : q^{v_q(n)}\mid |D(v)|\}=\{v\in V(\Gamma) : D(v)_q= \A_q\},
    \]
    and $E(B_q)=\emptyset$. Then there is a unique connected degree $n$ cyclic graph cover of $\Gamma$ with $\A$-dilation datum $D$ if and only if for every prime $q\mid n$ one of the following holds. 
    \begin{enumerate}[(i)]
         \item  $B_q\neq \emptyset$ and $H_1(\Gamma\backslash B_q)=0$, or \label{trivial homology}
         \item $B_q=\emptyset$ and $H_1(\Gamma)\iso \ZZ$. \label{cycle case}
     \end{enumerate}
\end{thm}
\begin{proof}

It suffices to prove that one of \eqref{unique trivial} or \eqref{unique nontrivial} in Theorem \ref{uniqueness vs dilated cohomology} holds if and only if, for each prime $q\mid n$, one of \eqref{trivial homology} or \eqref{cycle case} holds.
Consider the set of primes
 \[
 S:=\{q\mid n: q \text{ is a prime for which  \eqref{trivial homology} does not hold}\}.
 \]
For any prime $q\mid n$ such that $q\notin S$, Corollary \ref{superdilated} implies that $H^1(\Gamma, \A_{q,D})=0$. Then the decomposition of cohomology from Equation \eqref{decomposition of cohomology} can be taken over $q\in S$, i.e.
\begin{equation} \label{decomposition of cohomology less primes}
    H^1(\Gamma,\A_D)\iso \bigoplus_{q\in S} H^1(\Gamma, \A_{q,D}).
\end{equation}
If $S=\emptyset$, then $H^1(\Gamma, \A_D)=0$ and $B_q\neq \emptyset$ for all primes $q\mid n$, so $\sum_{v\in V(\Gamma)}D(v)=\A$. Conversely, suppose that $H^1(\Gamma, \A_D)=0$ and $\sum_{v\in V(\Gamma)}D(v)=\A$. Then $H^1(\Gamma,\A_{q,D})=0$ and $B_q\neq \emptyset$ for every prime $q\mid n$. If $\Gamma$ is not a tree, Corollary \ref{superdilated} implies that, for every prime $q\mid n$, every cycle in $\Gamma$ contains a vertex $v$ with $D(v)_q=\A_q$, so $H_1(\Gamma\backslash B_q)=0$. If $\Gamma$ is a tree, then $H_1(\Gamma)=0$, so $H_1(\Gamma\backslash B_q)=0$. Thus, the condition that every prime $q\mid n$ satisfies \eqref{trivial homology} is equivalent to Condition \eqref{unique trivial} in Theorem \ref{uniqueness vs dilated cohomology}. 

Next, we show that the set 
\[
S':=\{q\mid n: q \text{ is a prime for which \eqref{cycle case} holds} \}\subseteq S
\]
is nonempty and equal to $S$ if and only if Condition \eqref{unique nontrivial} in Theorem \ref{uniqueness vs dilated cohomology} holds. 
If $S'=S\neq \emptyset$, then $H_1(\Gamma)\iso\ZZ$, so $H^1(\Gamma,\A_D)\iso \ZZ/m\ZZ$ for some $m\geq 1$ by Proposition \ref{cohomology of a graph}. By the same proposition and Equation \eqref{decomposition of cohomology less primes}, a prime $q\mid m$ if and only if $q\in S'=S$. Since $S\neq \emptyset$, we have $m\geq 2$. Moreover, for every prime $q\mid m$, we have $B_q=\emptyset$, so $q^{v_q(n)}\nmid |D(v)|$ for every $v\in V(\Gamma)$. Therefore $\operatorname{lcm}_{v\in V(\Gamma)}(\frac{n}{q},|D(v)|)\neq n$, which implies that Condition \eqref{unique nontrivial} in Theorem \ref{uniqueness vs dilated cohomology} holds. Conversely, suppose Condition \eqref{unique nontrivial} in Theorem \ref{uniqueness vs dilated cohomology} holds. Then $S\neq \emptyset$ by the decomposition in \eqref{decomposition of cohomology less primes}, and $\Gamma$ has exactly one cycle $\Delta$, i.e. $H_1(\Gamma)\iso\ZZ$. Let $q\in S\backslash S'$, so $B_q\neq \emptyset$. Since $q$ does not satisfy \eqref{trivial homology}, we have $H_1(\Gamma\backslash B_q)\neq 0$, so $H_1(\Gamma\backslash B_q)\iso \ZZ$. This implies that $\Delta\subseteq \Gamma\backslash B_q$, so $\sum_{v\in V(\Delta)}D(v)_q\neq \A_q$. By Proposition \ref{cohomology of a graph} this implies that $q\mid m$. Then $\operatorname{lcm}_{v\in V(\Gamma)}(\frac{n}{q},|D(v)|)\neq n$, so we have $q^{v_q(n)}\nmid |D(v)|$ for every $v\in V(\Gamma)$, which contradicts that $B_q\neq \emptyset$. Therefore $S=S'$. \\
\end{proof}

\subsection{Reduction type covers} \label{sec: defn of reduction type covers}
\mbox{}
\vspace{5pt}

 In this section, we relate cyclic covers of graphs to \emph{reduction type covers}, which we will use to prove our main results. Let $\C/O_K$ be the minimal regular model with strict normal crossings (msnc) of $C$. Let $\Gamma$ be the dual graph of the special fibre $\C_k$. We recall that $\C_n$ denotes the normalisation of $\C\times _{O_K}O_{K(n)}$ and 
\[
\varphi: \C_n\to \C
\]
the canonical morphism. We denote by 
\[
\rho: \C(n)\to \C_n
\]
the minimal desingularisation of $\C_n$, and by $\C(n)_{\operatorname{min}}$ the msnc model of $C(n)=C\times_K K(n)$ over $O_{K(n)}$. We denote by $\varphi_k$ and $\rho_k$, respectively, the induced maps on special fibres. In the proof of Corollary \ref{unique reduction type}, we will show that $\C(n)_k \to \C_k$ gives rise to a \emph{reduction type cover}.  

\subsubsection{A dual graph $\Gamma_n$ for $(\C_n)_k$} \label{subsec: dual graph}
\mbox{}
\vspace{5pt}

Let us define a \emph{dual graph} $\Gamma_n$ for $(\C_n)_k$. Since $\C_n$ is not necessarily regular, we cannot use intersection theory to define $\Gamma_n$, but we can construct it explicitly as follows. Let the vertices in $\Gamma_n$ correspond to irreducible components in $(\C_n)_k$ and let each point $x$ in the intersection of two components $v\neq u$ in $(\C_n)_k$ give rise to an edge between the corresponding vertices. By Proposition \ref{Halle-Nicaise 1.3.2} and Proposition \ref{number of components divides}, for each component $v$ in $\C_k$ there are $\ell_v$ disjoint components in $(\C_n)_k$ lying above it under the map $\varphi_k: (\C_n)_k\to \C_k$, where $\ell_v \mid \gcd (\weight(v),n)$. Moreover, if $v$ has $g(v)=0$, then $\ell_v=\gcd (\weight(v),n)$ by Proposition \ref{number of components g=0}. For each intersection of components $u,v\in V(\Gamma)$, corresponding to an edge $e\in E(\Gamma)$, there are $\gcd(\weight(e),n)$ points in $(\C_n)_k$ lying above it by Fact \ref{number of points above intersection point}. These points lie in $\varphi_k^{-1}(v)\cap\varphi_k^{-1}(u)$, so each of them is an intersection point of two distinct irreducible components lying above $v$ and $u$, respectively.

 We note that, if $\C_n$ is already a regular model with normal crossings, then this recovers the usual dual graph defined via intersection theory. 
 
\subsubsection{The map $\Gamma_n\to \Gamma$ is a cyclic dilated cover}
\mbox{}
\vspace{5pt}

Let $v\in V(\Gamma)$, and let $e\in E(\Gamma)$ be an edge incident to $v$. Let $v_1,\dots,v_{\ell_v}\in V(\Gamma_n)$ be the vertices lying above $v$ under the map $\varphi_k$. By Lemma \ref{lem: Galois action on components}, the group $\Gal(K(n)/K)$ acts transitively on the set $\{v_1,\dots,v_{\ell_v}\}$. Then the number of edges in $\Gamma_n$ incident to each $v_{i}$ and lying above $e$ is
\[
\frac{\operatorname{gcd}(\operatorname{weight}(e),n)}{\ell_v}.
\] 
We note that this number is indeed an integer, since $ \operatorname{weight}(v) \mid \operatorname{weight}(e)$ by definition. 
Thus, the morphism of graphs $f:\Gamma_n\to \Gamma$ induced by $\varphi_k$ is a harmonic morphism of degree $n$ with dilation assignment 
\[
d_f(v')=\frac{n}{\ell_v} \ \text{ and } \ d_f(e')=\frac{n}{\gcd(\weight(e),n)},
\]
for every $v'\in V(\Gamma_n)$ and $e'\in E(\Gamma_n)$ lying above any $v\in V(\Gamma)$ and any $e\in E(\Gamma)$, respectively. Moreover, the group $\Gal(K(n)/K)\iso \ZZ/n\ZZ$ preserves the fibres of $\Gamma_n\to \Gamma$ and acts transitively on them. Moreover, since $C$ is assumed to be geometrically connected, the graph $\Gamma_n$ is connected (see e.g. \cite[Corollary 8.3.6(b)]{Liu}). Therefore, we obtain the following result.

\begin{lemma}\label{lem: cyclic cover of dual graphs}
    The graph morphism $\Gamma_n\to \Gamma$ induced by $(\C_n)_k\to \C_k$ is a connected cyclic dilated cover of degree $n$.
\end{lemma}

\subsubsection{Definition of reduction type cover}
\mbox{}
\vspace{5pt}

The lemma above, together with Proposition \ref{desingularisation}, motivates the following definitions. First, we recall that a \emph{graph subdivision} of $\Gamma$ is a graph $\Gamma'$ resulting from replacing some edges in $\Gamma$ by chains of vertices of degree 2, i.e. vertices attached exactly to 2 edges. The induced map $\Gamma'\to \Gamma$, which we call the \emph{subdivision map}, is a \emph{non-finite morphism} of graphs, i.e. a more general notion of graph morphism in which edges are allowed to contract to vertices.

\begin{defn} \label{resolution graph}
    Let $\Gamma$ be a dual graph and $f:\Gamma_n \to \Gamma$ a connected cyclic cover of graphs of degree $n$. The \emph{resolution} of $f$ is a subdivision map $\Gamma(n)\to \Gamma_n$ as follows.
    For every edge $e'$ connecting to 2 vertices $v'_1,v'_2\in \Gamma_n$ lying above vertices $v_1,v_2\in \Gamma$ of multiplicities $d_1,d_2$, respectively, compute: 
    \[
    n'=\frac{n\operatorname{gcd}(d_1,d_2,n)}{\operatorname{gcd}(d_1,n)\operatorname{gcd}(d_2,n)}.
    \]
    If $n'=1$, set $\lambda =0$. Otherwise, let $r$ be the unique solution to the equation $rd_1'+d_2'=0 \ \operatorname{mod} n'$ such that $0<r<n'$. Let $[b_1,\dots,b_\lambda]_{\operatorname{HJ}}$ be the Hirzebruch-Jung continued fraction expansion of $n'/r$. Let $(\mu_1,\dots,\mu_\lambda)$ be the multiplicity vector associated to the tuple $(n',r,d_1',d_2')$. 
    Then replace the edge $e'$ by $\lambda$ vertices of degree 2 with arithmetic and geometric genera $p=g=0$ and attach to them the multiplicities $\{\mu_1,\dots, \mu_\lambda\}$ in order from $v'_1$ to $v'_2$. We call $\Gamma(n)$ the \emph{resolution graph} of $\Gamma_n\to \Gamma$. To every vertex $v'\in \Gamma(n)$ that is also in $\Gamma_n$, we attach its multiplicity and genus given as follows. Let $v'\in \Gamma_n$ be a vertex lying above a vertex $v\in \Gamma$, and let $m_v$ and $g(v)$ denote its associated multiplicity and geometric genus.  
    \begin{itemize}
        \item The multiplicity of $v'$ is 
        \(
        m_{v'}=m_v/\gcd(m_v,n).
        \)
        \item The geometric genus of $v'$ is given by $g(u_j)$ in Equation \eqref{eqn: genus formula} (Genus formula).
    \end{itemize}
\end{defn}

\begin{defn} \label{defn: compatible dilation}
    Let $R$ be an rnc reduction type with dual graph $\Gamma$ and let $n\geq 1$. Let $\A=\ZZ/n\ZZ$. An $\A$-dilation datum $D$ is said to be \emph{compatible} with $R$ is it satisfies the following conditions.\begin{itemize}
        \item For any $x\in E(\Gamma)\cup \{v\in V(\Gamma): g(v)=0\}$, the group $D(x)\subseteq \A$ is the unique subgroup of order 
        \[
        \frac{n}{\gcd(\weight(x),n)}.
        \]
        \item For any $v\in V(\Gamma)$ with $g(v)\geq 1$, the group $D(x)\subseteq \A$ has order $n/d$, where 
        \[
        d \ \vert \ {\gcd(\weight(v),n)}.
        \]
    \end{itemize}
\end{defn}

\begin{defn} \label{defn: reduction type cover}
    Let $R$ be an rnc reduction type with dual graph $\Gamma$. A \emph{reduction type cover} of $R$ \emph{of degree} $n$ is an rnc reduction type $R'$ together with a map $\Gamma''\to \Gamma$ obtained as follows. 
    \begin{itemize}
        \item  The dual graph of $R'$ is the resolution graph $\Gamma''$ of a connected cyclic graph cover $\Gamma'\to \Gamma$ of degree $n$ with dilation datum $D$ compatible with $R$.
        \item The map $\Gamma''\to \Gamma$ is the morphism obtained as the composition of $\Gamma'\to\Gamma$ with its resolution.
        \item  The multiplicities $m_v$ and geometric genera $g(v)$ of the components of $R'$ are induced by those of the corresponding vertices in the resolution graph. The arithmetic genus $p(v)$ is given by the sum of $g(v)$ and the number of loops attached to $v$ in the dual graph.
    \end{itemize}
\end{defn}

Suppose that $R$ is a mrnc reduction type. Then, a reduction type cover $R'$ of $R$ need not be mrnc. 
Similarly, if $R$ is an msnc reduction type, then $R'$ need not be msnc. If $R'$ is snc, then we denote by $R'_{\text{min}}$ the msnc reduction type obtained by contracting every component $v$ in $R'$ that satisfies 
\begin{equation} \label{eqn: contractible}
  g(v)=0, \ \ v\cdot v=-1 \ \text{ and } \ |I(v)|<3, \ I(v)\neq \{v',v'\} \text{ for some } v'\neq v.  
\end{equation}

\subsection{Results on reduction type covers}\label{sec: results on reduction type covers}
\mbox{}
\vspace{5pt}

In this section we prove the main results of the paper. We will use  results from Section \ref{sec: results on cyclic covers}, where it is assumed that $\Gamma$ has no loops or multi-edges. Let us first explain why this assumption does not affect our results. 

Given a reduction type $R$ with dual graph $\Gamma$, we can always obtain an rnc reduction type $R'$ via a minimal number of blow ups such that its dual graph $\Gamma'$ has no loops or multi-edges. Suppose that, in order to obtain $R'$ from $R$, we blow up a point $e$ in the intersection of two components $v_1,v_2$ in $R$ with multiplicities $m_{v_1},m_{v_2}$. Then we obtain a new component $u$ intersecting $v_1$ and $v_2$ transversally, with $g(u)=0$ and multiplicity $m_u=m_{v_1}+m_{v_2}$. Note that the weight of $v_1,v_2$ or any other components in $R$ remains the same, and $\weight(u)=\weight(e)$, where we identify the point $e$ with the corresponding edge in $\Gamma$. Lastly, there is a correspondence between mrnc reduction types $R$ and minimal reduction types $R'$ as above, by \cite[Theorem 4.6]{Jakab}. Hence, in the proofs of this section we may assume that $\Gamma$ has no loops or multi-edges, but the statements hold in full generality. 

\begin{thm} \label{thm: reduction type cover}
    Let $n\geq 1$. A reduction type $R$ has a unique reduction type cover of degree $n$ if and only if for each prime $q\mid n$, one of the following holds: 
    \begin{enumerate}[(a)]
        \item $q\nmid \operatorname{weight}(R)$, and any component $v$ with geometric genus $g(v)\geq1$ and any cycle $Z$ in $R$ satisfy
        \[
        q\nmid \operatorname{weight}(v)  \text{ and }   q\nmid \operatorname{weight}(Z).\]
        
        \item  \(q\mid \operatorname{weight}(R)\), and $R$ has no cycles and exactly one component with geometric genus $g\geq 1$.
        \item \( q\mid \operatorname{weight}(R)\), and $R$ has exactly one cycle and no components with geometric genus $g\geq 1$.
    \end{enumerate}
\end{thm}
\begin{proof}
Let $\A=\ZZ/n\ZZ$, and let $\Gamma$ be the dual graph of a reduction type $R$. Since the resolution of a given cyclic dilated cover $\Gamma'\to \Gamma$ is unique, there is a unique reduction type cover of degree $n$ if and only if there is a unique $\A$-dilation datum $D$ compatible with $R$ and a unique connected cyclic cover of $\Gamma$ of degree $n$ with dilation datum $D$. 

Assume there is at least one $\A$-dilation datum $D$ giving rise to a \emph{connected} cyclic dilated cover of $\Gamma$. We claim that there is a unique such dilation datum compatible with $R$ if and only if, for every prime $q\mid n$, either every component $v$ in $R$ with $g(v)\geq 1$ has $q\nmid \weight(v)$, or \eqref{b} holds. Let $q\mid n$ be a prime. Note that it suffices to show the claim for $\A_q$-dilation data $D_q$. If for any $v\in V(\Gamma)$ such that $g(v)\geq 1$ we have $q\nmid \weight(v)$, then 
\[
\gcd(\weight(v),q^{v_q(n)})=1,
\]
so for any $x\in X(\Gamma)$, the group $D(x)_q\subseteq \A_q$ is the unique subgroup of order 
\[
\frac{q^{v_q(n)}}{\gcd(\weight(x),q^{v_q(n)})}.
\]
Thus, by Definition \ref{defn: compatible dilation}, the datum $D_q$ is uniquely determined by $R$.
Now suppose that $q$ satisfies \eqref{b}. Then $\Gamma$ is a tree, so $H^1(\Gamma,\A_D)=0$. Let $u\in V(\Gamma)$ denote the unique component in $R$ such that $g(u)\geq 1$. For any $v\in V(\Gamma)\backslash\{u\}$ we have $g(v)=0$ and $q\mid \weight(v)$, so $D(v)_q\subsetneq \A_q$. Suppose that $D(u)_q\subsetneq \A_q$. Then $\sum_{v\in V(\Gamma)}D(v)\neq \A$, so Proposition \ref{non-existence of a connected cover} implies that there is no connected cyclic cover of $\Gamma$ of degree $n$ with dilation datum $D$. 
Therefore $D(u)_q=\A_q$, so $D_q$ is uniquely determined by $R$. 

For the converse implication of the claim, suppose $u_1,\dots,u_r$ are the components in $R$ with $g(u_i)\geq 1$ and $q\mid \weight(u_i)$. Assume furthermore that $r\geq 1$ and \eqref{b} does not hold. If $r\geq 2$, then there are at least two dilation data compatible with $R$ given e.g. by $D(u_1)_q=\A_q$ and $D(u_2)_q=0$ or $\A_q$ by Proposition \ref{non-existence of a connected cover}. Thus, assume $r=1$. If $q\nmid \weight(R)$, then there is $v\in V(\Gamma)$ such that $q\nmid \weight(v)$, so $D(v)_q=\A_q$. Then by Proposition \ref{non-existence of a connected cover} there are at least two dilation data compatible with $R$ given by setting $D(u_1)_q=0$ or $\A_q$. Lastly, assume $q\mid \weight(R)$. Setting $D(u_1)_q=0$ gives $\sum_{v\in V(\Gamma)} D(v)_q\neq \A_q$. Then by Corollary \ref{superdilated} and since trees have trivial dilated cohomology, we have $H_1(\Gamma)=0$ if and only if $H^1(\Gamma,\A_{q,D})=0$. 
Thus, if $\Gamma$ has at least one cycle, i.e. $H_1(\Gamma)\neq 0$, setting either $D(u_1)_q=0$ or $\A_q$ again gives rise to two dilation data compatible with $R$ by Proposition \ref{non-existence of a connected cover}. This proves the claim.

Suppose that each prime $q\mid n$ satisfies one of \eqref{a}, \eqref{b} or \eqref{c}. Then the claim implies that any reduction type cover of $R$ of degree $n$ has the same dilation datum $D$. For any prime $q\mid n$, consider the subgraph $B_q$ of $\Gamma$ such that
\[
V(B_q) = \{v\in V(\Gamma) : q^{v_q(n)}\mid |D(v)|\}=\{v\in V(\Gamma) : D(v)_q= \A_q\},
\]
and $E(B_q)=\emptyset$. Then, if a prime $q\mid n$ satisfies either \eqref{a} or \eqref{b}, then $B_q\neq \emptyset$ and $H_1(\Gamma\backslash B_q)=0$, so \eqref{trivial homology} in Theorem \ref{unique cyclic graph cover} holds. 
If a prime $q\mid n$ satisfies \eqref{c}, then $H_1(\Gamma)\iso\ZZ$, all components in $R$ have geometric genus 0 and $q\mid \weight(v)$ for every $v\in V(\Gamma)$, so we have $B_q=\emptyset$, and \eqref{cycle case} in Theorem \ref{unique cyclic graph cover} holds. Therefore, Theorem \ref{unique cyclic graph cover} implies that there is a unique connected cyclic cover of $\Gamma$ of degree $n$ with dilation datum $D$, so there is a unique reduction type cover of degree $n$. 

For the converse, suppose that $R$ has a unique reduction type cover of degree $n$. Then there is a unique $\A$-dilation datum $D$ compatible with $R$ and a unique connected cyclic cover of $\Gamma$ of degree $n$ with dilation datum $D$. Let $q\mid n$ be a prime. Then the claim above implies that either \eqref{b} holds, or that any component $u$ in $R$ with $g(u)\geq 1$ has $q\nmid \weight(u)$, so $D(u)_q=\A_q$. We assume from now on that the latter holds. Since there is a unique connected cyclic cover of degree $n$ with dilation datum $D$, either \eqref{trivial homology} or \eqref{cycle case} in Theorem \ref{unique cyclic graph cover} hold. 
If \eqref{trivial homology} holds, then $B_q\neq \emptyset$, so there is $v\in V(\Gamma)$ such that $q^{v_q(n)}\mid |D(v)|$, so $q\nmid \weight(v)$. This implies that $q\nmid \weight(R)$. Furthermore, $H_1(\Gamma\backslash B_q)=0$ implies that every cycle $\Delta$ in $R$ contains a vertex $v_\Delta\in V(B_q)$, so $q\nmid \weight(v_\Delta)$. Thus, we have $q\nmid \weight(\Delta)$, and so \eqref{a} holds. 
If \eqref{cycle case} holds, then $B_q=\emptyset$, so by the proof of the claim, every $v\in V(\Gamma)$ has $g(v)=0$ and $q^{v_q(n)}\nmid |D(v)|$, so 
$q\mid \weight(R)$. Moreover, $H_1(\Gamma)\iso\ZZ$ implies that there is exactly one cycle in $R$, so \eqref{c} holds. 

\end{proof}

In the proof above, we have also shown the following result.

\begin{cor} \label{one comp lying above}
    Let $n\geq 1$. Suppose that, for each prime $q\mid n$, a reduction type $R$ satisfies either \eqref{a}, \eqref{b}, or \eqref{c} of Theorem \ref{thm: reduction type cover}. Let $\Gamma$ be the dual graph of $R$. Then the unique connected cyclic cover of $\Gamma$ of degree $n$ with dilation datum compatible with $R$ has exactly one component lying above each $v$ with $g(v)\geq 1$. 
\end{cor}

\begin{cor} \label{cor: reduction type}
        Let $p\nmid n$, and let $R$ be a reduction type. Suppose that each prime $q\mid n$ satisfies either \eqref{a}, \eqref{b}, or \eqref{c} of Theorem \ref{thm: reduction type cover}. 
    Then for any curve $C$ with reduction type $R$, the reduction type of $C(n)$ is uniquely determined, and it depends only on $R$ and $n$, and not on $C$ or $K$. 
    In particular, the result holds for curves with abelian and toric ranks zero. 
\end{cor}
\begin{proof}
Since there is a one-to-one correspondence between mrnc reduction types and msnc reduction types \cite[Theorem 4.6]{Jakab}, it suffices to prove the result for msnc reduction types. 
Let $C/K$ be a curve with msnc reduction type $R$, and let $\C$ be the msnc model of $C$. Let $R'$ be the msnc reduction type of the base change $C(n)/K(n)$. 
By Proposition \ref{Halle-Nicaise 1.3.2}, the model $\C(n)$ of $C(n)$ is also snc. Thus, the special fibre $\C(n)_k$ gives rise to an snc reduction type $R''$ of $C(n)$. Then the msnc reduction type $R'$ can be obtained by contracting those components $v$ satisfying Condition \eqref{eqn: contractible}. If we show that $R''$ is a reduction type cover of $R$ of degree $n$, then $R'=R''_\text{min}$ and the result follows by Theorem \ref{thm: reduction type cover}.

Let $\Gamma$, $\Gamma_n$ and $\Gamma(n)$ denote the dual graphs of $R$, $(\C_n)_k$ and $R''$, respectively. By Lemma \ref{lem: cyclic cover of dual graphs} the map $f:\Gamma_n \to \Gamma$ induced by $\varphi_k$ is a connected cyclic cover of degree $n$. 
Moreover, the graph $\Gamma(n)$ is the resolution graph of $f:\Gamma_n\to \Gamma$ (in the sense of Definition \ref{resolution graph}) by Proposition \ref{desingularisation}. Let $u$ be a component in $R''$. Since $R''$ is snc, the arithmetic and geometric genera of $u$ coincide, i.e. $p(u)=g(u)$, and they are given by the genus formula (Proposition \ref{genus formula}). Furthermore, by Propositions \ref{Halle-Nicaise 1.3.2} and \ref{desingularisation}, the multiplicity of $u$ in $R''$ is as in Definition \ref{resolution graph}. Therefore $R''$ is a reduction type cover of $R$ of degree $n$.

Lastly, if a curve has abelian and toric ranks zero, then by \cite[Corollary 7.4]{Liu-Lorenzini-Raynaud} (Corollary \ref{cor: Liu-Lorenzini-Raynaud}) we have $q\nmid \weight(R)$. Thus, Condition \eqref{a} in Theorem \ref{thm: reduction type cover} holds trivially for any prime $q\mid n$. \\ 
\end{proof}

\begin{remark}
    See Section \ref{algorithm} for how to compute explicitly the reduction type after tame base change.
\end{remark}

\subsection{Other consequences} \label{sec: other consequences}
\mbox{}
\vspace{5pt}

This framework also allows us to give a new proof of the following result by Liu, Lorenzini and Raynaud \cite[Corollary 7.4]{Liu-Lorenzini-Raynaud}. 

\begin{cor} \label{cor: Liu-Lorenzini-Raynaud}
Suppose a curve $C/K$ has abelian and toric ranks zero. Let $d_1,\dots,d_r$ be the multiplicities of the components in the reduction type $R$ of $C$. If $d:=\gcd(d_1,\dots,d_r) =\weight(R)> 1$, then $p>0$ and $d$ is a power of $p$.
\end{cor}
\begin{proof}
Write $d=\gcd(d_1,\dots,d_r)=np^s$ for some $n\geq1$ with $p\nmid n$ and $s\geq 0$. Let $K(n)/K$ be the unique extension of degree $n$. Let $R$ and $R'$ denote the msnc reduction types of $C$ and $C(n)=C\times_K K(n)$, respectively. In the proof of Corollary \ref{cor: reduction type} we showed that $R'=R''_\text{min}$ for a reduction type cover $R''$ of $R$ of degree $n$. Let $\A=\ZZ/n\ZZ$, and let $\Gamma$ denote the dual graph of $R$. 
The assumption that the toric rank is zero implies that $\Gamma$ is a tree, so $H^1(\Gamma,\A_D)=0$. Moreover, since the abelian rank is zero, every $v\in V(\Gamma)$ has $g(v)=0$, and by assumption we have $n\mid \weight(v)$. Therefore, the dilation datum associated to any connected cyclic dilated cover of $\Gamma$ of degree $n$ satisfies $D(x)=0$ for every $x\in X(\Gamma)$. In particular, $\sum_{v\in V(\Gamma)} D(v)=0$, so Proposition \ref{non-existence of a connected cover} implies that there is no connected cyclic dilated cover of $\Gamma$ of degree $n$ if $n\geq 2$. This implies that, if $n\geq 2$, there is no reduction type cover of $R$ of degree $n$. Thus $n=1$, so $s\geq 1$ and $p>0$. 

\end{proof}

\begin{remark}
    A perhaps more intuitive argument is as follows. If $\C$ is the msnc model of $C$, then the dual graph $\Gamma_n$ of $(\C_n)_k$, defined in Subsection \ref{subsec: dual graph}, is a cyclic dilated cover of $\Gamma$ of degree $n$. If $D(x)=0$ for every $x\in X(\Gamma)$, then any harmonic $\A$-cover of $\Gamma$ with dilation datum $D$ is a topological (unramified) covering. Thus, if $\Gamma$ is a tree, its fundamental group is trivial, so any topological covering of degree $n$ is given by $n$ disconnected copies of $\Gamma$. However, since $C$ is geometrically connected, the special fibre $(\C_n)_k$ is connected, hence so is $\Gamma_n$. 
\end{remark}

\begin{remark}
    Corollary \ref{cor: Liu-Lorenzini-Raynaud} implies that the condition $q\nmid \weight(R)$ in part \eqref{a} of Corollary \ref{unique reduction type} can be removed. 
\end{remark}

The following corollary follows immediately by Remark \ref{weight of a chain}.

\begin{cor} \label{non-zero genus}
    Suppose the reduction type of $C/K$ has only one principal component $v$ and no cycles. If $n\mid \operatorname{weight}(v)$ for some $n\geq 1$ such that $p\nmid n$, then $g(v)\geq 1$. 
\end{cor} 

\begin{cor} \label{tame g=1}
    If $C$ is a tame curve with $g(C)=1$, bad reduction, and toric rank zero, then $C$ is an elliptic curve.
\end{cor}
\begin{proof}
    The assumptions that $C$ has bad reduction and $g(C)=1$ imply that the abelian rank of $C$ is zero. The result now follows from Corollary \ref{cor: Liu-Lorenzini-Raynaud}, Theorem \ref{Halle} and the classification of reduction types in genus 1 (see e.g. \cite{Saito} or \cite[Corollary 9.1]{Tim2025}). \\
\end{proof}

\begin{thm} \label{thm: tame curve}
    A tame curve of genus $g\geq 1$ has potentially good reduction if and only if its reduction type has exactly one principal component and no cycles.
\end{thm}
\begin{proof}
     Suppose $C/K$ is a tame curve and let $K(n)/K$ denote the minimal extension such that $C(n)$ has semistable reduction. Since $C$ is tame, we have $p\nmid n$. Let $\C$ and $\C(n)_{\operatorname{min}}$ be the msnc models of $C$ and $C(n)$, respectively. Recall that $\C(n)$ is a snc model by Proposition \ref{Halle-Nicaise 1.3.2}, so there is a map $\C(n)\to \C(n)_{\operatorname{min}}$. 
     
     Suppose that $C$ has potentially good reduction. This implies that the toric rank of $C(n)$ is $t(C(n))=0$. Since $t(C(n))\geq t(C)$, the reduction type of $C$ has no cycles. Thus, the msnc reduction type and the (mrnc) reduction type of $C$ coincide, and similarly for $C(n)$. Moreover, the reduction type of $C(n)$ has only one principal component $u$ with $g(u)=g(C)\geq 1$. By abuse of notation, denote also by $u$ the unique component in $\C(n)_k$ mapping surjectively to $u$ under the map 
     \[
     \C(n)_k\to (\C(n)_{\operatorname{min}})_k.
     \]
     Since $g(u)\geq 1$, the component $u$ in $\C(n)_k$ is also principal. By Propositions \ref{desingularisation} and \ref{Halle-Nicaise 1.3.2}(\ref{HN 2}), the image of $u$ under the map
     \[
     \varphi_k\circ\rho_k:\C(n)_k\to (\C_n)_k\to \C_k
     \]
     is also a principal component in the reduction type of $C$, and we denote it by $v$. 
     If $g(C)\geq 2$, Theorem \ref{non-contractibility of principal components} implies that any principal component in the reduction type of $C$ other than $v$ is a D-tail. This also holds if $g(C)=1$ and $t(C)=0$ by the classification of reduction types in genus 1. Suppose $v$ has a D-tail $v'$ of weight $m$ for some $m\geq 1$. By Theorem \ref{Halle}, we have $2m\mid n$, so by Proposition \ref{number of components g=0}, there are $m$ components in the $(\C_n)_k$ lying above $v'$, and $2m$ inner chains lying above the inner chain from $v'$ to $v$. Recall that $\Gal(K(n)/K)$ permutes the components in $(\C_n)_k$ lying above a given component in $\C_k$. Since there is only one component $v_1$ in $(\C_n)_k$ lying above $v$, each component in $(\C_n)_k$ lying above $v'$ is connected to $v_1$ by 2 such inner chains. This implies that there are at least $m$ cycles in the dual graph of $(\C_n)_k$. By Proposition \ref{desingularisation} and since contracting $(-1)$-curves does not change the number of cycles, the number of cycles in the reduction type of $C(n)$ is the same as in $\Gamma_n$. This contradicts that $C(n)$ has good reduction. Therefore the reduction type of $C$ has no D-tails and $v$ is the only principal component in it. 

     Conversely, suppose that the reduction type $R$ of $C$ has exactly one principal component $v$ and no cycles. Then $R$ is also a msnc reduction type. By Corollary \ref{non-zero genus}, either $n\nmid \weight(v)$, or $g(v)\geq 1$. In either of these cases, Proposition \ref{number of components divides} and Corollary \ref{one comp lying above} show that there is exactly one component $u$ in $(\C_n)_k$ lying above $v$. Moreover, since $R$ has only one principal component and no loops, it has no inner chains. Thus, Propositions \ref{Halle-Nicaise 1.3.2} and \ref{desingularisation} imply that the special fibre of the desingularisation $\C(n)_k$ is snc and consists also of exactly one principal component, which we also denote by $u$, and outer chains. Since $C(n)$ has semistable reduction, all components in $(\C(n)_{\operatorname{min}})_k$ have multiplicity 1. Thus, there are no outer chains in $(\C(n)_{\operatorname{min}})_k$, as this would contradict the minimality of $\C(n)_{\operatorname{min}}$. Since contracting outer chains does not create any loops, the reduction type of $C(n)$ consists only of one component $u$ with no loops, so $C(n)$ has good reduction. \\
\end{proof} 

\section{Computing the reduction type after tame base change} \label{sec: computing}

In this section, we give an algorithm for computing the set of reduction type covers of degree $n$ of a given reduction type. We also give a formula for the toric and abelian ranks of the base changed curve $C(n)$ in terms of the reduction type of $C$. To simplify the desingularisation step in the algorithm, we specify how outer and inner chains behave under tame base change in Section \ref{sec: chains}. For this, we need to first study the \emph{length} of a chain. Throughout this section, we will often denote $\gcd(a_1,\dots,a_r)$ by $(a_1,\dots,a_r)$.

\subsection{Length of chains} \label{sec: depth}
\mbox{}
\vspace{5pt}

Let $Z$ be a chain of type $d_0\frac{d_1}{}$ or $d_0\edge{d_0}{d_r}{n}d_{r+1}$. We call $d_0,d_1,\dots,d_{r+1}$ the consecutive multiplicities of $Z$.

\begin{defn}
     We define the \emph{length} of $Z$ to be \[\ell(Z)=\sum_{i=0}^{r} \frac{1}{d_{i}d_{i+1}}.\]
\end{defn}

\begin{lemma} \label{lengths under base chage}
    Let $Z$ be a chain in $\C_k$ and let $w=\operatorname{weight}(Z)$. Let $Z'$ be a chain in $\C(n)_k$ lying above $Z$ under the map $\C(n)\to \C$. Then
    \(
    \ell(Z')=n(w,n)\ell(Z).
    \)
    \end{lemma} 
    \begin{proof}
    Let $e_i$ denote the multiplicity of a component in $Z'$ lying above a component in $Z$ of multiplicity $d_i$. By Proposition \ref{Halle-Nicaise 1.3.2}, we have $e_i=d_i/(d_i,n)$. Note that there might be components in $Z'$ lying above points in $Z$, arising from the desingularisation described in Proposition \ref{desingularisation}. Suppose $x_i$ is a point in $(\C_n)_k$ lying above the intersection of two components $v_i,v_{i+1}$ of respective multiplicities $d_i,d_{i+1}$. The minimal desingularisation of the local germ $\operatorname{Spec}(\mathcal{O}_{\C_n,x})$ give rise to a chain of components with multiplicities $\mu_1,\mu_2,\dots,\mu_\lambda$ between the components in $\C(n)_k$ lying above $v_i,v_{i+1}$ (if $x_i$ is not singular, this chain is empty). Let 
        \[
        n_i=\frac{n(d_i,d_{i+1},n)}{(d_i,n)(d_{i+1},n)}.
        \]
        By \cite[Lemma 2.9.1]{Padma} we have
        \[
        \frac{n_i}{e_ie_{i+1}}=\frac{1}{e_i\mu_1}+\frac{1}{\mu_1\mu_2}+\cdots+\frac{1}{\mu_{\lambda-1}\mu_\lambda}+\frac{1}{\mu_\lambda e_{i+1}}.
        \]
        This implies that
        \[
        \ell(Z')=\sum_{i=0}^{r}\frac{n_i}{e_ie_{i+1}}=\sum_{i=0}^{r}\frac{n_i(d_i,n)(d_{i+1},n)}{d_id_{i+1}}=\sum_{i=0}^{r}\frac{n(d_i, d_{i+1},n)}{d_id_{i+1}}=n(w,n)\ell(Z),
        \]
        where the last equality holds because $w=\operatorname{weight}(Z)=(d_i,d_{i+1})$ for all $0\leq i\leq r$ (see Remark \ref{weight of a chain}). \\
    \end{proof}

\begin{remark} \label{lengths under contraction}
    Note that, if $\C(n)$ is not minimal, the chain $Z'$ might have contractible components. If so, let $Z''$ denote the (possibly empty) chain obtained by contracting all contractible components in $Z'$. If a contractible component has multiplicity $m_i$ and intersects two components of multiplicities $m_{i-1},m_{i+1}$, this does not affect the length. Indeed, in this case $m_i=m_{i-1}+m_{i+1}$, and thus
    \[
    \frac{1}{m_{i-1}m_i}+\frac{1}{m_im_{i+1}}=\frac{m_{i+1}+m_{i-1}}{m_{i-1}m_im_{i+1}}=\frac{1}{m_{i-1}m_{i+1}}.
    \]
    In particular, if $Z'$ is an inner chain, then $\ell(Z'')=\ell(Z')$. 
    If $Z'$ is an outer chain and a contractible component lies at the end of it, i.e. it only intersects one other component in $Z'$ once, then both of these components have multiplicity $w_n:=\operatorname{weight}(Z')$. Therefore, $\ell(Z'')\equiv \ell(Z')=n(w,n)\ell(Z) \text{ mod } \frac{1}{w_n^2}$. 
\end{remark}

\begin{prop} \label{length proposition} Let $Z$ be a chain with consecutive multiplicities $m,d_1,\dots,d_r,m'$ in the special fibre of a mrnc model $\C$. Let $w=\operatorname{weight}(Z)(=(m,d_1))$, and $s=\operatorname{depth}(Z)$. Then the length of $Z$ is given by
    \[
    \ell(Z)=\frac{i_{d_1,m}}{mw}+\frac{i_{d_r,m'}}{m'w}+\frac{s}{w^2},
    \]
    where $i_{d_1,m}\geq 0$ is the smallest integer such that $i_{d_1,m}d \equiv (d_1,m) \text{ mod }m$, and similarly for $i_{d_r,m'}$. In particular, if $Z$ is an outer chain, then
    \[
    \ell(Z)=\frac{i_{d_1,m}}{mw}.
    \]
\end{prop}
\begin{proof}
    Assume $w=1$. By \cite[Definition 5.9 and Theorem 5.14]{Tim2025}, the sequence $m=d_0,d_1,\dots,d_r,d_{r+1}=m'$ corresponds to a 1-path, i.e. a sequence of fractions in lowest terms for some $r\geq 0$,
    \begin{equation*} \label{1-path}
         a=\frac{h_0}{m}=\frac{h_0}{d_0}>\frac{h_1}{d_1}>\dots>\frac{h_r}{d_r}>\frac{h_{r+1}}{d_{r+1}}=\frac{h_{r+1}}{m'}=a',
    \end{equation*}
    such that 
    \begin{equation}\label{1-path condition}
        h_id_{i+1}-d_ih_{i+1}=1,
    \end{equation} 
    for $0\leq i\leq r$ \cite[Appendix A.4]{ObusWewers}. This implies that 
    \[
   \frac{h_i}{d_i}-\frac{h_{i+1}}{d_{i+1}}=\frac{1}{d_id_{i+1}},
    \]
    and therefore,
    \[
    \ell(Z)=\sum_{i=0}^r \frac{1}{d_{i}d_{i+1}}=\sum_{i=0}^{r}\frac{h_i}{d_i}-\frac{h_{i+1}}{d_{i+1}}=\frac{h_0}{m}+\frac{-h_{r+1}}{m'}.
    \]
    Let $q_0,q_1\in \ZZ$ be such that $h_0=q_0m+r_0$ and $h_{r+1}=q_1m'+r_1$, for some $0\leq r_0<m$ and $0\leq r_1<m'$. From Equation \eqref{1-path condition} we have
    \[
    q_0md_1+r_0d_1-h_1m=1, \text{ and } \ h_rm'-d_rq_1m'-d_rr_1=1,
    \]
    so $r_0d\equiv 1\equiv (d,m) \operatorname{mod } m$ and $-r_1d'\equiv 1\equiv (d',m') \operatorname{mod } m'$. Since $0\leq r_0<m$ and $0\leq r_1<m'$, we have $r_0=i_{d_1,m}$ and $m'-r_1=i_{d_r,m'}.$
    Then
    \[
    \ell(Z)=\frac{q_0m+r_0}{m}+\frac{-q_1m'-r_1}{m'}=\frac{i_{d_1,m}}{m}+\frac{i_{d_r,m'}}{m'}+q_0-q_1-1.
    \]
    By definition, the depth $s$ of the chain $Z$ satisfies $s=-1+\#\{d_i=1\}$. Then Equation \eqref{1-path condition} implies that $s=-1+\#\{\text{integers $n$ such that } a\leq n\leq a'\}=-1+q_0-q_1$. This proves the result for $w=1$.
    
    Suppose now that $w>1$. Define $e_i=d_i/w$ and let $Z'$ be some chain (not necessarily in $\C$) with consecutive multiplicities $e_0,\dots,e_{r+1}$. Then 
    \[
    \ell(Z')=w^2\ell(Z).
    \]
    Since $(d_i,d_{i+1})=\operatorname{weight}(Z)=w$ for all $0\leq i \leq r$, we have $\operatorname{weight}(Z')=(e_i,e_{i+1})=1$, so we can compute the length of $Z'$ as above to get 
    \[
    \ell(Z')=\frac{i_{e_1,e_0}}{e_0}+\frac{i_{e_r,e_{r+1}}}{e_{r+1}}+s,
    \]
    where $s=\operatorname{depth}(Z)=\operatorname{depth}(Z')$.
    Since \[1=i_{e_1,e_0}e_1-je_0=i_{e_1,e_0}\frac{d_1}{w}-j\frac{d_0}{w}\] for some $j\in \ZZ$, and $w=(d_0,d_1)$, we have $i_{e_1,e_0}=i_{d_1,d_0}$. Similarly, $i_{e_r,e_{r+1}}=i_{d_r,d_{r+1}}$. 
    Putting everything together gives
    \[
    \ell(Z)=\frac{1}{w^2}\ell(Z')=\frac{1}{w^2}\left(\frac{i_{d_1,d_0}}{d_0/w}+\frac{i_{d_r,d_{r+1}}}{d_{r+1}/w}+s\right)=\frac{i_{d_1,m}}{mw}+\frac{i_{d_r,m'}}{m'w}+\frac{s}{w^2}.
    \]

    If $Z$ is an outer chain, then $d_{r+1}=w$ and so $i_{d_r,d_{r+1}}=0$. And since $\C$ is minimal, there is exactly one component (the last component) of multiplicity $w$ (see e.g. \cite[Corollary 5.6]{Tim2025}). \\
\end{proof}

\subsection{Behaviour of chains under tame base change} \label{sec: chains}
\mbox{}
\vspace{5pt}

By Remark \ref{weight of a chain}, Proposition \ref{Halle-Nicaise 1.3.2}, and Fact \ref{number of points above intersection point}, we know the number of chains lying above each chain after tame base change. Therefore we are left to study what each chain lying above a chain looks like. As explained in Section \ref{sec: reduction type}, an outer chain in an mrnc model is determined by its initial multiplicity. An inner chain from $v$ to $u$ is determined by the initial multiplicities at $v$ and $u$ and its depth. In this section, we determine how the initial multiplicities and the depth of chains behave under tame base change. Let $\C$ be the mrnc model of $C$, and let $\C_n$ and $\C(n)$ be as before. 

\begin{prop} \label{initial multiplicity of chain after base change}
 Let $n$ be prime. Let $v$ be a principal component in $\C_k$ of multiplicity $m$. Let $Z$ be an outer chain from $v$, and denote by $d$ the initial multiplicity of $Z$ at $v$. Suppose $Z'$ is the image under the map $\C(n)_k\to (\C(n)_{\operatorname{min}})_k$ of a chain in $\C(n)_k$ lying above $Z$. Then the initial multiplicity $d'$ of $Z'$ is given by
 \[
   d'=\left\{
\begin{array}{lr}
n^{-1}d \  \operatorname{mod} m & \text{if } n\nmid m,\\
d \ \operatorname{mod} m/n & \text{if } n\mid m, \  n\nmid d, \\
d/n \  & \text{if } n\mid m, \ n\mid d,\\
\end{array}
\right.
\]
where we let $d'=0$ if $Z'$ is empty.
\end{prop}
\begin{proof}
    Let $v''$ and $Z''$ be a component and a chain in $\C(n)_k$ lying above $v$ and $Z$ under the map $\C(n)\to \C$, respectively. Let $Z'$ and $v'$ be the respective images of $Z''$ and $v''$ under the map $\C(n)\to \C(n)_{\operatorname{min}}$. Recall that, by Proposition \ref{Halle-Nicaise 1.3.2}, the multiplicity of $v''$ is $m/(m,n)$. Moreover, since $v''$ is a principal component, it does not get contracted by Theorem \ref{non-contractibility of principal components}, so $v'$ is also a component of multiplicity $m/(m,n)$. Let $d'$ (resp. $d_1$) denote the initial multiplicity of $Z'$ (resp. $Z''$), i.e. the multiplicity of the component $v_1'$ in $Z'$ (resp. $v_1''$ in $Z''$) intersecting $v'$ (resp. $v''$). Since $Z'$ is obtained by contracting all the $(-1)$-curves in $Z''$, it is possible that no such $v_1'$ exists, in which case we let $d'=0$. 
    Since $Z$ is an outer chain, $Z'$ is also an outer chain. This implies that $0\leq d'<m/(m,n)$ (see e.g. \cite[Corollary 5.6]{Tim2025}). Therefore, it is enough to know the value of $d' \text{ mod } m/(m,n)$. 

    Suppose that $n\nmid m$. By Lemma \ref{weight lemma}, the weight of $Z''$ is $w/(w,n)=w/(m,d,n)=w$. Then Lemma \ref{lengths under base chage} and Remark \ref{lengths under contraction} imply that 
    \[
    \ell(Z')\equiv \ell(Z'')\equiv n\ell(Z) \text{ mod } \frac{1}{w^2}.
    \]
    Since the multiplicity of $v''$ is $m/(m,n)=m$, from Proposition \ref{length proposition} we obtain that
    \[
    \frac{i_{d',m}}{mw}\equiv n\frac{i_{d,m}}{mw} \text{ mod } \frac{1}{w^2},
    \]
    or, equivalently, 
    \begin{equation} \label{inverse equation}
        i_{d',m}\equiv ni_{d,m} \text{ mod } \frac{m}{w}.
    \end{equation}
    Recall that
    \[
    i_{d',m}\frac{d'}{w}\equiv 1\equiv i_{d,m}\frac{d}{w}\text{ mod } \frac{m}{w},
    \]
    by definition of $i_{d',m}$ and $i_{d,m}$. Then, multiplying by $d/w$, $d'/w$ and $n^{-1}$ on both sides of  \eqref{inverse equation}, we obtain
    \[
    n^{-1}\frac{d}{w} \equiv \frac{d'}{w} \text{ mod } \frac{m}{w}, 
    \]
    or, equivalently,
    \[
    d' \equiv n^{-1}d  \text{ mod } m. 
    \]

    Suppose now that $n\mid m$. Then, regardless of whether $n\mid d$ or $n\nmid d$, we have that  
    \[
    n'=\frac{n(m,d,n)}{(m,n)(d,n)}=1.
    \]
    By Proposition \ref{desingularisation}, this implies that the intersection of any two components in the normal model $\C_n$ lying above $v$ and the initial component $v_1$ of the chain $Z$, respectively, is already regular. By Proposition \ref{Halle-Nicaise 1.3.2}, the multiplicity of the component lying above $v_1$ in the special fibre $(\C_n)_k$ (and hence, in $\C(n)_k$) is $d_1=d/(d,n)$. 
    
    We next show that $d'\equiv d_1 \text{ mod } m/n$. We begin by contracting all $(-1)$-curves in $Z''$ other than $v_1''$. Let $d_1,d_2,d_3,\dots, d_r$ denote the consecutive multiplicities of the chain from $v_1''$ obtained by such contraction. 
    We write $d_2=0$ if this chain consists only of $v_1''$
    . The self-intersection of $v_1''$ is
    \[
    (v_1'')^2=-\frac{\frac{m}{n}+d_2}{d_1} \in \ZZ.
    \]
    Thus, $v_1''$ is contractible if and only if 
    \[
    d_1=\frac{m}{n}+d_2.
    \]
    Note that $d_1,d_2,\dots,d_r$ is an outer sequence, so it is strictly decreasing \cite[Definition 5.1 and Corollary 5.6]{Tim2025}. Therefore, $v_1''$ is contractible if and only if $d_1\geq m/n$. 
    If $v_1''$ is not contractible, then $d'=d_1$. If it is contractible, then we contract it and repeat the same argument with the component $v_2''$ of multiplicity $d_2$. Note that, if $v_1''$ is contractible, then $d_2\equiv d_1 \text{ mod } m/n$. If $d_2<m/n$, then $d'=d_2$, and otherwise we contract $v_2''$ and repeat the argument with the component of multiplicity $d_3$. Iterating this argument gives that the initial multiplicity $d'$ of $Z'$ satisfies $d'\equiv d_1 \text{ mod } m/n$.

\end{proof}

\begin{remark} \label{inner chain remark}
    Since the multiplicity $d_1$ can be obtained via the local computation described in Proposition \ref{desingularisation}, its value does not depend on whether $Z$ is an inner or an outer chain. Moreover, a similar argument to the last part of the above proof shows that $d'\equiv d_1 \text{ mod } m/(m,n)$, though note that, if $Z$ is an inner chain, $d'$ need not satisfy $d'<m/(m,n)$. 
\end{remark}

\begin{prop} \label{depth}
    Let $n$ be prime. Let $Z$ be an inner chain of type $m_1\edge{d_1}{d_2}{s}m_2$ in $\C_k$. Let $w=\operatorname{weight}(Z)=(m_i,d_i)$. Then the depth of an inner chain in $(\C(n)_{\operatorname{min}})_k$ lying above $Z$ is 
    \[
    s'= \left\{
    \begin{array}{lr}
    s  & \text{if } n\mid m_i, \  n\mid d_i,\\
     sn+s_1+s_2 & \text{otherwise,}
    \end{array}
    \right.
    \] 
    where for $i=1,2$, we define
    \[
    s_i=\left\lfloor{\frac{i_{d_i,m_i}nw}{m_i}}\right\rfloor.
    \]
\end{prop}
\begin{proof}
    Note that $n\mid m_1$ and $n\mid d_1$ holds if and only if $n\mid m_2$ and $n\mid d_2$, since $w=(m_1,d_1)=(m_2,d_2)$. 
    
    \textbf{Case 1:} $n\mid m_i$ and $n\mid d_i$ for $i=1,2$.
    
    If $n\mid m_i$ and $n\mid d_i$ for some $i=1,2$, then $n\mid w$. By Remark \ref{weight of a chain}, for any two consecutive components $u,v$ in the chain $Z$ with multiplicities $a,b$ we have $w=(a,b)$, so $n\mid a$ and $n\mid b$. Thus,
    \[
    n'=\frac{n(a,b,n)}{(a,n)(b,n)}=1.
    \]
    By Proposition \ref{desingularisation} this implies that the intersection of any two components in the normal model $\C_n$ lying above the intersection of $u$ and $v$ is already regular. In particular, the depth is preserved, i.e. $s'=s$.  

    \textbf{Case 2:} $n\nmid w=(m_1,d_1)=(m_2,d_2)$.
    
    Let $Z''$ be a chain in $\C(n)_k$ lying above $Z$ under the map $\C(n)\to \C$. Let $Z'$ be the image of $Z''$ under the map $\C(n)\to \C(n)_{\operatorname{min}}$.
    Since $Z$ is an inner chain, Lemma \ref{lengths under base chage} and Remark \ref{lengths under contraction} imply that 
    \[
    \ell(Z'')=\ell(Z')=n\ell(Z).
    \]
    For $i=1,2$, let $m_i'=m_i/(m_i,n)$ and $d_i'$ be the initial multiplicity of $Z'$. From Proposition \ref{length proposition} we obtain that
    \begin{equation} \label{comparing lengths equation}
    \frac{i_{d_1',m_1'}}{m_1'w}+\frac{i_{d_2',m_2'}}{m_2'w}+\frac{s'}{w^2}= n\frac{i_{d_1,m_1}}{m_1w}+n\frac{i_{d_2,m_2}}{m_2w}+n\frac{s}{w^2}.
    \end{equation}
    By Remark \ref{inner chain remark}, the value of $d_i' \text{ mod } m_i'$ is given by Proposition \ref{initial multiplicity of chain after base change}. Defining
    \[
    s_i:=\frac{i_{d_i,m_i}nw}{m_i}-\frac{i_{d_i',m_i'}w}{m_i'},
    \]
    we can rewrite $s'=sn+s_1+s_2$.  We next show that $s_i\in \ZZ$ and give a formula for $s_i$ that does not depend on $d_i'$ or $m_i'$.
   
    \textbf{Case 2a:} $n\mid m_i$ and $n\nmid d_i$.

    By Proposition \ref{initial multiplicity of chain after base change}, we have $d_i'\equiv d_i \text{ mod } m_i/n$. Then 
    \[
    i_{d_i', m_i'} \frac{d_i}{w}\equiv 1 \equiv i_{d_i,m_i}\frac{d_i}{w} \text{ mod } \frac{m_i}{nw}.
    \]
    Note that we can divide by $d_i/w$ since it is a unit mod $m_i/w$. Thus we get that
    \[
    s_i=\frac{i_{d_i,m_i}-i_{d_i',m_i'}}{m_i/(nw)}\in \ZZ.
    \]
    Since $0\leq i_{d_i',m_i'}<m_i'/w=m_i/(nw)$, and $s_i\in \ZZ$, we can write 
    \[
    s_i=\left\lfloor{\frac{i_{d_i,m_i}}{m_i/(nw)}}\right\rfloor.
    \]

    \textbf{Case 2b:} $n\nmid m_i$.
    
    By Proposition \ref{initial multiplicity of chain after base change}, we have $d_i'\equiv n^{-1}d_i \text{ mod } m_i$. Then 
    \[
    i_{d_i', m_i'} \frac{n^{-1}d_i}{w}\equiv 1 \equiv i_{d_i,m_i}\frac{d_i}{w} \text{ mod } \frac{m_i}{w}.
    \]
    Note that we can divide by $d_i/w$ and multiply by $n$ since both are units mod $m_i/w$. Therefore we get that
    \[
    s_i=\frac{i_{d_i,m_i}n-i_{d_i',m_i'}}{m_i/w}\in \ZZ.
    \]
    Since $0\leq i_{d_i',m_i'}<m_i/w$, and $s_i\in \ZZ$, we can write 
    \[
    s_i=\left\lfloor{\frac{i_{d_i,m_i}n}{m_i/w}}\right\rfloor.
    \] 
\end{proof}

By \cite[Theorems 5.3 and 5.14(3)]{Tim2025}, an inner chain of type $m_1\edge{d_1}{d_2}{s}m_2$ is uniquely determined by its depth and the values of $d_i \text{ mod } m_i$. Thus, Remark \ref{inner chain remark} and Proposition \ref{depth} determine the behaviour of inner chains under tame base change.

\subsection{Algorithm} \label{algorithm}
\mbox{}
\vspace{5pt}

In this section, we explain how to compute explicitly the possibilities for the reduction type of a curve after tame base change. We follow the method explained in the proof of Corollary \ref{unique reduction type}. We will identify components in a reduction type with the corresponding vertices in the dual graph. Let $R$ be an (mrnc) reduction type, and $\Gamma$ its dual graph. Let $n$ be an integer such that $p\nmid n$, and set $\A=\ZZ/n\ZZ$. Let $C/K$ be a curve with reduction type $R$, and $C(n)$ its base change by the unique extension of degree $n$. Then we compute a set $\mathcal{R}(n)$ of possible reduction types of $C(n)$ by computing all the reduction type covers of $R$, and taking $R'_{\text{min}}\in \mathcal{R}(n)$ for each reduction type cover $R'\to R$ (see Section \ref{sec: defn of reduction type covers}). 
\vspace{4pt}

\textbf{Input:} the reduction type $R$ of $C$ and a degree $n$.

\textbf{Output:} a finite set $\mathcal{R}(n)$ of possible reduction types of $C(n)$.
\vspace{5pt}

Steps:
\begin{enumerate}[1.]  \setlength{\itemsep}{4pt} 
    \item Compute the set $\text{Cov}$ of connected cyclic covers of $\Gamma$ of degree $n$ with dilation compatible with $R$. To be precise,
    \begin{enumerate} 
        \item [1.1] Let $\mathcal{D}$ be the set of $\A$-dilation data satisfying the conditions in Definition \ref{defn: compatible dilation}.
        \item [1.2] For each $D\in \mathcal{D}$ compute the dilated cohomology group $H=H^1(\Gamma,\A_D)$. 
        \item [1.3] For each cyclic subgroup $G\subseteq H$, pick a generator $[\eta_G]\in G$. Then, compute the harmonic $\A$-cover $\Gamma'\to \Gamma$ associated to $[\eta_G]$ (see \cite[pages 9-10]{Len-Ulrisch-Zakharov}). 
        If $\Gamma'$ is disconnected, we discard it. Otherwise, $\Gamma'\to \Gamma$ is a connected cyclic cover of $\Gamma$ of degree $n$. This gives every element of $\text{Cov}$.   
    \end{enumerate}
    \item For each $\Gamma'\to \Gamma\in \text{Cov}$, let $R'$ and $\Gamma''\to \Gamma$ denote the associated reduction type cover. We obtain the reduction type $R'_{\text{min}}\in\mathcal{R}(n)$ as follows. 
    \begin{enumerate}
        \item [2.1] For each principal component $v$ in $R$, we consider the vertices $u_1,\dots,u_\ell$ in $\Gamma''$ lying above $v\in \Gamma$. Let $m_v$ be the multiplicity of $v$ in $R$. Then, each $u_i$ has multiplicity $m_v/\gcd(m_v,n)$ and geometric genus $g(u_i)$ given by Equation \eqref{eqn: genus formula} in $R'$.
        \item [2.2] For each outer (resp. inner) chain $Z$ in $R$, let $Z_1,\dots,Z_\ell$ denote the connected components in the preimage of $Z$ under the map $\Gamma''\to \Gamma$. Then $Z_i$ are outer (resp. inner) chains in $R'$. Let $Z_i'$ denote the image of $Z_i$ in $R'_\text{min}$. Then, the initial multiplicity of $Z_i'$ is determined by Proposition \ref{initial multiplicity of chain after base change} (see Remark \ref{inner chain remark} for inner chains). If $Z$ is an inner chain, the depth of $Z_i'$ is determined by Proposition \ref{depth}. We note that one has to go prime-by-prime to apply these propositions.
        \item [2.3] Compute the multiplicities and number of components in each chain in $R'_\text{min}$ from the initial multiplicities and depth using \cite[Corollary 5.6 and Algorithm 5.16]{Tim2025}.
        \end{enumerate}
\end{enumerate}

\begin{remark}
    In Step 1.3, we discard disconnected covers manually after constructing them. Alternatively, one can discard them before computing the covers using Proposition \ref{prop: connected covers}. 
\end{remark}

\begin{remark}
    An alternative for Steps 2.1--2.3 is as follows. Do Step 2.1 for every component, not just principal ones. Then apply Proposition \ref{desingularisation}. Lastly, contract all non-principal components $v$ with $v\cdot v=-1$.  
\end{remark}

\subsection{Abelian and toric rank after tame base change} \label{sec: abelian and toric ranks}
\mbox{}
\vspace{5pt}

Recall that the abelian and toric rank can be computed from the reduction type (see Section \ref{sec: reduction type}). In this section we use our results to give a formula for the abelian and toric ranks of a curve $C$ (or its Jacobian) after tame base change of degree $n$ in terms of the reduction type of $C$. 
More specifically, let $C(n)$ and $J(n)$ denote the base change of $C$ and its Jacobian $J$ by the unique tame extension of degree $n$. We give formulas for the abelian and toric ranks
\[
a(C(n))=a(J(n))\ \text{ and }\ t(C(n))=t(J(n)),
\] 
in terms of the mrnc reduction type $R$ of $C$. Let $\I$ denote the set of inner chains in $R$ (recall that this includes the intersection of principal components and nodes). We denote the set of components in $R$ by $S$, and by $\mathcal{P}\subseteq S$ the set of principal ones. We write $m_v$ for the multiplicity of $v\in S$ and $w(-):=\weight(-)$.

\begin{prop}\label{abelian and toric rank} Assume $g(C)\geq 2$, and that the reduction type $R$ of $C$ satisfies the conditions in Corollary \ref{unique reduction type}, or it has no components of geometric genus $g\geq 1$. Then the toric rank of $C(n)$ and $J(n)$ is given by the following sum.
\[
t(J(n))=t(C(n))=\sum_{Z \in \I}  (w(Z),n)-\sum_{v\in \mathcal{P},\  g(v)=0}  (w(v),n) - \#\{v\in S:g(v)\geq 1\} +1.
\]
For any principal component $v$ in $R$ intersecting the rest of $R$ at $r_v$ points, denote by $d_{v,1},\dots,d_{v,r_v}$ the multiplicities of the corresponding components. Then the abelian rank of $C(n)$ and $J(n)$ is given by 
\begin{align*}
  a(J(n))=a(C(n)) & =\sum_{v\in \mathcal{P}} \left( (m_v,n)(g(v)-1)+\frac{1}{2}\sum_{i=1}^{r_v}( (m_v,n)- (m_v,d_{v,i},n))\right) \\ &+\sum_{v\in \mathcal{P},  g(v)=0}  (w(v),n)+\#\{v\in S: g(v)\geq 1\}.
\end{align*}
In particular, if $n$ is coprime to $m_v$ for every $v\in \mathcal{P}$, then $t(J(n))=t(J)$ and $a(J(n))=a(J)$.
\end{prop}

\begin{proof}
    Since the abelian and toric rank of $C$ can be computed from the rnc reduction type arising from any rnc model, assume that $\C$ is the msnc model of $C$ and $R$ the reduction type arising from $\C_k$. Let $R_n$ be the msnc reduction type of $C(n)$. Given a component $v\in S$, we denote by $\ell_v$ the number of components in $\C(n)_k$ lying above $v$ under the map $ \varphi \circ \rho : \C(n)\to \C$, and let $u_v$ be any such component. If $v\notin \mathcal{P}$, then $v$ is a component in a (outer, resp. inner) chain, so Proposition \ref{Halle-Nicaise 1.3.2}\eqref{HN 2} implies that $u_v$ is part of a (outer, resp. inner) chain, so $g(u_v)=0$. Since the maps $\rho: \C(n)\to \C_n$ and $\C(n)\to \C(n)_\text{min}$ only add or contract non-principal components in the special fibre, any principal component in $R_n$ arises from some $v\in \mathcal{P}$. Moreover, if $v\in \mathcal{P}$, then $u_v$ does not get contracted under the map $\C(n)\to \C(n)_\text{min}$ by Theorem \ref{non-contractibility of principal components}. Thus, the genus formula gives 
    \begin{align*}
        a(C(n)) & =\sum_{v\in \mathcal{P}}\ell_vg(u_v)=\sum_{v\in \mathcal{P}}\ell_v\left(\frac{(m_v,n)}{\ell_v}(g(v)-1)+\frac{1}{2}\sum_{i=1}^{r_v}\frac{(m_v,d_{v,i},n)}{\ell_v}\left(\frac{(m_v,n)}{(m_v,d_{v,i},n)}-1\right)+1\right) \\
        & =\sum_{v\in \mathcal{P}}\left((m_v,n)(g(v)-1)+\frac{1}{2}\sum_{i=1}^{r_v}((m_v,n)-(m_v,d_{v,i},n))+\ell_v\right).
    \end{align*}
    Recall that Proposition \ref{number of components g=0} and Corollary \ref{one comp lying above} show that $\ell_v= (w(v),n)$ for components $v$ with $g(v)=0$, and $\ell_v=1$ whenever $g(v)\geq 1$ under the assumptions on $R$. This implies the result for $a(C(n))$. 
    
     Let $\Gamma$ be the dual graph of $R$ and let $\Gamma(n)$ be the dual graph of $\C(n)_k$. Since $C$ is geometrically connected, the graph $\Gamma(n)$ is also connected, so
    \[
    t(C(n))=|E(\Gamma(n))|-|V(\Gamma(n))|+1.
    \]
    We identify components and chains in $R$ with the corresponding vertices and subgraphs in $\Gamma$. Suppose $Z$ is an outer chain from a principal component $v$ in $R$, and let $Z$ also denote the subgraph (with legs) of $\Gamma$ corresponding to the chain, including the edge to the vertex $v$ but not $v$ itself. Then $|E(Z)|-|V(Z)|=0$. Therefore, outer chains in $R$ do not contribute to $t(C)$. If $Z$ is an inner chain from principal components $v_1,v_2$ in $R$, and $Z$ also denotes the subgraph of $\Gamma$ corresponding to $Z$ (with the edges to $v_i$, but without the vertices $v_i$), then $|E(Z)|-|V(Z)|=1$. Thus, we have
    \[
    t(C)=|E(\Gamma)|-|V(\Gamma)|+1=|\I|-|\mathcal{P}|+1.
    \]
    Note that we can apply the same argument to $\Gamma(n)$. As argued above, the number of principal components in $\C(n)_k$ is $\sum_{v\in \mathcal{P}}\ell_v$. The number of inner chains in $\C(n)_k$ is the same as in $(\C_n)_k$. Moreover, each inner chain in $\C(n)_k$ lies above an inner chain in $\C_k$. By Remark \ref{weight of a chain} and Proposition \ref{Halle-Nicaise 1.3.2}, for each inner chain $Z$ in $\C_k$, there are $(w(Z),n)$ inner chains in $\C(n)_k$ lying above it. Together with Proposition \ref{number of components g=0}, this implies that
    \begin{align*}
    t(C(n)) & =\sum_{Z\in \I}  (w(Z),n)-\sum_{v\in \mathcal{P}}\ell_v +1 \\
    & = \sum_{Z\in \I}  (w(Z),n)-\sum_{v\in \mathcal{P}, g(v)=0} (w(v),n) - \#\{v\in S: g(v)\geq 1\}+1.
    \end{align*}
    Since the msnc or mrnc reduction type of $C$ only differs in the number of components in each chain, both formulas are independent of this choice.

    If $n$ is coprime to $m_v$ for every $v\in \mathcal{P}$, the formulas above simplify as follows.
    \[
    a(C(n))=\sum_{v\in \mathcal{P}}(g(v)-1)+|\mathcal{P}|=a(C),
    \]
    and 
    \[
    t(C(n))= |\I|-|\mathcal{P}|+1=t(C).
    \]
\end{proof}

\begin{remark}
    If $n$ is coprime to $m_v$ for every $v\in \mathcal{P}$, it can be shown similarly that $t(J(n))=t(J)$ holds without any assumptions on the reduction type of $C$, since in this case we have
    \[
    \sum_{v\in \mathcal{P}}   (w(v),n)=|\mathcal{P}|.
    \]
    This recovers \cite[Lemma 4.3.1.3]{Halle-Nicaise}. 
\end{remark}

\begin{example}
    Suppose a curve $C$ (of genus 3) has reduction type $R$ as below. Let us compute the toric and abelian ranks of the base change of $C$ by a tame extension of degree 4. We can apply Proposition \ref{abelian and toric rank} since $R$ has no components of geometric genus $g\geq 1$ (or since $R$ satisfies condition \eqref{c} in Corollary \ref{unique reduction type} for $n=4$). 
    
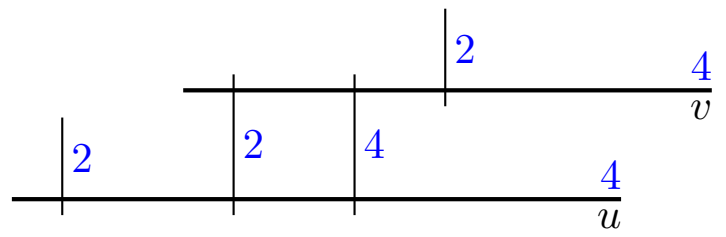
\begin{figure}[H]
    \vspace{-15pt}
        \centering
\begin{tikzpicture}[xscale=1,yscale=0.9]
\draw[thick,shorten >=2] (-1/3,0)--(113/30,0) node[inner sep=2pt,blue,above 
left,scale=0.8] {4} node[inner sep=2pt,below left,scale=0.8] {$u$};
\draw[shorten <=-3pt] (0,0)--node[inner sep=2pt,scale=0.8,blue,right] {2} 
(0,3/5);
\draw[shorten <=-3pt, shorten >=-3pt] (17/15,0)--node[inner 
sep=2pt,scale=0.8,blue,right] {2} (17/15,4/5);
\draw[shorten <=-3pt, shorten >=-3pt] (29/15,0)--node[inner 
sep=2pt,scale=0.8,blue,right] {4} (29/15,4/5);
\draw[thick,shorten >=2] (4/5,4/5)--(131/30,4/5) node[inner sep=2pt,blue,above 
left,scale=0.8] {4} node[inner sep=2pt,below left,scale=0.8] {$v$};
\draw[shorten <=-3pt] (38/15,4/5)--node[inner sep=2pt,scale=0.8,blue,right] {2} 
(38/15,7/5);
\end{tikzpicture} 
\vspace{-10pt}
    \end{figure}
    \noindent There are exactly 2 inner chains in $R$, one of weight 2 and one of weight 4. Then
    \[
    t(C(4))=(4,4)+(2,4)-(\weight(u),4)-(\weight(v),4)+1=4+2-2-2+1=3.
    \]
    To compute the abelian rank of $C(4)$, note that the reduction type is symmetric in $u$ and $v$. Then the formula simplifies to
    \[
    a(C(4))=2\left((4,4)(-1)+\frac{1}{2}\left((4,4)-(4,4,4)+2\left((4,4)-(4,2,4)\right)\right)\right)+2(\weight(v),4)=0.
    \]
    Alternatively, if the curve $C$ is tame (which holds if and only if $p\neq 2$ by a theorem of Saito \cite[Theorem 10.4.47]{Liu}, \cite[Theorem 3.11]{Saito}), then Theorem \ref{Halle} implies that $C(4)$ has semistable reduction. Thus 
    \[
    a(C(4))+t(C(4))=g(C)=3,
    \]
    so it is enough to compute one of $t(C(4))$, $a(C(4))$ to obtain the other one.
\end{example}

\section{Examples}\label{sec: examples}

In this section, given a reduction type $R$ and a degree $n$, we compute reduction type covers $R''\to R$ of degree $n$, and obtain $R'=R''_\text{min}$. Recall that if a curve has reduction type $R$, then its reduction type after tame base change of degree $n$ will be given by such an $R'$.

We start with an example of a reduction type satisfying the conditions in Corollary \ref{unique reduction type} for some $n$, and use the Algorithm in Section \ref{algorithm} to compute its unique reduction type cover of degree $n$. In Example 2, we look at a reduction type not satisfying the conditions for some $n$, and compute all its possible reduction type covers of degree $n$. In Example 3, we observe that computing the set of possible reduction type covers of degree $n$ by going prime-by-prime may give a larger set than if we compute it directly for composite $n$ using the algorithm. 

\subsection*{Example 1}
\mbox{}
\vspace{5pt}

Let us compute the reduction type $R'$ after tame base change by an extension of degree 3 of the reduction type $R$ below, which is of genus 4. We follow the Algorithm in Section \ref{algorithm}. Let $\A=\ZZ/3\ZZ$.
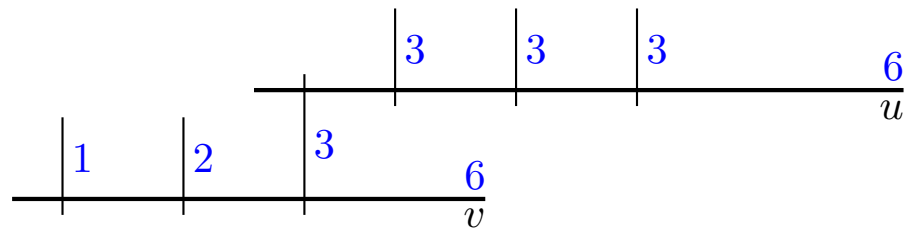
\begin{figure}[H]
\vspace{-5pt}
\begin{tikzpicture}[xscale=1,yscale=0.9]

\draw[thick,shorten >=2] (-1/3,0)--(43/15,0) node[inner sep=2pt,blue,above 
left,scale=0.8] {6} node[inner sep=2pt,below left,scale=0.8] {$v$};
\draw[shorten <=-3pt] (0,0)--node[inner sep=2pt,scale=0.8,blue,right] {1} 
(0,3/5);
\draw[shorten <=-3pt] (4/5,0)--node[inner sep=2pt,scale=0.8,blue,right] {2} 
(4/5,3/5);
\draw[shorten <=-3pt, shorten >=-3pt] (8/5,0)--node[inner 
sep=2pt,scale=0.8,blue,right] {3} (8/5,4/5);
\draw[thick,shorten >=2] (19/15,4/5)--(169/30,4/5) node[inner sep=2pt,blue,above 
left,scale=0.8] {6} node[inner sep=2pt,below left,scale=0.8] {$u$};
\draw[shorten <=-3pt] (11/5,4/5)--node[inner sep=2pt,scale=0.8,blue,right] {3} 
(11/5,7/5);
\draw[shorten <=-3pt] (3,4/5)--node[inner sep=2pt,scale=0.8,blue,right] {3} 
(3,7/5);
\draw[shorten <=-3pt] (19/5,4/5)--node[inner sep=2pt,scale=0.8,blue,right] {3} 
(19/5,7/5);
\end{tikzpicture}
\vspace{-10pt}
\end{figure}

\noindent Since all the components in $R$ have geometric genus 0, there is a unique $\A$-dilation datum $D$ compatible with $R$ given by $D(w)=\ZZ/({3/\gcd(\weight(w),3)})\ZZ$ for any component $w$ in $R$. Moreover, since the dual graph $\Gamma$ of $R$ is a tree, we have $H^1(\Gamma,\A_D)=0$. Therefore, there is a unique cyclic cover $\Gamma'$ of $\Gamma$ of degree 3 with dilation datum $D$. Below we have drawn the dilated cover $\Gamma\leftarrow \Gamma'$. Recall that, from the the dilation datum, we know that each vertex $w$ in $\Gamma$ has $\gcd(\weight(w),3)$ components lying above it in $\Gamma'$. For $i=1,2,3$, the vertices $u_i,v_1\in\Gamma'$ map to $u,v\in V(\Gamma)$, respectively. More generally, a vertex in $\Gamma'$ maps to the vertex in $\Gamma$ of the same colour. The numbers represent the multiplicities of the corresponding components. In the case of $\Gamma'$, we have only included the multiplicities of the vertices corresponding to principal components, which are $m_w/\gcd(m_w,3)$ for a component $w_i$ lying above $w$. 

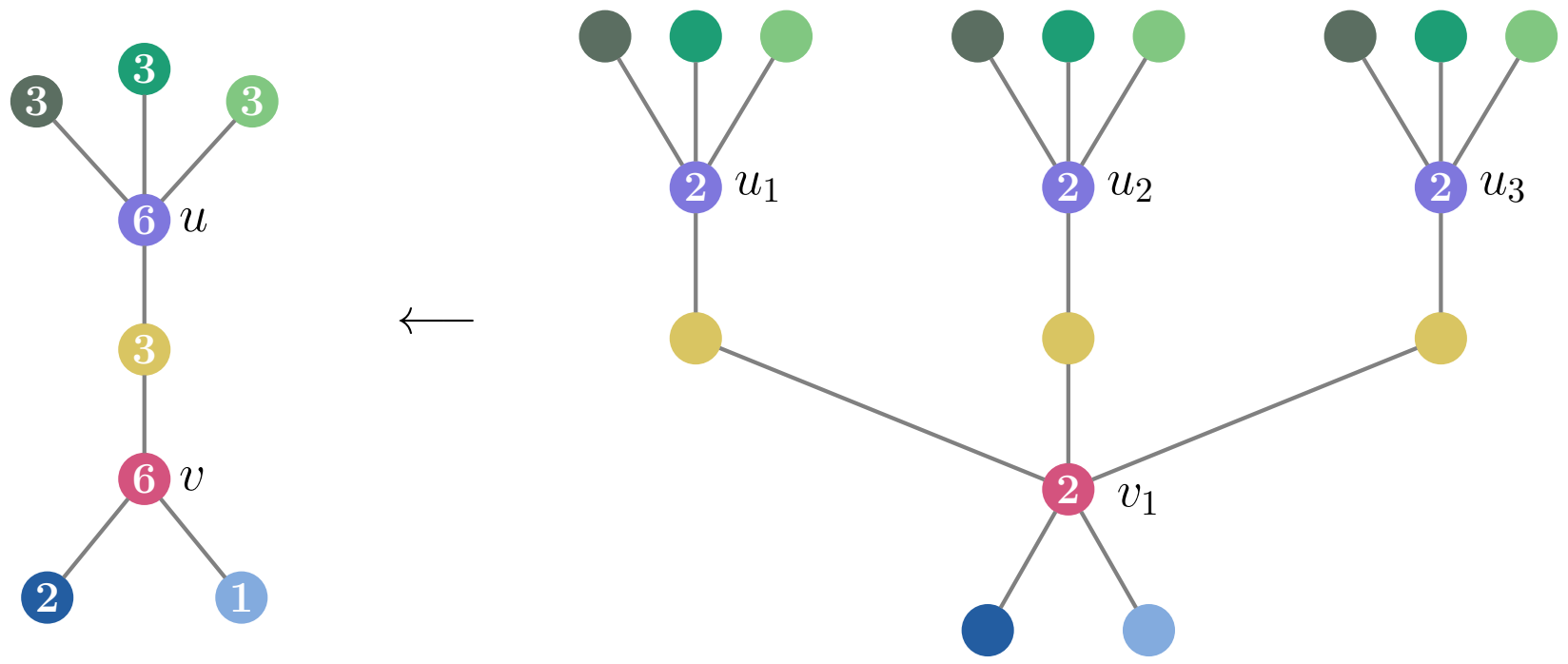
\definecolor{cpurple}{HTML}{7F77DD}
\definecolor{cteal}{HTML}{81C781}
\definecolor{cgreen}{HTML}{1D9E75}
\definecolor{csage}{HTML}{5B6E61}
\definecolor{ccoral}{HTML}{D85A30}
\definecolor{camber}{HTML}{D9C562}
\definecolor{cpink}{HTML}{D4537E}
\definecolor{cblue}{HTML}{235DA1}
\definecolor{cred}{HTML}{E24B4A}
\definecolor{cbabyblue}{HTML}{83ABDE}

\begin{figure}[H]
\centering
\begin{tikzpicture}[
  baseline=(current bounding box.center),
  scale=0.5, transform shape,
  vtx/.style={circle,draw=#1,fill=#1,minimum size=7mm,inner sep=0pt,
              text=white,font=\huge\bfseries},
  lbl/.style={text=black,font=\Huge},
  every edge/.style={draw,gray,thick}
]
  \node[vtx=cteal]   (a) at ( 1.5, 1.65) {3};
  \node[vtx=cgreen]  (b) at ( 0  , 2.1 ) {3};
  \node[vtx=csage]   (c) at (-1.5, 1.65) {3};
  \node[vtx=cpurple, label={[lbl]right:$u$}] (u) at (0,0) {6};
  \node[vtx=camber]  (w) at ( 0  ,-1.8 ) {3};
  \node[vtx=cpink,   label={[lbl]right:$v$}] (v) at (0,-3.6) {6};
  \node[vtx=cbabyblue] (p) at ( 1.35,-5.25) {1};
  \node[vtx=cblue]     (q) at (-1.35,-5.25) {2};

  \draw (u) edge (a) edge (b) edge (c) edge (w);
  \draw (w) edge (v);
  \draw (v) edge (p) edge (q);
\end{tikzpicture}
\qquad$\longleftarrow$\qquad
\begin{tikzpicture}[
  baseline=(current bounding box.center),
  scale=0.7, transform shape,
  vtx/.style={circle,draw=#1,fill=#1,minimum size=5mm,inner sep=0pt,
              text=white,font=\Large\bfseries},
  lbl/.style={text=black,font=\LARGE},
  every edge/.style={draw,gray,thick}
]
  
  \node[vtx=csage]  (a1) at (-4.6,4.5) {};
  \node[vtx=cgreen] (a2) at (-3.7,4.5) {};
  \node[vtx=cteal]  (a3) at (-2.8,4.5) {};
  \node[vtx=csage]  (b1) at (-0.9,4.5) {};
  \node[vtx=cgreen] (b2) at ( 0  ,4.5) {};
  \node[vtx=cteal]  (b3) at ( 0.9,4.5) {};
  \node[vtx=csage]  (c1) at ( 2.8,4.5) {};
  \node[vtx=cgreen] (c2) at ( 3.7,4.5) {};
  \node[vtx=cteal]  (c3) at ( 4.6,4.5) {};

  \node[vtx=cpurple, label={[lbl]right:$u_1$}] (u1) at (-3.7,3) {2};
  \node[vtx=cpurple, label={[lbl]right:$u_2$}] (u2) at ( 0  ,3) {2};
  \node[vtx=cpurple, label={[lbl]right:$u_3$}] (u3) at ( 3.7,3) {2};

  \node[vtx=camber] (m1) at (-3.7,1.5) {};
  \node[vtx=camber] (m2) at ( 0  ,1.5) {};
  \node[vtx=camber] (m3) at ( 3.7,1.5) {};

  \node[vtx=cpink, label={[lbl,xshift=1mm,yshift=-1mm]right:$v_1$}] (v1) at (0,0) {2};
  \node[vtx=cblue]     (q1) at (-0.8,-1.4) {};
  \node[vtx=cbabyblue] (q2) at ( 0.8,-1.4) {};

  \draw (u1) edge (a1) edge (a2) edge (a3) edge (m1);
  \draw (u2) edge (b1) edge (b2) edge (b3) edge (m2);
  \draw (u3) edge (c1) edge (c2) edge (c3) edge (m3);
  \draw (v1) edge (m1) edge (m2) edge (m3) edge (q1) edge (q2);
\end{tikzpicture}
\end{figure}

\noindent Using the genus formula, we obtain the geometric genera $g(u_i)=g(v_1)=0$. By Proposition \ref{Halle-Nicaise 1.3.2}, any other components in $R'$ also have geometric genus 0. To compute the number of components and the multiplicities corresponding to the outer chains, we use Proposition \ref{initial multiplicity of chain after base change}. Then each chain lying above a chain in $\Gamma$ with initial multiplicity 3 (i.e. the green chains) has initial multiplicity $3/3=1$. Each outer chain lying above one of the other two outer chains, which have multiplicities $d=1,2$, respectively, has initial multiplicity $d \text{ mod } 2$. For $d=2$, this implies that the chain lying above it gets contracted. An outer chain with initial multiplicity $1$ has exactly one component, so we are done computing the outer chains. Similarly, the inner chains $Z_i$ between $u_i$ and $v$, which lie above the unique inner chain $Z$ in $R$, have initial multiplicity 1 (both at $u_i$ and at $v_1$). Moreover, Proposition \ref{depth} implies that the depth of $Z_i$ is the same as the depth of $Z$. Therefore, the reduction type $R'$ is the following.

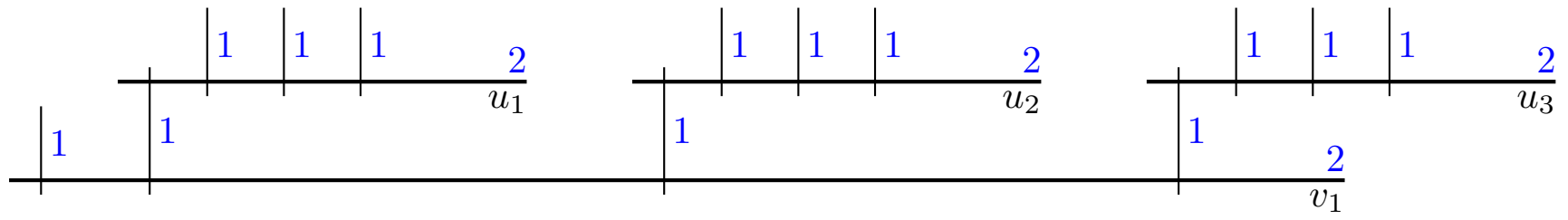
\begin{figure}[H]
\begin{tikzpicture}[xscale=0.7,yscale=0.9]
\draw[thick,shorten >=2] (-1/3,0)--(137/10,0) node[inner sep=2pt,blue,above 
left,scale=0.8] {2} node[inner sep=2pt,below left,scale=0.8] {$v_1$};
\draw[shorten <=-3pt] (0,0)--node[inner sep=2pt,scale=0.8,blue,right] {1} 
(0,3/5);
\draw[shorten <=-3pt, shorten >=-3pt] (17/15,0)--node[inner 
sep=2pt,scale=0.8,blue,right] {1} (17/15,4/5);
\draw[thick,shorten >=2] (4/5,4/5)--(31/6,4/5) node[inner sep=2pt,blue,above 
left,scale=0.8] {2} node[inner sep=2pt,below left,scale=0.8] {$u_1$};
\draw[shorten <=-3pt] (26/15,4/5)--node[inner sep=2pt,scale=0.8,blue,right] {1} 
(26/15,7/5);
\draw[shorten <=-3pt] (38/15,4/5)--node[inner sep=2pt,scale=0.8,blue,right] {1} 
(38/15,7/5);
\draw[shorten <=-3pt] (10/3,4/5)--node[inner sep=2pt,scale=0.8,blue,right] {1} 
(10/3,7/5);
\draw[shorten <=-3pt, shorten >=-3pt] (13/2,0)--node[inner 
sep=2pt,scale=0.8,blue,right] {1} (13/2,4/5);
\draw[thick,shorten >=2] (37/6,4/5)--(158/15,4/5) node[inner sep=2pt,blue,above 
left,scale=0.8] {2} node[inner sep=2pt,below left,scale=0.8] {$u_2$};
\draw[shorten <=-3pt] (71/10,4/5)--node[inner sep=2pt,scale=0.8,blue,right] {1} 
(71/10,7/5);
\draw[shorten <=-3pt] (79/10,4/5)--node[inner sep=2pt,scale=0.8,blue,right] {1} 
(79/10,7/5);
\draw[shorten <=-3pt] (87/10,4/5)--node[inner sep=2pt,scale=0.8,blue,right] {1} 
(87/10,7/5);
\draw[shorten <=-3pt, shorten >=-3pt] (178/15,0)--node[inner 
sep=2pt,scale=0.8,blue,right] {1} (178/15,4/5);
\draw[thick,shorten >=2] (173/15,4/5)--(159/10,4/5) node[inner 
sep=2pt,blue,above left,scale=0.8] {2} node[inner sep=2pt,below left,scale=0.8] 
{$u_3$};
\draw[shorten <=-3pt] (187/15,4/5)--node[inner sep=2pt,scale=0.8,blue,right] {1}
(187/15,7/5);
\draw[shorten <=-3pt] (199/15,4/5)--node[inner sep=2pt,scale=0.8,blue,right] {1}
(199/15,7/5);
\draw[shorten <=-3pt] (211/15,4/5)--node[inner sep=2pt,scale=0.8,blue,right] {1}
(211/15,7/5);
\end{tikzpicture}
\end{figure}

\subsection*{Example 2}
\mbox{}
\vspace{5pt}

The reduction type $R$ below is an example in genus 3 of a reduction type that does not satisfy the conditions in Corollary \ref{unique reduction type} for even $n$, so it has more than one reduction type cover. Let us look at the different possible reduction type covers of $R'$ of degree $n=2$.
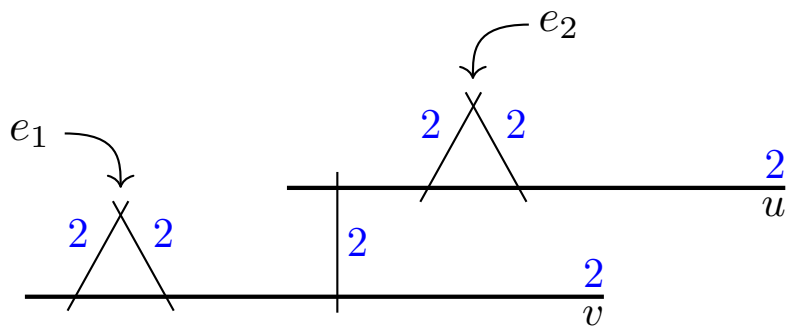
\begin{figure}[H]
\vspace{-10pt}
\begin{tikzpicture}[xscale=1,yscale=0.9]

\draw[->] (-0.07,1.2)
    .. controls (0.3,1.2) and (0.3,1.0) ..
    (0.3,0.8);
\node[scale=0.85] at (-0.3,1.2) {$e_1$};
\draw[->] (3,2.0)
    .. controls (2.63,2.0) and (2.63,1.8) ..
    (2.63,1.6);
\node[scale=0.85] at (3.2,2.0) {$e_2$};

\draw[thick,shorten >=2] (-1/3,0)--(107/30,0) node[inner sep=2pt,blue,above 
left,scale=0.8] {2} node[inner sep=2pt,below left,scale=0.8] {$v$};
\draw[shorten <=-3pt, shorten >=-3pt] (3/5,0)--node[inner 
sep=2pt,scale=0.8,blue,above right] {2} (3/10,3/5);
\draw[shorten <=-3pt, shorten >=-3pt] (0,0)--node[inner 
sep=2pt,scale=0.8,blue,above left] {2} (3/10,3/5);
\draw[shorten <=-3pt, shorten >=-3pt] (26/15,0)--node[inner 
sep=2pt,scale=0.8,blue,right] {2} (26/15,4/5);
\draw[thick,shorten >=2] (7/5,4/5)--(143/30,4/5) node[inner sep=2pt,blue,above 
left,scale=0.8] {2} node[inner sep=2pt,below left,scale=0.8] {$u$};
\draw[shorten <=-3pt, shorten >=-3pt] (44/15,4/5)--node[inner 
sep=2pt,scale=0.8,blue,above right] {2} (79/30,7/5);
\draw[shorten <=-3pt, shorten >=-3pt] (7/3,4/5)--node[inner 
sep=2pt,scale=0.8,blue,above left] {2} (79/30,7/5);
\end{tikzpicture}
\vspace{-5pt}
\end{figure}
\noindent Let $\A=\ZZ/2\ZZ$. Since all the components in $R$ have geometric genus 0, there is a unique $\A$-dilation datum $D$ compatible with $R$ given by $D(w)=\ZZ/({2/\gcd(\weight(w),2)})\ZZ$ for any component $w$ in $R$. Let $\Gamma$ be the dual graph of $R$, and let $e_1$ and $e_2$ be the edges in $\Gamma$ corresponding to the intersections of components highlighted in the picture above. Using Proposition \ref{cohomology of a graph}, we can compute that 
\[
H^1(\Gamma,\A_D)\iso \A/C(e_1)\oplus \A/C(e_2)=\ZZ/2\ZZ\oplus \ZZ/2\ZZ.
\]
Then the harmonic $\A$-cover of $\Gamma$ corresponding to $(0,0)$ would give rise to a disconnected reduction type as below, which contradicts the connectedness of reduction types.
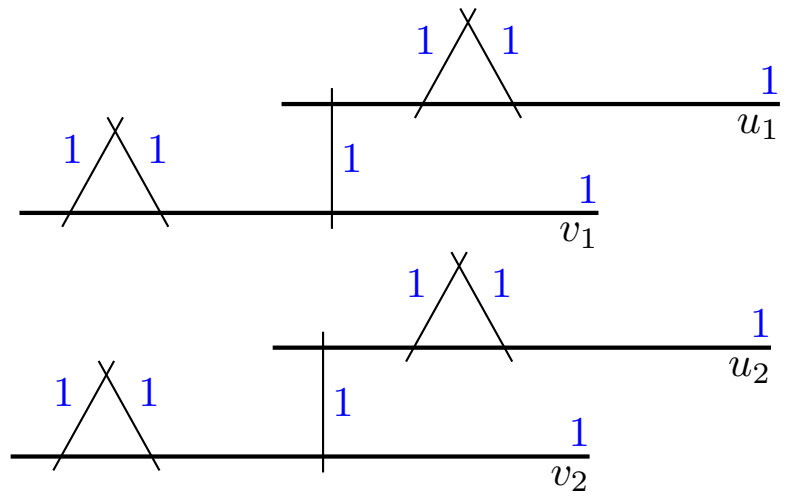
\begin{figure}[htbp]
    \centering
    \begin{subfigure}[b]{0.8\textwidth}
        \centering
    \begin{tikzpicture}[xscale=1,yscale=0.9]
\draw[thick,shorten >=2] (-1/3,0)--(107/30,0) node[inner sep=2pt,blue,above 
left,scale=0.8] {1} node[inner sep=2pt,below left,scale=0.8] {$v_1$};
\draw[shorten <=-3pt, shorten >=-3pt] (3/5,0)--node[inner 
sep=2pt,scale=0.8,blue,above right] {1} (3/10,3/5);
\draw[shorten <=-3pt, shorten >=-3pt] (0,0)--node[inner 
sep=2pt,scale=0.8,blue,above left] {1} (3/10,3/5);
\draw[shorten <=-3pt, shorten >=-3pt] (26/15,0)--node[inner 
sep=2pt,scale=0.8,blue,right] {1} (26/15,4/5);
\draw[thick,shorten >=2] (7/5,4/5)--(143/30,4/5) node[inner sep=2pt,blue,above 
left,scale=0.8] {1} node[inner sep=2pt,below left,scale=0.8] {$u_1$};
\draw[shorten <=-3pt, shorten >=-3pt] (44/15,4/5)--node[inner 
sep=2pt,scale=0.8,blue,above right] {1} (79/30,7/5);
\draw[shorten <=-3pt, shorten >=-3pt] (7/3,4/5)--node[inner 
sep=2pt,scale=0.8,blue,above left] {1} (79/30,7/5);
\end{tikzpicture} 
\end{subfigure}
    \begin{subfigure}[b]{0.8\textwidth}
        \centering
\begin{tikzpicture}[xscale=1,yscale=0.9]
\draw[thick,shorten >=2] (-1/3,0)--(107/30,0) node[inner sep=2pt,blue,above 
left,scale=0.8] {1} node[inner sep=2pt,below left,scale=0.8] {$v_2$};
\draw[shorten <=-3pt, shorten >=-3pt] (3/5,0)--node[inner 
sep=2pt,scale=0.8,blue,above right] {1} (3/10,3/5);
\draw[shorten <=-3pt, shorten >=-3pt] (0,0)--node[inner 
sep=2pt,scale=0.8,blue,above left] {1} (3/10,3/5);
\draw[shorten <=-3pt, shorten >=-3pt] (26/15,0)--node[inner 
sep=2pt,scale=0.8,blue,right] {1} (26/15,4/5);
\draw[thick,shorten >=2] (7/5,4/5)--(143/30,4/5) node[inner sep=2pt,blue,above 
left,scale=0.8] {1} node[inner sep=2pt,below left,scale=0.8] {$u_2$};
\draw[shorten <=-3pt, shorten >=-3pt] (44/15,4/5)--node[inner 
sep=2pt,scale=0.8,blue,above right] {1} (79/30,7/5);
\draw[shorten <=-3pt, shorten >=-3pt] (7/3,4/5)--node[inner 
sep=2pt,scale=0.8,blue,above left] {1} (79/30,7/5);
\end{tikzpicture} 
\end{subfigure}
\vspace{-5pt}
\end{figure}

\noindent Due to the symmetry of $R$ in $u$ and $v$, the elements $(0,1)$ and $(1,0)$ give rise to different reduction type covers, but the same reduction type.

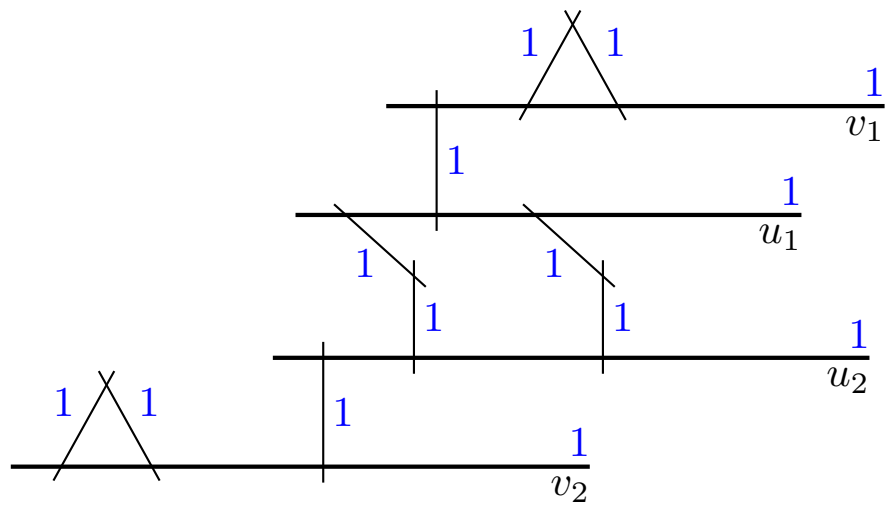
\begin{figure}[H]
\begin{tikzpicture}[xscale=1,yscale=0.9]
\draw[thick,shorten >=2] (-1/3,0)--(107/30,0) node[inner sep=2pt,blue,above 
left,scale=0.8] {1} node[inner sep=2pt,below left,scale=0.8] {$v_2$};
\draw[shorten <=-3pt, shorten >=-3pt] (3/5,0)--node[inner 
sep=2pt,scale=0.8,blue,above right] {1} (3/10,3/5);
\draw[shorten <=-3pt, shorten >=-3pt] (0,0)--node[inner 
sep=2pt,scale=0.8,blue,above left] {1} (3/10,3/5);
\draw[shorten <=-3pt, shorten >=-3pt] (26/15,0)--node[inner 
sep=2pt,scale=0.8,blue,right] {1} (26/15,4/5);
\draw[thick,shorten >=2] (7/5,4/5)--(65/12,4/5) node[inner sep=2pt,blue,above 
left,scale=0.8] {1} node[inner sep=2pt,below left,scale=0.8] {$u_2$};
\draw[shorten <=-3pt, shorten >=-3pt] (7/3,4/5)--node[inner 
sep=2pt,scale=0.8,blue,right] {1} (7/3,7/5);
\draw[shorten <=-3pt, shorten >=-3pt] (7/3,7/5)--node[inner 
sep=2pt,scale=0.8,blue,below left=-1pt and -1pt] {1} (113/60,37/20);
\draw[shorten <=-3pt, shorten >=-3pt] (43/12,4/5)--node[inner 
sep=2pt,scale=0.8,blue,right] {1} (43/12,7/5);
\draw[shorten <=-3pt, shorten >=-3pt] (43/12,7/5)--node[inner 
sep=2pt,scale=0.8,blue,below left=-1pt and -1pt] {1} (47/15,37/20);
\draw[thick,shorten >=2] (31/20,37/20)--(149/30,37/20) node[inner 
sep=2pt,blue,above left,scale=0.8] {1} node[inner sep=2pt,below left,scale=0.8] 
{$u_1$};
\draw[shorten <=-3pt, shorten >=-3pt] (149/60,37/20)--node[inner 
sep=2pt,scale=0.8,blue,right] {1} (149/60,53/20);
\draw[thick,shorten >=2] (43/20,53/20)--(331/60,53/20) node[inner 
sep=2pt,blue,above left,scale=0.8] {1} node[inner sep=2pt,below left,scale=0.8] 
{$v_1$};
\draw[shorten <=-3pt, shorten >=-3pt] (221/60,53/20)--node[inner 
sep=2pt,scale=0.8,blue,above right] {1} (203/60,13/4);
\draw[shorten <=-3pt, shorten >=-3pt] (37/12,53/20)--node[inner 
sep=2pt,scale=0.8,blue,above left] {1} (203/60,13/4);
\end{tikzpicture} 
\end{figure}

\noindent Lastly, the element $(1,1)$ gives rise to the reduction type below.

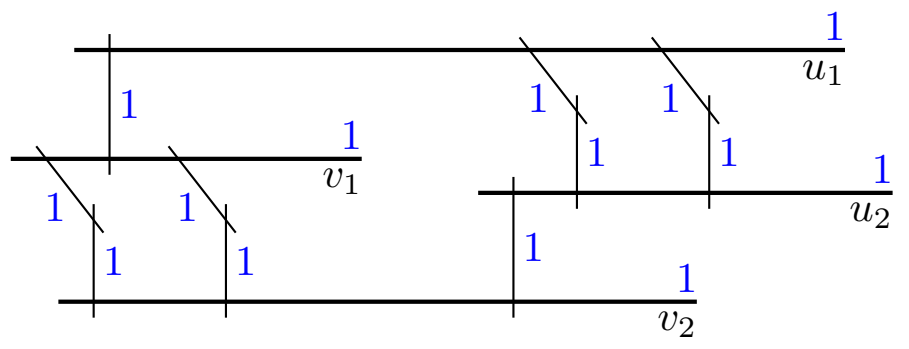
\begin{figure}[H]
\begin{tikzpicture}[xscale=0.7,yscale=0.9]
\draw[thick,shorten >=2] (7/60,0)--(25/4,0) node[inner sep=2pt,blue,above 
left,scale=0.8] {1} node[inner sep=2pt,below left,scale=0.8] {$v_2$};
\draw[shorten <=-3pt, shorten >=-3pt] (9/20,0)--node[inner 
sep=2pt,scale=0.8,blue,right] {1} (9/20,3/5);
\draw[shorten <=-3pt, shorten >=-3pt] (9/20,3/5)--node[inner 
sep=2pt,scale=0.8,blue,below left=-1pt and -1pt] {1} (0,21/20);
\draw[shorten <=-3pt, shorten >=-3pt] (17/10,0)--node[inner 
sep=2pt,scale=0.8,blue,right] {1} (17/10,3/5);
\draw[shorten <=-3pt, shorten >=-3pt] (17/10,3/5)--node[inner 
sep=2pt,scale=0.8,blue,below left=-1pt and -1pt] {1} (5/4,21/20);
\draw[thick,shorten >=2] (-1/3,21/20)--(37/12,21/20) node[inner 
sep=2pt,blue,above left,scale=0.8] {1} node[inner sep=2pt,below left,scale=0.8] 
{$v_1$};
\draw[shorten <=-3pt, shorten >=-3pt] (3/5,21/20)--node[inner 
sep=2pt,scale=0.8,blue,right] {1} (3/5,37/20);
\draw[shorten <=-3pt, shorten >=-3pt] (53/12,0)--node[inner 
sep=2pt,scale=0.8,blue,right] {1} (53/12,4/5);
\draw[thick,shorten >=2] (49/12,4/5)--(81/10,4/5) node[inner sep=2pt,blue,above 
left,scale=0.8] {1} node[inner sep=2pt,below left,scale=0.8] {$u_2$};
\draw[shorten <=-3pt, shorten >=-3pt] (301/60,4/5)--node[inner 
sep=2pt,scale=0.8,blue,right] {1} (301/60,7/5);
\draw[shorten <=-3pt, shorten >=-3pt] (301/60,7/5)--node[inner 
sep=2pt,scale=0.8,blue,below left=-1pt and -1pt] {1} (137/30,37/20);
\draw[shorten <=-3pt, shorten >=-3pt] (94/15,4/5)--node[inner 
sep=2pt,scale=0.8,blue,right] {1} (94/15,7/5);
\draw[shorten <=-3pt, shorten >=-3pt] (94/15,7/5)--node[inner 
sep=2pt,scale=0.8,blue,below left=-1pt and -1pt] {1} (349/60,37/20);
\draw[thick,shorten >=2] (4/15,37/20)--(153/20,37/20) node[inner 
sep=2pt,blue,above left,scale=0.8] {1} node[inner sep=2pt,below left,scale=0.8] 
{$u_1$};
\end{tikzpicture} 
\end{figure}

\subsection*{Example 3}
\mbox{}
\vspace{5pt}

One may think that Theorem \ref{unique reduction type cover} could be phrased only for prime $n$, and compute the reduction type covers of composite degree by going prime-by-prime. However, this does not work: going prime-by-prime might give more possibilities for reduction type covers than there are. We exhibit an example of this phenomenon. Consider the (genus 3) reduction type $R$ below. 
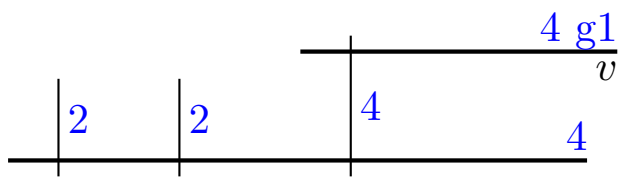
\begin{figure}[H]
\begin{tikzpicture}[xscale=1,yscale=0.9]
\draw[thick,shorten >=2] (-1/3,0)--(107/30,0) node[inner sep=2pt,blue,above 
left,scale=0.8] {4};
\draw[shorten <=-3pt] (0,0)--node[inner sep=2pt,scale=0.8,blue,right] {2} 
(0,3/5);
\draw[shorten <=-3pt] (4/5,0)--node[inner sep=2pt,scale=0.8,blue,right] {2} 
(4/5,3/5);
\draw[shorten <=-3pt, shorten >=-3pt] (29/15,0)--node[inner 
sep=2pt,scale=0.8,blue,right] {4} (29/15,4/5);
\draw[thick,shorten >=2] (8/5,4/5)--(113/30,4/5) node[inner sep=2pt,blue,above 
left,scale=0.8] {4 \smash{g}1} node[inner sep=2pt,below left,scale=0.8] 
{$v$};
\end{tikzpicture}  
\end{figure}
\noindent Then, by Corollary \ref{unique reduction type}, any curve with reduction type $R$ has the same reduction types $R_2$ and $R_4$ after tame base change of degree 2 and 4, respectively. Using the algorithm in Section \ref{algorithm}, we obtain that the reduction type $R_2$ is the following.

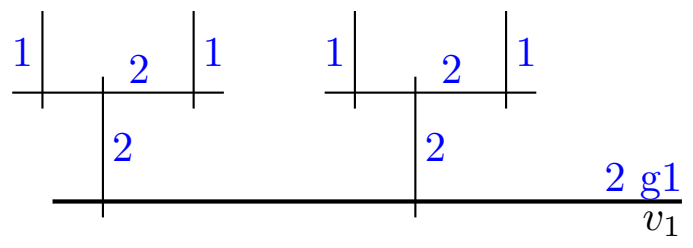
\begin{figure}[H]
\begin{tikzpicture}[xscale=1,yscale=0.9]
\draw[thick,shorten >=2] (4/15,0)--(9/2,0) node[inner sep=2pt,blue,above 
left,scale=0.8] {2 \smash{g}1} node[inner sep=2pt,below left,scale=0.8] 
{$v_1$};
\draw[shorten <=-3pt, shorten >=-3pt] (3/5,0)--node[inner 
sep=2pt,scale=0.8,blue,right] {2} (3/5,4/5);
\draw (0,4/5)--node[inner sep=2pt,scale=0.8,blue,above,pos=0.6] {2} (7/5,4/5);
\draw[shorten <=-3pt] (1/5,4/5)--node[inner sep=2pt,scale=0.8,blue,left] {1} 
(1/5,7/5);
\draw[shorten <=-3pt] (6/5,4/5)--node[inner sep=2pt,scale=0.8,blue,right] {1} 
(6/5,7/5);
\draw[shorten <=-3pt, shorten >=-3pt] (8/3,0)--node[inner 
sep=2pt,scale=0.8,blue,right] {2} (8/3,4/5);
\draw (31/15,4/5)--node[inner sep=2pt,scale=0.8,blue,above,pos=0.6] {2} 
(52/15,4/5);
\draw[shorten <=-3pt] (34/15,4/5)--node[inner sep=2pt,scale=0.8,blue,left] {1} 
(34/15,7/5);
\draw[shorten <=-3pt] (49/15,4/5)--node[inner sep=2pt,scale=0.8,blue,right] {1} 
(49/15,7/5);
\end{tikzpicture} 
\end{figure}
\noindent Using the algorithm again to compute $R_4$ from $R$ we get the following reduction type.

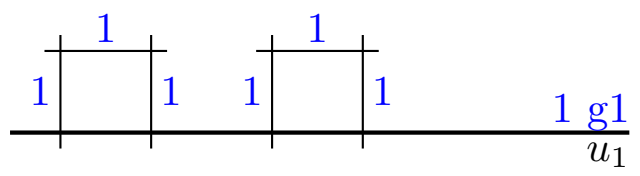
\begin{figure}[H]
\begin{tikzpicture}[xscale=1,yscale=0.9]
\draw[thick,shorten >=2] (-1/3,0)--(23/6,0) node[inner sep=2pt,blue,above 
left,scale=0.8] {1 \smash{g}1} node[inner sep=2pt,below left,scale=0.8] 
{$u_1$};
\draw[shorten <=-3pt, shorten >=-3pt] (0,0)--node[inner 
sep=2pt,scale=0.8,blue,left] {1} (0,3/5);
\draw[shorten <=-3pt, shorten >=-3pt] (0,3/5)--node[inner 
sep=2pt,scale=0.8,blue,above] {1} (3/5,3/5);
\draw[shorten <=-3pt, shorten >=-3pt] (3/5,0)--node[inner 
sep=2pt,scale=0.8,blue,right] {1} (3/5,3/5);
\draw[shorten <=-3pt, shorten >=-3pt] (7/5,0)--node[inner 
sep=2pt,scale=0.8,blue,left] {1} (7/5,3/5);
\draw[shorten <=-3pt, shorten >=-3pt] (7/5,3/5)--node[inner 
sep=2pt,scale=0.8,blue,above] {1} (2,3/5);
\draw[shorten <=-3pt, shorten >=-3pt] (2,0)--node[inner 
sep=2pt,scale=0.8,blue,right] {1} (2,3/5);
\end{tikzpicture} 
\end{figure}
\noindent However, if we input the reduction type $R_2$ and the degree $2$, the algorithm outputs the reduction type $R_4$ and also the reduction type $R_4'$ below.

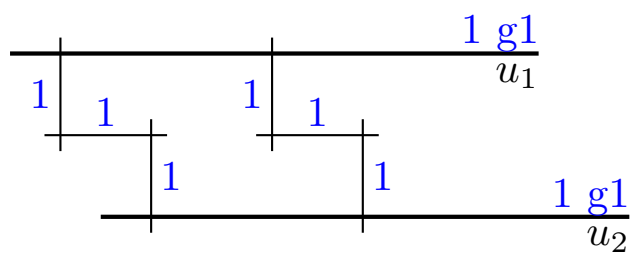
\begin{figure}[H]
\begin{tikzpicture}[xscale=1,yscale=0.9]
\draw[thick,shorten >=2] (4/15,0)--(23/6,0) node[inner sep=2pt,blue,above 
left,scale=0.8] {1 \smash{g}1} node[inner sep=2pt,below left,scale=0.8] 
{$u_2$};
\draw[shorten <=-3pt, shorten >=-3pt] (3/5,0)--node[inner 
sep=2pt,scale=0.8,blue,right] {1} (3/5,3/5);
\draw[shorten <=-3pt, shorten >=-3pt] (3/5,3/5)--node[inner 
sep=2pt,scale=0.8,blue,above] {1} (0,3/5);
\draw[shorten <=-3pt, shorten >=-3pt] (0,3/5)--node[inner 
sep=2pt,scale=0.8,blue,left] {1} (0,6/5);
\draw[shorten <=-3pt, shorten >=-3pt] (2,0)--node[inner 
sep=2pt,scale=0.8,blue,right] {1} (2,3/5);
\draw[shorten <=-3pt, shorten >=-3pt] (2,3/5)--node[inner 
sep=2pt,scale=0.8,blue,above] {1} (7/5,3/5);
\draw[shorten <=-3pt, shorten >=-3pt] (7/5,3/5)--node[inner 
sep=2pt,scale=0.8,blue,left] {1} (7/5,6/5);
\draw[thick,shorten >=2] (-1/3,6/5)--(97/30,6/5) node[inner sep=2pt,blue,above 
left,scale=0.8] {1 \smash{g}1} node[inner sep=2pt,below left,scale=0.8] 
{$u_1$};
\end{tikzpicture} 
\end{figure}

\noindent Therefore, if there is a curve with reduction type $R_2$ that has reduction type $R_4'$ after a tame extension of degree 2, this curve cannot be the base change of a curve with reduction type $R$ by a tame extension of degree 2.

\printbibliography

@book{Halle-Nicaise,
  author = {Halle, Lars Halvard and Nicaise, Johannes},
    title = {Néron Models and Base Change},
    publisher = {Springer Cham},
    year = {2016},
}

@article{Halle,
  title={Stable reduction of curves and tame ramification},
  author={Halle, Lars Halvard},
  journal={Mathematische Zeitschrift},
  volume={265},
  number={3},
  pages={529--550},
  year={2010},
  publisher={Springer}
}

@book{Liu,
  TITLE = {{Algebraic Geometry and Arithmetic Curves}},
  AUTHOR = {Liu, Qing},
  PUBLISHER = {{Oxford University Press,}},
  PAGES = {576},
  YEAR = {2006}
}

@article{Saito,
  title={Vanishing cycles and geometry of curves over a discrete valuation ring},
  author={Saito, Takeshi},
  journal={American Journal of Mathematics},
  volume={109},
  number={6},
  pages={1043--1085},
  year={1987},
  publisher={JSTOR}
}

@misc{Tim2025,
  title={A classification of reduction types of curves},
  author={Dokchitser, Tim},
  eprint={2512.09082},
  archivePrefix={arXiv},
  year={2025}
}

@misc{Jakab,
  title={A classification of reduction types of curves},
  author={Schrettner, Jakab},
  eprint={2607.16159},
  archivePrefix={arXiv},
  year={2026}
}

@article{Liu-Lorenzini-Raynaud,
  title={N{\'e}ron models, Lie algebras, and reduction of curves of genus one.},
  author={Liu, Qing and Lorenzini, Dino and Raynaud, Michel},
  journal={Inventiones mathematicae},
  volume={157},
  number={3},
  year={2004}
}

@article{Lorenzini2000,
author = {Lorenzini, Dino},
year = {2000},
pages = {117-149},
title = {Reduction of points in the group of components of the Néron model of a jacobian},
journal = {Journal für die reine und angewandte Mathematik},
}

@misc{Padma,
    author={Padmavathi Srinivasan},
    title={PhD thesis: Explicit computation of Tamagawa numbers and the sizes of component groups of Néron models of Jacobians using matrix-tree theorems},
    year = {2016},
    url = {https://padmask.github.io/thesis.pdf},
}

@article {ObusWewers,
    AUTHOR = {Obus, Andrew and Wewers, Stefan},
     TITLE = {Explicit resolution of weak wild quotient singularities on
              arithmetic surfaces},
   JOURNAL = {J. Algebraic Geom.},
  FJOURNAL = {Journal of Algebraic Geometry},
    VOLUME = {29},
      YEAR = {2020},
    NUMBER = {4},
     PAGES = {691--728},
}

@article {Tim-Vladimir2013,
author = {Dokchitser, Tim and Dokchitser, Vladimir},
journal = {Acta Arithmetica},
number = {1},
pages = {95-100},
title = {A remark on Tate's algorithm and Kodaira types},
volume = {160},
year = {2013},
}

@article{Len-Ulrisch-Zakharov, 
title={Abelian tropical covers}, 
volume={176}, 
number={2}, 
journal={Mathematical Proceedings of the Cambridge Philosophical Society}, 
author={Len, Yoav and Ulirsch, Martin and Zakharov, Dmitry}, 
year={2024}, 
pages={395–416}
}

@misc{v1, 
title={Abelian tropical covers}, 
version = {1},
year = {2019},
eprint={1906.04215v1},
archivePrefix={arXiv},
author={Len, Yoav and Ulirsch, Martin and Zakharov, Dmitry}, 
}

@article {Lorenzini90,
    AUTHOR = {Lorenzini, Dino},
     TITLE = {Dual graphs of degenerating curves},
   JOURNAL = {Math. Ann.},
  FJOURNAL = {Mathematische Annalen},
    VOLUME = {287},
      YEAR = {1990},
    NUMBER = {1},
     PAGES = {135--150},
}

@article{Namikawa-Ueno,
author = {Ueno, Kenji and Namikawa, Y.},
journal = {Manuscripta mathematica},
pages = {143-186},
title = {The Complete Classification of Fibres in Pencils of Curves of Genus Two},
volume = {9},
year = {1973},
}

@article{Winters,
 author = {Gayn B. Winters},
 journal = {American Journal of Mathematics},
 number = {2},
 pages = {215--228},
 title = {On the Existence of Certain Families of Curves},
 volume = {96},
 year = {1974}
}

@article{Gambheera-Vallieres, 
title={Iwasawa theory for branched $\mathbb{Z}_{p}$-towers of finite graphs and Ihara zeta and L-functions}, 
journal={Mathematical Proceedings of the Cambridge Philosophical Society}, 
author={Gambheera, Rusiru and Vallières, Daniel}, year={2026}, 
pages={1–53}
}

@article{Kataoka,
title = {Kida's formula for graphs with ramifications},
journal = {Acta Arithmetica},
author = {Kataoka, Takenori},
year = {2025},
pages = {271-291},
volume = {221},
number = {2},
}

@misc{Mizuno, 
title={Brauer-Kuroda relations for ramified graph covers}, 
version = {1},
year = {2026},
eprint={2609.25111},
archivePrefix={arXiv},
author={Mizuno, Kosuke}, 
}

@unpublished{monodromy,
  author       = {Jakovác, Gergely and Martínez-Marín, Júlia and Perez, Allan and Shi, Jia and Vilà-Casadevall, Jordi},
  title        = {On the Semistability Defect of Hyperelliptic Curves via Cluster Pictures},
  date         = {2026},
  pubstate     = {inpreparation},
}
\end{document}